\documentclass[reqno]{amsart}

\usepackage{amsmath,amsthm,amssymb,amsfonts,mathtools,mathrsfs,enumerate,amscd}
\usepackage[usenames,dvipsnames]{color}
\usepackage{graphicx,psfrag,epsfig,ifpdf}
\usepackage{mathabx}
\usepackage{pstricks}
\usepackage{pst-bezier}
\usepackage{comment}
\usepackage{epstopdf}
\usepackage{hyperref}

\newcounter{parag}[subsection]

\newcounter{paraga}[subsection]
\renewcommand{\theparaga}{{\bf\arabic{paraga}.}}
\newcommand{\paraga}{\medskip \addtocounter{paraga}{1}
\noindent{\theparaga\ } }

\newcounter{pparag}

\newtheorem{thm}{Theorem}

\newtheorem{lemma}{Lemma}
\newtheorem{prop}{Proposition}

\newtheorem{cor}{Corollary}
\newtheorem{Def}{Definition}
\newtheorem{rem}{Remark}

\def\text#1{\,\hbox{#1}\;}

\def\al{\alpha}
\def\be{\beta}
\def\ga{\gamma}
\def\Ga{{\Gamma}}
\def\de{\delta}

\def\De{\Delta}
\def\eps{{\varepsilon}}
\def\ka{\kappa}
\def\la{\lambda}

\def\om{\omega}
\def\Om{\Omega}

\def\sig{{\sigma}}
\def\Sig{{\Sigma}}

\def\th{{\theta}}
\def\Th{\Theta}

\def\ph{\varphi}
\def\ze{{\zeta}}

\def\jA{{\mathscr A}}
\def\jB{{\mathscr B}}
\def\jC{{\mathscr C}}

\def\jN{{\mathscr N}}
\def\jO{{\mathscr O}}
\def\jP{{\mathscr P}}

\def\jT{{\mathscr T}}

\def\jZ{{\mathscr Z}}

\def\A{{\mathbb A}}
\def\B{{\mathbb B}}

\def\N{{\mathbb N}}
\def\R{{\mathbb R}}
\def\T{{\mathbb T}}
\def\Z{{\mathbb Z}}

\def\bA{{\bf A}}
\def\bB{{\bf B}}

\def\Im{{\rm Im\,}}

\def\Sup{\mathop{\rm Sup\,}\limits}
\def\Inf{\mathop{\rm Inf\,}\limits}
\def\Max{\mathop{\rm Max\,}\limits}

\def\dist{{\rm dist\,}}

\def\Fr{{\rm Fr}}

\def\Int{{\rm Int\,}}
\def\diam{{\rm diam\,}}

\def\setm{\setminus}

\def\ov{\overline}
\def\til{\widetilde}
\def\ha{\widehat}
\def\d{\partial}
\def\inv{^{-1}}

\def\pdemi{{\tfrac{1}{2}}}

\def\norm#1{\Vert#1\Vert}
\def\setm{\setminus}

\def\cR{{\mathcal R}}

\def\cU{{\mathcal U}}
\def\cR{{\mathcal R}}

\def\bu{{\bullet}}
\def\e{{\bf e}}

\def\bA{{\bf A}}
\def\Homt{{\bf Homt}}
\def\Hett{{\bf Hett}}
\def\Domt{{\rm{Domt}}}

\def\beq{\begin{equation}}
\def\eeq{\end{equation}}

\def\bu{{\bullet}}

\def\trans{\pitchfork}

\def\dma{\d_\bu\bA}
\def\CC{{\rm cc}}

\def\vide{\varnothing}

\def\Graph{{\rm Graph\,}}

\def\AA{{\bf A}}

\def\cl{{\rm cl\,}}
\def\Fr{{\rm Fr\,}}
\def\Adm{{\rm Coh}}

\def\dma{\Ga(a)}

\def\bmu{{\boldsymbol \mu}}
\def\bga{{\boldsymbol \ga}}

\def\trans{\pitchfork}

\def\Dom{{\rm Dom\,}}

\def\Ess{{\rm Ess}}

\def\Supp{{\rm Supp}}

\def\bU{{\mathcal U}}
\def\Gr{{\rm Graph\,}}
\def\Im{{\rm Im\,}}
\def\g{{\bf g}}
\def\hc{{\,\ha{\!\jC}}}
\def\inte{{\rm int\,}}
\def\bW{{\bf W}}
\def\cR{{\mathcal R}}
\def\Rot{{\rm Rot}}
\def\bI{{\bf I}}
\def\bph{{\boldsymbol \ph}}
\def\bpsi{{\boldsymbol \psi}}
\def\bPhi{{\boldsymbol \Phi}}
\def\bxi{{\boldsymbol \xi}}
\def\bGa{{\boldsymbol \Ga}}
\def\bsig{{\boldsymbol \sig}}
\def\bI{{\bf I}}
\def\bx{{\bf x}}
\def\bz{{\bf z}}

\newcommand{\twist}{{tilt} }

\begin{document}

\title{Diffusing orbits along chains of cylinders}

\author{Marian Gidea}
\address{Department of Mathematical Sciences, 215 Lexington Ave, Yeshiva University, New York, NY 10016, USA}
\email{Marian.Gidea@yu.edu}
\author{Jean-Pierre Marco}
\address{Sorbonne Universit\'e (Paris 6) and  Institut de Math\'ematiques de Jussieu-Paris Rive Gauche,
4 Place Jussieu, 75005 Paris cedex 05, France}
\email{jean-pierre.marco@imj-prg.fr}

\maketitle

\begin{abstract}
We develop a geometric mechanism  to prove the existence of
orbits that drift along a prescribed sequence of  cylinders, under some general conditions on the dynamics.
This mechanism can be used to prove the existence of Arnold diffusion for large families of perturbations of Tonelli Hamiltonians on $\A^3$. Our approach can also be applied to more general Hamiltonians that are not necessarily convex.

The main geometric objects in our framework are
$3$--dimensional invariant cylinders with boundary
(not necessarily hyperbolic), which are assumed to  admit center-stable and center-unstable
manifolds. These enable us to define \emph{chains of cylinders},  i.e., finite, ordered families
of cylinders where each cylinder admits homoclinic connections,
and any two consecutive cylinders in the chain admit heteroclinic connections.

Our main result is on the existence of diffusing orbits which drift along such chains of cylinders,
under precise conditions on the dynamics on the cylinders --  i.e., the existence of  Poincar\'e sections
with the return maps satisfying a {\em tilt condition} -- and on
the geometric properties of the intersections of the center-stable and center-unstable
manifolds of the cylinders  -- i.e., certain compatibility conditions between the  \emph{tilt map} and the
homoclinic maps attached to its essential invariant circles.

We give two proofs of our result, a very short and abstract one, and a more constructive one,
aimed at possible applications to concrete systems.
\end{abstract}

\vskip2cm

%%%%%%%%%%%%%%%%%%%%%%%%%%%%%%%%%%%%%%%%%%%%%%
%%%%%%%%%%%%%%%%%%%%%%%%%%%%%%%%%%%%%%%%%%%%%%
%%%%%%%%%%%%%%%%%%%%%%%%%%%%%%%%%%%%%%%%%%%%%%

\section{Introduction and main result}

The present work  originates in a geometric proof of the Arnold conjecture
in the Mather {\em a priori stable} setting\footnote{We  pay homage to the fundamental role played by  Mather in the development of the field, and
in providing many ideas and inspiration \cite{Ma04,Ma10,Ma12}.} for 3-degree-of-freedom Hamiltonian systems (see \cite{Mar1,Mar2} which give the necessary
geometric framework and the application to the proof).
Our paper is however self-contained, and our results can also be applied to other contexts, as it will be detailed below.

\paraga  Let us first informally recall the convex setting  for Arnold diffusion.
Given $n\geq1$, we denote by $\A^n=\T^n\times\R^n$ the cotangent bundle of the torus $\T^n$, equipped
with its natural angle-action coordinates $(\th,r)$ and its exact-symplectic form $\Om=\sum_{i=1}^ndr_i\wedge d\th_i$.
Consider a Hamiltonian of class $C^\ka$ on $\A^3$, of the form
\beq\label{eq:hampert}
H(\th,r)=h(r)+ f(\th,r), \qquad (\th,r)\in\A^3,
\eeq
where $\ka$ is large enough, and $h$ is strictly convex and with superlinear growth at infinity w.r.t. the actions\footnote{The strict convexity assumption is not present in Arnold's original formulation of the problem, while the superlinear growth is assumed only to get compact energy levels}.
We think of the Hamiltonian $H$ as a small perturbation of the Hamiltonian $h$. This setting of the  Arnold diffusion problem, where the dynamics of $h$ is described in terms of angle-action variables only, is referred to as the {\em a priori stable} case, whereas the setting where the dynamics of $h$ is  the product of an angle-action factor with
a hyperbolic one is referred to as the  {\em a priori unstable} case.

The main problem is the {\em a priori} stable case. Namely, fix a regular level set $h\inv(\e)$ of the unperturbed energy in the action space $\R^3$,
consider an arbitrary family $(\ha O_i)_{1\leq i\leq m}$ of small open subsets in
$\R^3$ which intersect $h\inv(\e)$, and set $O_i=\T^3\times \ha O_i$.
The diffusion problem  is on proving that the system $H$ admits orbits visiting each
$O_i$, for a large class of perturbations $f$
in $C^\ka(\A^3)$. Such orbits will be referred to as  {\em diffusing orbits}.  The {\em a priori} unstable case is an intermediate but very interesting problem
which serves as a guide for developing the methods.

Three independent approaches of the diffusion problem in the {\em a priori} stable setting were developed recently:
see \cite{C19,KZ,Mar2} and references therein.
A common feature of the works \cite{C19,KZ,Mar2} is the use of ``normally hyperbolic cylinders'' which form
``chains'' intersecting the open sets  $O_i$. Once the existence of such chains of cylinders is shown
(under appropriate nondegeneracy conditions on $f$), proving the existence
of diffusing orbits amounts to proving the existence of orbits ``drifting along the cylinders as well as from one cylinder to the next'',
possibly under additional nondegeneracy conditions on $f$.

The very definition of normally hyperbolic  cylinders and chains  is {\em not} the same in the three approaches,
mainly regarding the invariance condition of the cylinders under the Hamiltonian flow.
In \cite{Mar2} the cylinders are $3$-dimensional,
{\em genuinely invariant}, normally hyperbolic and compact (with boundary), which ensures the existence of  well-defined asymptotic manifolds,
namely the center-stable and center-unstable manifolds consisting of orbits that approach the cylinder in forward and backward time, respectively.

\paraga In this paper we consider a more general situation, that of {\em tame cylinders}, which possess asymptotic manifolds and admit a
lambda-lemma type property, see below.
Crossings between these asymptotic manifolds enables us to define  homoclinic connections for each cylinder, as well as heteroclinic connections
between consecutive cylinders, which allows us to formulate a natural geometric definition of a {\em chain of cylinders}. The main interest in this notion is that non-convex
integrable systems yield invariant manifolds which are, as a rule, not normally hyperbolic but only tame.

Once a chain of tame cylinders is given, the central idea of our approach is to perform a systematic reduction to two dimensional dynamics,
which allows us to go back to very
simple arguments to detect the existence of diffusing orbits. Here again, we have non-convex situations in mind and we assume that the dynamics
inside our cylinders admit a section whose return map is a {\em tilt map} (instead of the more usual assumption of being a {\em twist map}).
Then, we will not need more than the usual Birkhoff theory of twist maps\footnote{
with its now well-known extensions to tilt maps},  properly generalized here following an initial idea by Moeckel (see \cite{M02}), to produce
diffusing pseudo-orbits along the chain (see the definition below). Thanks to the lambda-lemma property
and the Poincar\'e recurrence theorem, a suitable shadowing process then easily yields the existence of true orbits of the perturbed system
following those pseudo-orbits, and therefore drifting along the chain.

\paraga We now briefly describe our set-up and the main results.
We consider a $C^\ka$  Hamiltonian function $H$ on $\A^3$, with $\ka\geq 2$.
Although it is not absolutely necessary, we find it convenient to restrict our study to regular levels
$H\inv(\e)$ on which the Hamiltonian vector field is complete. In this case we say that
$\e$ (or the level $H\inv(\e)$) is {\em completely regular}.

In this paper, a {\em cylinder}  is a compact submanifold with boundary of $\A^3$, diffeomorphic to
$\T^2\times[0,1]$. We consider cylinders invariant under the flow
generated by a Hamiltonian function $H$ on $\A^3$, and contained in a given, completely regular energy level $H\inv(\e)$ of~$H$.
%For such invariant cylinders, the notion of normal hyperbolicity {\em relative  to the energy level} can be
%well-defined.
In this paper we consider  more general conditions on the cylinders
rather than the usual normal hyperbolicity, and call them {\em tame cylinders};  details
will be given in Section~\ref{Sec:cylinders}.  In particular, a tame cylinder $\jC$ is $3$-dimensional,
admits $4$-dimensional center-stable and center-unstable manifolds (with boundary) $W^\pm(\jC)$,
which are themselves foliated by the $1$-dimensional stable
and  unstable manifolds $W^\pm(x)$, respectively, associated to the points $x\in\jC$.
Moreover, a tame cylinder admits a lambda-lemma property\footnote
{A normally hyperbolic invariant cylinder is a particular case of a tame cylinder.}.

A {\em chain of cylinders} for $H$ is a finite ordered sequence $(\jC_k)_{1\leq k\leq k_*}$ of tame cylinders contained
in the same energy level $H\inv(\e)$, such that
 each cylinder $\jC_k$, for $1\leq k\leq k_*$, admits homoclinic connections, that is
\[
W^-(\jC_k)\cap W^+(\jC_{k })\neq \emptyset,
\]
and each consecutive pair  of cylinders  $\jC_k$  and  $\jC_{k+1}$ in the chain, for
$1\leq k\leq k_*-1$,
admits heteroclinic connections,  that is
$$
W^-(\jC_k)\cap W^+(\jC_{k+1})\neq \emptyset.
$$

\begin{figure}
\centering
\includegraphics[width=0.9\textwidth]{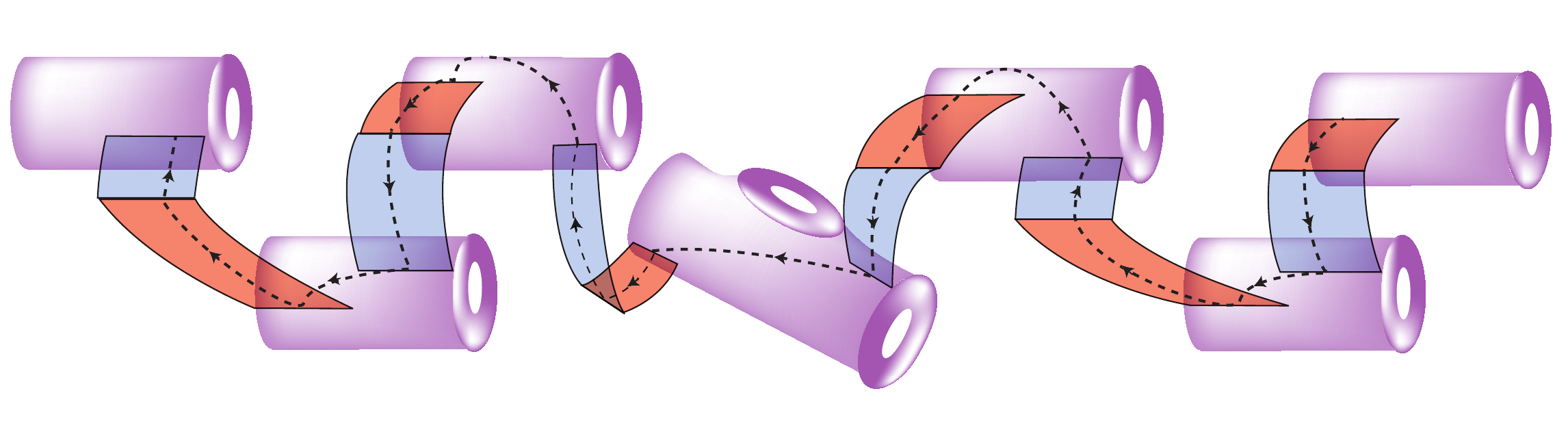}
\caption{An illustration of a cylinder chain, consisting of a sequence of  tame cylinders  with boundary,
connected via homoclinic and heteroclinic orbits. This chain also contains a  singular cylinder (shown in the center of the figure);  see Section~\ref{sec:singular} for details.}
\label{fig:cylinder_chain}
\end{figure}

To ensure the existence of orbits  drifting along a chain,  we will require additional conditions on the dynamics on the
cylinders, and their homoclinic and heteroclinic connections. In this paper, a cylinder (resp. a chain) satisfying those additional properties is called a {\em good cylinder} (resp. a {\em good chain}, see Figure \ref{fig:cylinder_chain});
see Section~\ref{Sec:goodcyl} for the definitions.

Given a good cylinder $\jC$,
we define an {\em essential subtorus} of  $\jC$ as a $2$--dimensional torus
(not necessarily differentiable) contained in $\jC$ and  invariant under the Hamiltonian flow,
which intersects a certain Poincar\'e section $\Sig\sim\T\times[a,b]$ along an essential circle\footnote{that is,
homotopic to $\T\times\{a\}$ inside $\Sigma$}
(see Definition~\ref{def:twistsec}).
Examples of essential subtori are the components of the boundary $\d \jC$;
any essential subtorus is homotopic in $\jC$ to each of these components. The family of essential
subtori will serve us as a guide to building our drifting orbits.

\begin{Def}\label{def:admissible}
Consider a good chain $(\jC_k)_{1\leq k\leq k_*}$ contained in some completely regular energy level $H\inv(\e)$.
Given $\de>0$,  we say that an orbit of the Hamiltonian flow is {\em $\de$-admissible for the
chain} when it intersects the $\de$-neighborhood in $H\inv(\e)$ of any essential subtorus of
$\jC_k$, $1\leq k\leq k_*$.
\end{Def}

The notion of $\delta$-admissible orbits reveals itself to be useful in \cite{Mar2} to obtain
the existence of diffusing orbits for systems (\ref{eq:hampert}) once  a chain of cylinders
that intersects  each open set $O_i$ is constructed\footnote{Under suitable assumptions, each open set $O_i$
contains an essential torus}.
A related notion is that of a $\de$-good chain, which involves a quantitative control on the homoclinic and heteroclinic connections
of the cylinders. This will be introduced in Definitions~\ref{def:bounded} and~\ref{def:goodchains}.
%Indeed, for small enough perturbations $f$, it can be proved
%that a cylinder intersecting $O_i$ necessarily contains  an essential
%torus which is itself contained in~$O_i$. So any $\de$-admissible orbit for  the chain
%intersects each set $O_i$ too, provided  that $\de$ is small enough.

\vskip2mm

The first main result of this paper is the following.

\begin{thm}\label{thm:main1} Let $H$ be a $C^\ka$ Hamiltonian on $\A^3$, with $\ka\geq 2$,
 and let $\e$ be a completely
regular value of~$H$.
Fix $\de>0$.
Then, for any $\de$-good chain of cylinders contained in $H\inv(\e)$,
there exists a $\de$-admissible orbit for the chain.
\end{thm}

\vskip2mm

{\em Otherwise explicitly mentioned, the energies $\e$ and the energy levels we consider in the following will be assumed to be completely regular.}

\vskip2mm

The proof of {\bf Theorem~\ref{thm:main1}} is divided into two parts. In the first part, we introduce
a {\em polysystem} of maps and homoclinic or heteroclinic correspondences associated with the chain of cylinders, and prove the existence
of ``drifting pseudo-orbits'' for this polysystem. This is the content of {\bf Theorem~\ref{thm:pseudoorb1}}, which is
stated and proved in Section~\ref{Sec:pseudoorb} (with an algorithmic
version in
Section~\ref{Sec:constructive}).  In the second part we derive a  shadowing
process (similar to but simpler than that in \cite{GLS}) to prove the existence of genuine
orbits of the Hamiltonian system which intersect arbitrarily small neighborhoods of each point
of the pseudo-orbits of {\bf Theorem~\ref{thm:main1}}; see also the related shadowing results in
\cite{BT99,DLS00,DGR13,GR13,GT17}.
This is the content of {\bf Theorem~\ref{thm:shadowing}}, stated
and proved in Section~\ref{Sec:shadowing}.
We now describe informally
both theorems.

\paraga  %We begin with {\bf Theorem~\ref{thm:pseudoorb1}}, and
We first describe the main objects involved in the construction of drifting orbits and pseudo-orbits.
Here we call a polysystem on some space $X$ a dynamical system formed by a finite family of maps $(f_i)_{i\in I}$, each defined on some subset of $X$ (possibly the whole set $X$),  which can be iterated in any order for which the composition of maps is well defined.

Consider a good chain $(\jC_k)_{1\leq k\leq k_*}$ contained in a level $H\inv(\e)$.
Via the shadowing process developed in Section~\ref{Sec:shadowing},
we will reduce the dynamics in the neighborhood of the cylinders $\jC_k$
and along their homoclinic and heteroclinic connections to that of a  polysystem on
$2$-dimensional annuli $(\Sig_k)_{1\leq k\leq k_*}$, which are diffeomorphic to $\T\times[0,1]$. Let us describe this polysystem.

In the neighborhood of a single cylinder $\jC:=\jC_k$ of the chain,
we will consider a polysystem $(f_i)_{i\in I}$
defined on a global section $\Sig$
of the Hamiltonian flow restricted to $\jC$.
The polysystem consists of a pair $(\ph,\psi)$,
where $\ph$ is a diffeomorphism on $\jC$,  and $\psi$ is a homoclinic {\em correspondence}\footnote{In Section~\ref{Sec:goodcyl} we will consider homoclinic correspondences consisting of families of local diffeomorphisms on $\jC$, whose domains are not necessarily mutually disjoint} on $\jC$,
whose existence relies on additional conditions imposed on good cylinders.
By construction, the global section $\Sig$ is endowed with symplectic coordinates $(\th,r)\in\T\times[a,b]$
for some $0<a<b$.
The diffeomorphism $\ph$ is the Poincar\'e return map associated to $\Sig$, so $\ph$ is symplectic,
We assume $\ph$ to tilt all vertical lines $\th=\textrm{cst.}$  to the right (or to the left).
Other mild conditions on $\ph$ will also be required (see Definition~\ref{def:good}).
The essential subtori of $\jC$ as those which intersect $\Sig$ along essential
invariant circles of $\ph$.

The homoclinic correspondence $\psi$ is determined by the homoclinic orbits to $\jC$.
A first natural homoclinic correspondence is defined on $\jC$ rather than on $\Sig$,
and associates to each element $x\in\jC$ the subset of
all $y\in\jC$ such that  $W^-(x)$ intersects $W^+(y)$ (with additional transversality requirements)\footnote{This type of correspondence was previously studied in \cite{DLS00,DLS06a}, where the authors use an appropriate restriction of the correspondence, which they refer to as the {\em scattering map}}.

To derive the  homoclinic correspondence $\psi$ on $\Sig$  from the  homoclinic correspondence  on $\jC$,
we perform a ``reduction''
to the section $\Sig$, which is achieved by a suitable transport by the Hamiltonian flow inside~$\jC$.
The correspondence $\psi$ turns out to be a collection of local diffeomorphisms of $\Sig$, whose
domains may intersect one another.

This way, each cylinder $\jC_k$ of the chain is equipped with a polysystem $(\ph_k,\psi_k)$. Finally, the
complete polysystem attached to the chain $(\jC_k)_{1\leq k\leq k_*}$ is formed by the collection of the polysystems attached to each cylinder, together with additional  local diffeomorphisms corresponding to
the heteroclinic maps or to the transitions maps (defined along the flow)  between consecutive cylinders.

\paraga We  now describe briefly {\bf Theorem~\ref{thm:pseudoorb1}}, which is given in Section~\ref{Sec:pseudoorb}.
This theorem  asserts the existence, for every $\de>0$,  of a {\em $\de$-admissible pseudo-orbit}
for the polysystems $(\ph_k,\psi_k)$, that is, an orbit of the polysystem which passes $\de$-close to every  essential invariant circle
contained in each of the sections $\Sig_k\subset\jC_k$. In particular, such pseudo-orbits drift along the chain in a natural sense.

Let us restrict ourselves first to a single cylinder $\jC:=\jC_k$ of the chain, equipped with a section  $\Sig\subset\jC$
and a polysystem $(\ph,\psi)$. To prove the existence of $\de$-admissible pseudo-orbits,
 we will require additional compatibility conditions between
$\ph$ and $\psi$, as explained below.

Our result generalizes the approach introduced by Moeckel in
\cite{M02} and developed by Le Calvez in \cite{LC07}, see also \cite{NP12} and references therein.
Their main result is that ``generically'', any polysystem formed
by an area-preserving  \twist  map\footnote{We recall here that a tilt map is a map that tilts every vertical line always to the right (or always to the left);
for instance any twist map, or any finite composition of twist maps (all being right-twists or left-twists), is a tilt map.}
 $\ph$ on $\bA=\T\times[a,b]$ and a globally defined area-preserving diffeomorphism
$\psi$ of $\bA$ admits a  finite ``connecting pseudo-orbit'' whose first point is arbitrarily
close to $\T\times\{a\}$ and whose last point is arbitrarily close to $\T\times\{b\}$.
The single map $\ph$ would not in general admit a connecting orbit, due to the existence of essential invariant circles that separate the  annular section.
The role of the diffeomorphism $\psi$ is to allow the pseudo-orbits of the polysystem to ``jump''  over  the essential
invariant circles of $\ph$ -- under the crucial assumption that such an invariant circle is {\em not} invariant under
 $\psi$.

In our case,
the main difficulty is that $\psi$ is {\em not} everywhere defined (and, in general, is multivalued). The
compatibility conditions that we require generalize in some natural sense  the previous
``no simultaneous invariant circle condition.''
We now require that the range of $\psi$ intersects every
essential invariant circle $\Ga$ of $\ph$ and that $\psi$ ``breaks'' this circle, in the sense that
$\psi^{-1}(\Ga)$
  admits a ``topologically transverse'' intersection{\footnote{Specifically, we require that
 $\psi^{-1}(\Ga)$ contains  certain arcs emanating from  $\Ga$ which are below $\Ga$;
we can think of $\psi^{-1}(\Ga)$ as a union of arcs,
   since $\psi$ is a collection of local
diffeomorphisms} with $\Ga$ (see Definition \ref{def:splittingarc}).
Under this mild condition, we prove that the polysystem $(\ph,\psi)$ admits connecting pseudo-orbits
similar to those in \cite{M02} and \cite{LC07}. A further restriction on the size of the domain of the correspondence
$\psi$  (the $\de$-bounded property, see Definition \ref{def:bounded}) guarantees that connecting pseudo-orbits are $\de$-admissible.

Once the existence of $\de$-admissible pseudo-orbits along a single cylinder is proved,
the existence of $\de$-admissible pseudo-orbits for the full polysystem along the whole chain of cylinders
follows immediately, once we require some additional quantitative conditions on the heteroclinic maps  between consecutive sections.

The proof of {\bf Theorem~\ref{thm:pseudoorb1}} is given in Section \ref{sec:pseudoorbchain}.
It relies on {\bf Proposition \ref{prop:tamecyl}}.
We will give two different proofs of {\bf Proposition \ref{prop:tamecyl}}.

The first  proof of {\bf Proposition ~\ref{prop:tamecyl}}, given in Section \ref{sec:pseudoorbcyl}, is non-constructive
and extends quite directly the methods introduced in \cite{M02} to the case of polysystems of correspondences.
It requires very mild nondegeneracy conditions, which facilitate  the applications to
occurrence of diffusion under ``generic'' conditions on the perturbation.

The second proof of {\bf Proposition ~\ref{prop:tamecyl}} relies on  {\bf Proposition~\ref{prop:tamecyl2}}, whose statement and proof are given in Section~\ref{Sec:constructive}. The assumptions on the polysystem in {\bf Proposition~\ref{prop:tamecyl2}} are  slightly more stringent, but enable us to use an iterative ``Birkhoff procedure.''
Due to the  ``algorithmic'' nature of the underlying construction,
we expect this approach to be applicable to specific examples.

\paraga We now  briefly describe   {\bf Theorem~\ref{thm:shadowing}} from Section \ref{Sec:shadowing}, on the shadowing of pseudo-orbits of a polysystem.

We first consider the case of a
single good cylinder. Given
a pseudo-orbit $(x_n)_{1\leq n\leq n_*}$ of the polysystem $(\ph,\psi)$ on $\Sig$, where $\psi=(\psi_i)_{i\in I}$,
with $\psi_i$ being a local diffeomorphism of $\Sig$ for $i\in I$.  In particular, $\psi$ can consist of  a single global diffeomorphism.
By the definition of a pseudo-orbit, either
$x_{n+1}=\ph(x_n)$ or there exists $i_n\in I$ such that $x_{n+1}=\psi_{i_n}(x_n)$,  for $1\leq n\leq n_*-1$.

In the first case, since $\ph$ is a flow-induced return map,
there is a time $T_n\geq 0$ such that $x_{n+1}=\Phi_H^{T_n}(x_n)$ (where $\Phi_H^T$ stands for the time-$T$ map
of the Hamiltonian flow).

In the second case, using the definition of $\psi$ as a reduced homoclinic correspondence,
 together with the Poincar\'e recurrence theorem applied to $\ph$, it turns out that there exist
a point $\xi\in H\inv(\e)$ and a time $T_n>0$ such that  $\xi$ is arbitrarily close to $x_n$ and
$\Phi_H^{T_n}(\xi)$ is arbitrarily close to $x_{n+1}$. The possibility of using the recurrence properties
of $\ph$ follows directly from the symplectic nature of our setting and the compactness
of $\Sig$ (hence, its finite measure),  and therefore is intimately related to our definition of cylinders.
The idea of using the recurrence property of the dynamics in such a context was  used in \cite{GLS} (see also \cite{NP12}).

One can expect that slightly perturbing and ``gluing together'' the previous pieces of Hamiltonian orbits
$\Phi_H\big([0,T_n]\times\{x_n\}\big)$ yield a genuine orbit of the Hamiltonian flow which contains points located arbitrarily close
to each point $x_n$ of the initial orbit.
This is the main statement of {\bf Theorem~\ref{thm:shadowing}}.
The proof uses the
inclination property (a $\lambda$-lemma for arcs in the setting of tame cylinders) to construct a positively invariant sequence
of balls~$B_n$ centered on the unstable manifolds of the points $x_n$ (and arbitrarily close to them),
such that $\Phi_H^{T_n}(B_n)\subset B_{n+1}$ for a suitable $T_n>0$. This proves our claim;
for similar approaches, see also \cite{BT99,DLS00,FM03,GLS}.

Finally, the shadowing process for the complete polysytem along the chain of cylinders is based on the same idea and
involves similar considerations regarding  the heteroclinic transitions.

With the constructive method used in the proof of {\bf Proposition~\ref{prop:tamecyl2}},
one can expect to obtain a quantitative control on the recurrence times of the Hamiltonian flow
(in the specific zones under consideration), and avoid the use of the Poincar\'e recurrence theorem
in specific models.
%at least in generic situations, in view of estimating diffusion times.
%This question will be more thoroughly investigated in further work, see the comments in Section~\ref{Sec:constructive}.

{\bf Theorem~\ref{thm:main1}} immediately follows from {\bf Theorem~\ref{thm:pseudoorb1}} and
{\bf Theorem~\ref{thm:shadowing}}.

\paraga
As mentioned above, the present paper considers the problem of generic drift in an abstract setting. Let us now describe
its possible applications and extensions.

\vskip1mm

\noindent{\bf 7.1.} The conditions we impose to our chains are also satisfied by the usual examples of {\em a priori unstable}
systems, limited to a single cylinder.
In this respect, the present work is an abstract approach of the problems considered in
\cite{ChierchiaG94,Bessi1996,Berti2003,B08,BKZ13,BT99,CY09,DLS00,DLS06a,DLS06b,FM03,
GT17,GL06,GR07,GR09,GR13,GLS,M02,T04,DS17}, amongst others. However, our assumptions
are less stringent than the usual ones.

\vskip1mm

\noindent{\bf 7.2.} Our approach can be applied to Hamiltonian systems that do not satisfy the strict convexity assumption mentioned earlier, but instead it satisfies the weaker condition that the return map to the sections is a twist map (hence, a tilt map).
Consider the following Hamiltonian on $\A^3$:
\begin{equation}\label{eqn:H-example-1} H(\th,r) = \frac{r_1^4}{4}+r_2+\frac{r_3^2}{2}+\eps\cos(2\pi\th_3)+\eps\mu f(\th),\end{equation}
with $\mu\ll 1$.
When $\mu=0$, the unperturbed Hamiltonian $h$ has a normally hyperbolic annulus ${\bf A}=\{(\theta,r)\,|\, \th_3=0, r_3=0\}$ diffeomorphic to $\A^2$.
The intersection of ${\bf A}$ with a level $H^{-1}({\bf e})$ is transverse and yields a cylinder $\jC$ diffeomorphic to $\T^2\times\R$, with coordinates $(\th_1,\th_2,  r_1)$.
The return map to the section $\Sig=\{\th_2=0\}$ inside this cylinder has the twist property but the global system is not
not convex (nor quasi-convex)  in the actions.

This system would not in principle be accessible to weak KAM theory or usual variational methods.

\vskip1mm

\noindent{\bf 7.3.}  Our method can be applied to Hamiltonian systems that satisfy the even weaker condition that the return map to the section is a  tilt map (rather than a twist map).
Consider the system \eqref{eqn:H-example-1} when $h$ is changed to
\[h(\theta,r) =\frac{r_1^4}{4}+r_2 +\frac{r_3^2}{2}+ a \cos (2\pi \th_1)+b\cos(2\pi\th_3).\]
The annulus ${\bf A}$ defined above is normally hyperbolic if  $a\ll b$.
Now consider the dynamics inside the cylinder $\jC$  as above.
The return map to the section $\Sig=\{\th_2=0\}$ inside $\jC$  is a tilt map but not a twist map when $a$ is large enough.  Such systems occur for (large) perturbations of integrable Hamiltonians
at the crossing of two independent simple resonances. Our method would still apply in this case, which is
``far from convex and far from perturbative''.

\vskip1mm

\noindent{\bf 7.4.}   We mainly developed the present method in view of a proof of the  Mather version of the Arnold conjecture (see \cite{Ma04,Ma12}) for a precise formulation of the problem). The main difficulty in this setting is to get ride of any quantitative estimate on the so-called ``splitting of separatrices'', since one cannot expect it to be bounded from below in generic situations. In particular, usual transversality arguments are not adapted to this setting and this is the main reason to shift to measure theoretic ones, as explained above.
See \cite{Mar4} for more details and ``simple'' examples.

%
%Before stating the result, we introduce some notation. For $\ka>0$, denote by $\mathscr{S}^\ka$ and
%$\mathscr{B}^\ka$ the unit sphere and the unit ball, respectively,  in $C^\ka_b(\A^3,\R)$, the Banach algebra of bounded, $C^\ka$-differentiable functions, endowed with the standard $\|\cdot\|_{\ka}$ norm.
%
%Given a ``threshold function'' $\eps_0:\mathscr{S}^\ka \to [0,+\infty[$, define the ``generalized ball''
%\[\mathscr{B}^\ka_{\eps_0}=\{\eps g\,|\, g\in \mathscr{S}^\ka,\,\eps\in ]0,\eps_0(g)[\, \}.\]
%The generalized ball $\mathscr{B}^\ka_{\eps_0}$  is open in $C^\ka_b(\A^3,\R)$
%provided that $\eps_0$ is lower-semicontinuous.
%
%\vskip3mm
%
%\noindent{\bf Mather version of the Arnold conjecture.} {\em
%Consider a $C^\ka$  integrable
%Tonelli Hamiltonian $h$ on $\A^3$, with $\ka \geq \ka_0$ large enough.
%Fix $e \geq \textrm{Min}(h)$ and a finite family of open sets $\hat{O}_1,\ldots,\hat{O}_m$ which intersect $h^{-1}(\e)$.
%
%Then there exists a lower semicontinous
%function
%\[\eps_0 : \mathscr{S}^\ka \to [0,+\infty[\]
%which takes strictly positive values on an open dense subset $\mathscr{D}$ of  $\mathcal{S}^\ka$,  and an open dense  $\mathscr{D}_{\eps_0}$ of $\mathscr{B}^\ka_{\eps_0}$, such that, for each $f\in\mathscr{D}_{\eps_0}$, the system
%\[H(\theta,r)=h(r)+f(\theta, r)\]
%admits an orbit which intersects each set $O_i=\T^3\times\hat{O}_i$, $i=1,\ldots,m$.}
%
%\vskip3mm

\vskip 1mm

\noindent{\bf 7.5.}  Non-hyperbolic cylinders naturally appear in quantitative questions related to diffusion. The main results in quantitative perturbation theory
have been obtained in the case of perturbations of convex or quasi-convex systems (e.g., optimality of the Nekhoroshev stability times in various functional
classes ($C^k$, $G^{\al,L}$, $C^\om$), optimality of the splitting of invariant tori). However, a much more natural geometric class for the unperturbed systems would be the steep one,
introduced by Nekhoroshev in \cite{N77}. The optimality of Nekhoroshev stability times for steep unperturbed systems is not proved yet, and a general study
clearly involves more general objects than normally hyperbolic ones. As for another example, regarding the problem of the optimality of the measure of the wandering
domains in exact-symplectic perturbations of integrable diffeomorphisms on the annulus $\A^n$, it turns out that the introduction of non-hyperbolic objects is
unavoidable (see \cite{LMS}). We hope that our work will motivate a more thorough study of these objects, probably using blowing-up methods from symplectic geometry.
\vskip1mm
{\em We emphasize that non-hyperbolic objects are in general difficult to control from the point of view of persistence theory, so that the aforementioned applications are (up to now) limited to the construction of particular examples for which the perturbations can be suitably chosen - as is the case for instance in the original Arnold example.}

\paraga We summarize the distinguished features of our approach, and emphasize the differences from other approaches.
\begin{itemize}
\item it uses homoclinic and heteroclinic ``maps'' that are only {\em locally defined}, and requires only very mild
topological transversality assumptions;
between the stable and unstable manifolds of the essential tori contained in the cylinders; this makes the method
well-adapted to passing to the singular limit where the ``hyperbolicity'' of the cylinders tends to $0$;
\item no requirement on the ``lower bounds of the splitting of tori'', nor that they for countable families, is required;
\item this latter singular limit generates important difficulties due to the fact that the cylinders have ``wild embeddings''
in the {\em a priori} stable case when the perturbation parameter tends to $0$, this makes the variational methods
difficult to apply, while the symplectic geometric ones (closely linked to our approach in this paper) are easier;
\item our method does not require the system restricted to the cylinder to be close to integrable,
since we only assume the return map to have a tilt property, in particular
it does not require any precise knowledge on the invariant objects for the dynamics
restricted to the cylinders (apart from their homology);
\item it does not require the cylinders to be normally hyperbolic, nor the inner dynamics to possess strong torsion property,
which allows us to expect easy applications to non-convex settings (as well as to more quantitative questions, see \cite{LMS})
\item our construction of pseudo-orbits can be made in %some sense
``algorithmic'' from a numerical point of view, which makes
it accessible to computer-assisted methods (some numerical approaches showing the existence of diffusing orbits along normally hyperbolic cylinders,
via computer assisted proofs, can be found in \cite{CG18,CGMMJ21}).
\end{itemize}

\paraga The structure of the paper is the following. Section~\ref{Sec:setting} is devoted to the description of the general setting
and the definition good cylinders and good chains of cylinders.
In Section~\ref{Sec:pseudoorb} we state and prove {\bf Theorem~\ref{thm:pseudoorb1}}.
{\bf Theorem~\ref{thm:shadowing}} and {\bf Theorem~\ref{thm:main1}} are proved in Section~\ref{Sec:shadowing}.
Section~\ref{Sec:constructive} is devoted to the algorithmic method to prove  {\bf Theorem~\ref{thm:pseudoorb1}}.
Finally, we  recall  in Appendix~\ref{App:twistmaps} some basic results on \twist maps, and in particular
a strong form of the Birkhoff transition lemma through a Birkhoff zone.

%%%%%%%%%%%%%%%%%%%%%%%%%%%%%%%%%%%%%%%%%%%%
%%%%%%%%%%%%%%%%%%%Section2%%%%%%%%%%%%%%%%%%%%%
%%%%%%%%%%%%%%%%%%%%%%%%%%%%%%%%%%%%%%%%%%%%

\section{The setting: good chains of cylinders}\label{Sec:setting}
In this section we make precise the dynamical features of  good cylinders and good chains.

%%%%%%%%%%%%%%%%%%%%%%%%%%%%%%%%%%%%%%%%%%%%%%
%%%%%%%%%%%%%%%%%%%%%%%%%%%%%%%%%%%%%%%%%%%%%%

\subsection{Tame cylinders}\label{Sec:cylinders}

The main objects in  \cite{Mar1,Mar2} are normally hyperbolic invariant cylinders  with boundary\footnote{{or normally
hyperbolic {\em singular} cylinders,
whose dynamical study will be reduced here to the non-singular case, see Section~\ref{sec:singular}}}, satisfying
some additional conditions. In this paper we adopt a slightly more general setting, relaxing the
normal hyperbolicity assumption but preserving all essential features of this situation
(existence of center-stable and center-unstable manifolds foliated by stable and unstable manifolds of points,
$\lambda$-lemma for arcs, and existence of an invariant measure).
We call here {\em tame cylinders} the resulting objects, for which we give a
formal definition gathering together the necessary properties.

Some examples of tame cylinders that are not normally hyperbolic, as well as
applications to the problem of wandering domains, can be found in \cite{LMS}.
Other examples  appear in concrete systems from Celestial Mechanics, see \cite{GMSS}.
Also, an example in the case of  a rotator-pendulum type system is presented at the end of this section.

\paraga  Let $X$ be a complete vector field on a smooth manifold $M$.
We say that $\jC\subset M$ is an {\em invariant cylinder} for $X$ if $\jC$ is a $C^1$--submanifold  of
$M$,  $C^1$--diffeomorphic to $\T^2\times [0,1]$, such that $X$ is everywhere tangent
to $\jC$ and is moreover tangent to $\partial \jC$ at each point of~$\partial \jC$.
Note that $\jC$ is compact and invariant under the flow of $X$.

In the following we consider  $M=\A^3$, a  $C^\ka$-Hamiltonian function $H:\A^3\to\R$, with $\ka\geq 2$.  We denote the
Hamiltonian vector field associated to $H$ by $X_H$, and its Hamiltonian flow by $\Phi_H$.
Recall that we only consider energies  $\e$ and levels $H^{-1}(\e)$ which are completely regular.
The following definition is ``minimal'' in the sense that the requirements on the regularity of the various objects
are the less stringent ones in view of our future constructions. As a consequence, our definition has no relation
with one could expect form the point of view of persistence theory applied in the normally hyperbolic case, where
the regularity of the invariant manifolds and foliations are directly related to those of the system and spectral data of the invariant object at hand.

\begin{Def}\label{def:invcyl}
Let $H$ be a $C^\ka$ Hamiltonian function on $\A^3$,  with $\ka\geq 2$, and fix an energy~$\e$. We denote by $d$ the usual distance on $\A^3$.
We define a {\em tame cylinder at energy~$\e$} for $H$ as a cylinder $\jC$
contained in $H\inv(\e)$ and invariant under $\Phi_H$, which satisfies the following properties:
\begin{itemize}
\item There exists a 5-dimensional manifold with boundary $U\subset H\inv(\e)$ containing $\jC$,
with $\inte\jC\subset\inte U$ and $\d\jC\subset \d U$,
such that the subsets
\beq
W_U^\pm(\jC)=\Big\{x\in U\mid \Phi^t_H(x)\in U,\  \forall t\in\R^\pm,\ \text{and}\ \lim_{t\to\pm\infty} d\big(\Phi_H(x), \jC\big)=0\Big\}
\eeq
are $4$-dimensional $C^{1}$ embedded submanifolds of $U$, with boundaries contained in $\d U$. We fix $U$ once and for all
and set
\beq
W^-(\jC)=\bigcup_{t\geq0}\Phi_H^t\big(W_U^-(\jC)\big),\qquad
W^+(\jC)=\bigcup_{t\leq0}\Phi_H^t\big(W_U^+(\jC)\big),
\eeq
which are $\Phi_H$-invariant $C^{1}$ immersed submanifolds of $H\inv(\e)$.
\item There exist $C^1$-diffeomorphisms
\beq
J^\pm:\jC\times \,]-1,1[\ \to W_U^\pm(\jC)
\eeq
such that $J^\pm(x,0)=x$ for $x\in\jC$,  and, setting
\beq
W_U^\pm(x)=J^\pm\big(\{x\}\times\,]-1,1[\big),\qquad x\in\jC,
\eeq
then $J^\pm(x,\cdot):\,]-1,1[\,\to W_U^\pm(x)$ is a $C^1$ diffeomorphism, and
\beq
\forall y\in W_U^\pm(x), \quad \lim_{t\to\pm\infty}d\big(\Phi_H^t(x), \Phi_H^t(y)\big)=0.
\eeq
\item We set
\beq\label{eq:defcentman}
W^-(x)=\bigcup_{t\geq0}\Phi_H^t\big(W_U^-(\Phi_H^{-t}(x))\big),\qquad
W^+(x)=\bigcup_{t\leq0}\Phi_H^t\big(W_U^+(\Phi_H^{-t}(x))\big),
\eeq
which are $C^1$ immersed $1$-dimensional submanifolds of $H\inv(\e)$.

\item The manifolds  $W^\pm(x)$  satisfy
the equivariance property
\beq
\Phi_H^t\big(W^\pm(x)\big)=W^\pm\big(\Phi_H^t(x)\big),\qquad \forall x\in\jC,\ \forall t\in\R.
\eeq

\item There is a  negatively invariant neighborhood $\jN^-\subset W_U^-(\jC)$ of $\jC$ in $W^-(\jC)$ such that for
$x\in\jC$ and  $y\in  \jN^-\cap W^-(x)$,
\beq\label{eq:decreasing}
d\big(\Phi_H^{-t}(x), \Phi_H^{-t}(y)\big)\leq d(x,y),\qquad \forall t\geq0,
\eeq
and a positively invariant neighborhood $\jN^+\subset W_U^+(\jC)$ of $\jC$ in $W^+(\jC)$ such that for
$x\in\jC$ and  $y\in  \jN^+\cap W^+(x)$,
\beq\label{eq:decreasing}
d\big(\Phi_H^{t}(x), \Phi_H^{t}(y)\big)\leq d(x,y),\qquad \forall t\geq 0.
\eeq

\item {\bf The $\la$-property.}  For $x\in\jC$, set $J^-_x=J^-(x,\cdot):\,]-1,1[\ \to W_U^-(x)$.
Fix $x\in\jC$.  Then for any $1$-dimensional $C^1$ submanifold $\De$ of $H\inv(\e)$ which intersects
$W^+(\jC)$ transversely in $H\inv(\e)$ at  $\xi\in W^+(x)$,
there exist a family $(\De_t)_{t\geq t_0}$ of submanifolds of $\De$ containing $\xi$
and $C^1$ parametrizations $\ell_t:\,]-1,1[\,\to H^{-1}(\bf e)$ of the images $\Phi_H^t(\De_t)$ such that:
\beq
\lim_{t\to+\infty}\norm{\ell_t-J^-_{\Phi_H^t(x)}}_{C^0(]-1,1[)}=0.
\eeq
\item There exists a $\Phi_H$-invariant Borel measure $\bmu$ on $\jC$ such that $\bmu(O)>0$
for any nonempty open subset $O$ of $\jC$.

\item The cylinder $\jC$ is contained in the interior $\inte \hc$ of a cylinder $\hc$ which satisfies the previous
six conditions, with natural continuation conditions (that is, with obvious notation, $U\subset \ha U$,
$\ha J^\pm_{\vert U}=J^\pm$, $\jN^-\subset\ha\jN^-$). Any such cylinder $\hc$ is said to be a
continuation of $\jC$.
\end{itemize}
\end{Def}

Normally hyperbolic invariant cylinders with boundary are particular cases of tame cylinders
\footnote{In particular, the cylinders consider in \cite{Mar1} are readily seen to be tame cylinder, by usual normal hyperbolicity arguments.}.

\vskip1mm

Here we use the terminology of stable/unstable manifolds to refer to the leaves $W^\pm(x)$ associated to points $x\in \jC$, and of center-stable/unstable  manifolds\footnote{This terminology is used, for instance, in \cite{BKZ00}; it is however more common to refer
to the global manifolds as stable/unstable  manifolds, rather than center-stable/unstable  manifolds}
to refer to the invariant manifolds $W^\pm(\jC)$ associated to the cylinder $\jC$.

\vskip 1mm

\vskip1mm
In the setting of Definition~\ref{def:invcyl}, the global manifolds $W^\pm(\jC)$ could depend on the choice of $U$, but this will
cause no trouble in what  follows.

\paraga We set out now  some remarks and conventions. Given a subset $A$ of a tame cylinder $\jC$, we set
\beq
W^\pm(A)=\bigcup_{x\in A} W^\pm(x).
\eeq

Thus $W^\pm(A)$ are invariant when $A$ is invariant. This will be the case in particular when $A=\inte \jC$.

\vskip1mm

Observe that since the parametrization $J^+$ is a $C^1$  diffeomorphisms, the (stable) characteristic projection
\beq\label{eq:charproj}
\Pi^+: W^+(\jC)\to \jC
\eeq
which to a point $y\in W^+(\jC)$ associate the unique point $x$ of $\jC$ such that $y\in W^+(x)$,
is also $C^1$. A similar property holds  for the (unstable) characteristic projection $\Pi^-: W^-(\jC)\to \jC$.

\vskip1mm

\noindent{\bf Convention.} {\em In the following, given a tame cylinder $\jC$, we will choose once and for all one continuation
of $\jC$, which we always denote by $\hc$.}

\paraga {\bf An example of a tame cylinder.} Let us give an example, similar  to that in \cite{LMS},
of a system  that has a tame cylinder which is not normally hyperbolic.
Consider the following variant of the Arnold example on $\A^3$:
\beq\label{eq:ArnoldHam}
H_\eps(\th,r)=\pdemi(r_0^2+r_1^2+r_2^2)-(\cos(2\pi\th_2)-1)^2+\eps f(\th),\qquad (\th,r)\in\A^3,
\eeq
with $\eps\geq 0$ and where $\Supp f\subset \big\{\th_2\in[1/4,3/4]\big\}$.
When $\eps=0$, the system is the uncoupled product of the integrable system generated by $\pdemi(r_0^2+r_1^2)$
on $\A^2$ and the pendulum-like system $\pdemi r_2^2-(\cos(2\pi\th_2)-1)^2$ on $\A$. This latter system admits
the  fixed point $O$ of coordinates $(\th_2,r_2)=(0,0)$. Clearly $O$ is degenerate since  the expansion of the Hamiltonian in the
neighborhood of $O$ reads $\pdemi r_2^2-4\pi^4\th_2^4+\cdots$.
Note that $O$ admits well-defined stable and unstable manifolds, and the pendulum-like system satisfies the usual
$\lambda$-lemma by easy $2$-dimensional arguments.

The annulus $\jA=\A^2\times\{O\}$ is  invariant under the flow, but not normally hyperbolic. This is also the case
when $\eps>0$, thanks to the condition on the support of $f$.

For $\eps$ small enough, the level $H_\eps\inv(0)$ is compact and regular, and its intersection with $\jA$ is
the $3$-dimensional torus
$$
\jT=\big\{(\th,r)\in\A^3\mid (\th_2, r_2)=O,\ \pdemi(r_0^2+r_1^2)=1\big\},
$$
which is invariant in $H_\eps\inv(1)$ but not normally hyperbolic. However $\jT$  clearly admits global invariant manifolds, which are the product of $\jT$ with the stable and unstable manifolds of $O$. In the same way, the
stable and unstable manifolds of the points of $\jT$ are the product of those points with the stable and unstable manifolds of $O$.
They are uniquely defined.

Moreover,
given a point ${\bf a}=(a_0,a_1)$ on the circle $C$ of equation $\pdemi(r_0^2+r_1^2)=1$, the $2$-dimensional torus
$$
T_{{\bf a}}=\big\{(\th,r)\in\A^3\mid (\th_2, r_2)=O,\  (r_0, r_1)={\bf a}\big\}\subset \jT
$$
is invariant under the flow. Varying $\bf a$ yields a foliation of $\jT$ by invariant tori. Given now any compact connected
arc $C_{\bf a}^{{\bf a}'}\subset C$  with extremities ${\bf a}$ and ${\bf a}'$, the subset
$$
\jC=\big\{(\th,r)\in\A^3\mid (\th_2, r_2)=O,\ (r_0,r_1)\in C_{\bf a}^{{\bf a}'}\big\}
$$
immediately satisfies our definition of a tame cylinder at energy $1$ for $H_\eps$: the two
components of its boundary are the  invariant tori
$T_{\bf a}$ and $T_{{\bf a}'}$. The existence of continuations follows from the choice of larger compact arcs
$C_{\bf b}^{{\bf b}'}\supset C_{\bf a}^{{\bf a}'}$. The $\la$-property follows from easy two-dimensional considerations
around $O$. The invariant measure $\bmu$
is the Liouville measure on $H\inv(0)\cap \jA$ deduced from the induced symplectic form on $\jA$.

\vskip1mm

The interest in such examples comes from the fact that the distance between two consecutive periodic invariant curves  in the ``pendulum'' part of \eqref{eq:ArnoldHam} is polynomially small with respect to the period, while it is exponentially small in the case of the usual simple pendulum. This is of crucial importance in constructing examples of perturbed initially stable systems with ``large'' wandering domains (see \cite{LMS}).

%%%%%%%%%%%%%%%%%%%%%%%%%%%%%%%%%%%%%%%%%%%%%%%%%%%%%%%%%%%%%%%%%%%%
%%%%%%%%%%%%%%%%%%%%%%%%%%%%%%%%%%%%%%%%%%%%%%%%%%%%%%%%%%%%%%%%%%%%

\subsection{Tilt sections, homoclinic correspondences and splitting arcs}

We introduce in this section the main ingredients to
reduce our problem to a two-dimensional one.

\paraga {\em Two-dimensional tilt sections inside a cylinder.} We first introduce the notion of a {\em tilt section}, which will be our main tool
to reduce our problem to the Birkhoff
theory in dimension 2.
Given $a<b$ we set:
$$
\bA(a,b)=\T\times[a,b],\qquad
\Ga(a)=\T\times\{a\},\qquad
\Ga(b)=\T\times\{b\}.
$$
We often abbreviate $\bA(a,b)$ to $\bA$ when there is no risk of confusion.

\vskip2mm

We refer to Appendix~\ref{App:twistmaps}  for the usual definition  of a \twist map, which we always assume to  tilt the vertical to the right (or always to the left).
We denote by $\Ess(\ph)$ the set of essential invariant circles of an area-preserving \twist map $\ph:\bA\righttoleftarrow$.
Elements $\Gamma\in\Ess(\ph)$ are graphs of uniformly Lipschitz functions $\ell_{\Gamma}:\T\to[a,b]$, by the Birkhoff theorem \cite{HF,Ma84,Y92}.
We endow $\Ess(\ph)$ with the Hausdorff topology or, equivalently, with the uniform $C^0$ topology on the corresponding functions $\ell_\Gamma$.

\begin{Def}\label{def:good}
We say that an area-preserving \twist map $\ph$ of $\bA$ is {\em good} if
\begin{itemize}
\item  $\ph$ does not admit any essential invariant circle with rational rotation number;
\item  The boundaries $\Ga(a)$ and $\Ga(b)$ are both dynamically minimal;
\item  Each boundary $\Ga(a)$, $\Ga(b)$ is accumulated by a sequence  of dynamically minimal elements of $\Ess(\ph)$.
\end{itemize}
\end{Def}

\vskip1mm

The last two conditions are always satisfied in most cases of interest (twist maps on annuli with Diophantine boundaries, under the KAM theorem conditions).

\vskip1mm

Given $\Ga\in \Ess(\ph)$,
$\Ga^-$ denotes the connected component of $\bA\setm\Ga$ located below $\Ga$
in $\bA$, and $\Ga^+$ denotes the connected component of $\bA\setm\Ga$ located above $\Ga$
in $\bA$.
In the following we will crucially use the following result, which is proved in Appendix~\ref{App:twistmaps}.

\begin{lemma}\label{lem:accum}
Let $\ph$ be a good area-preserving \twist map $\ph$ of $\bA$.
Then the following properties hold true:
\vskip1mm
(i) Any two distinct elements of $\Ess(\ph)$ are  disjoint, so that the set $\Ess(\ph)$
admits a natural order $\leq$ given by $\Gamma_1\leq\Gamma_2$ if $\ell_{\Gamma_1}(\theta)\leq \ell_{\Gamma_2}(\theta)$  for all $\theta\in \mathbb{T}$, where $\Gamma_j$ is the graph of $\ell_{\Gamma_j}$, $j=1,2$;
\vskip1mm
(ii) Given a nonempty subset $E\subset\Ess(\ph)$, the greatest lower bound $\Inf E$ and the least upper bound $\Sup E$ relative to $\leq$
exist;
\vskip1mm
(iii) Every invariant essential circle $\Ga\subset\big(\bA\setm\Ga(a)\big)$ (resp., $\Ga\subset\big(\bA\setm\Ga(b)\big)$)
is either  the upper (resp., lower) boundary of a Birkhoff zone of $\ph$, or is accumulated
by a sequence of elements of $\Ess(\ph)$ located in $\Ga^-$ (resp. $\Ga^+$).
\end{lemma}

We can now introduce our main notion. We consider a $C^2$ Hamiltonian function $H$ on $\A^3$, we fix an energy
$\e$ and a tame cylinder $\jC$ at energy $\e$. Given a topological space $X$ and two subsets $A\subset B$
of $X$ with $A$ connected, we denote by
$$
\CC(B,A)
$$
the connected component of $B$ which contains $A$.

\begin{Def}\label{def:twistsec}
A {\em \twist section} for $\jC$ is a quadruple $(\Sig,\bA,\chi,\ph)$ such that:
{\begin{itemize}
\item $\Sig$ is a global Poincar\'e section for  the flow $(\Phi_H)_{\vert\jC}$, with return map $\ph$;
\item $\chi$ is a $C^1$ embedding of $\bA$ in $\jC$ with image $\Sig=\chi(\bA)$;
\item $\chi\inv\circ\ph\circ\chi$ is a good area-preserving \twist map of $\bA$.
\end{itemize}}
Given such a \twist section on $\jC$, we set
\beq
\d_\bu\Sig=\chi\big(\Ga(a)\big), \qquad \d^\bu\Sig=\chi\big(\Ga(b)\big).
\eeq
These two circles are contained in the two boundary components of
$\d\jC$. We set
\beq
\d_\bu\jC=\CC\big(\d\jC,\d_\bu\Sig\big),\qquad \d^\bu\jC=\CC\big(\d\jC,\d^\bu\Sig\big).
\eeq
We finally define a continuation of $(\Sig,\bA,\chi,\ph)$ for the continuation $\hc$
as a \twist section $(\ha\Sig,\ha\bA,\ha\chi,\ha\ph)$
for $\hc$ which continues the previous one in the natural way.
\end{Def}

\paraga {\em Homoclinic correspondences.}
The homoclinic correspondences that we consider will always be families of locally defined diffeomorphisms
(with possibly intersecting domains). The closely related  {\em scattering maps}
(see \cite{DLS08} and references therein) are local  diffeomorphisms on cylinders,
while our homoclinic diffeomorphisms are local diffeomorphisms on tilt sections (see also \cite{DGR13}).
This is a crucial step in order to reduce our approach to two-dimensional dynamics.

\vskip1mm

$\bu$ We first recall the usual context for scattering maps and examine the transversality properties of
the center-stable and center-unstable manifolds of the cylinder.  Given any manifold with boundary $L$, we write
$\inte L=L\setm\d L$.

\begin{Def} We consider a $C^2$ Hamiltonian function $H$ on $\A^3$, we fix an energy~$\e$ and consider a tame cylinder $\jC$ at energy $\e$.
We define the {\em transverse homoclinic intersection}
of the continuation $\hc$ as the set
\beq
\Homt(\hc)\subset W^+(\inte \hc)\cap W^-(\inte\hc)
\eeq
formed by the points $\xi$ such that
\beq\label{eqn:trans_fibers}
W^-\big(\Pi^-(\xi)\big)\trans_\xi W^+(\inte\hc)
\quad\textrm{and}\quad W^+\big(\Pi^+(\xi)\big)\trans_\xi W^-(\inte\hc),
\eeq
where $\Pi^\pm$ are the characteristic projections defined in (\ref{eq:charproj}) and where
 $\trans_\xi$ stands for ``intersects transversely at $\xi$ relatively to $H\inv(\e)$''.
\end{Def}

Note that $W^\pm(\inte\hc)=\inte W^\pm(\hc)$, following our definitions.
The following lemma is an immediate consequence of the
implicit function theorem.

\begin{lemma}\label{lemma:scatteringmap}
Fix $\xi\in \Homt(\hc)$ and
write $x^\pm=\Pi^\pm(\xi)\in\hc$.
Then the following hold:
\begin{itemize}
\item[(i)] $W^+(\hc)$ and $W^-(\hc)$ intersect transversely
at $\xi$ in $H\inv(\e)$;
\item[(ii)] There exists a $3$-dimensional open neighborhood $\jO$ of $\xi$ in $W^+(\inte\hc) \cap W^-(\inte\hc)$
such that the condition  \eqref{eqn:trans_fibers} is satisfied at all points $\eta\in \jO$,
and, moreover,
\begin{equation}\label{eqn:homoclinic_channel}
\begin{split}
&W^-\big(\Pi^-(\xi)\big)\trans_\eta \jO \textrm { relative to } W^-(\inte\hc),
 \textrm{ and } \\
&W^+\big(\Pi^-(\xi)\big)\trans_\eta \jO \textrm { relative to } W^+(\inte\hc).
\end{split}
\end{equation}
\item[(iii)] Provided that the open neighborhood $\jO$ of
$\xi$ in $W^+(\inte \hc)\cap W^-(\inte\hc)$ is chosen sufficiently small,
there exist open neighborhoods $O^\pm$ of $x^\pm$ in $\inte \hc$ such that the restrictions $\Pi^\pm_{\vert \jO}$ are $C^1$
diffeomorphisms from $\jO$ onto $O^\pm$.
\end{itemize}
\end{lemma}

In \cite{DLS08} the open set $\jO$  is referred to as a {\em homoclinic channel}, and the map
\[S=\Pi^+\circ(\Pi^-_{\jO})^{-1}:O^-\to O^+\]
as the  {\em scattering map} associated to  $\jO$.

\vskip2mm

Our main variation relative to the usual definition of a scattering map  consists in a {\em two-dimensional reduction}
by transport to the section $\Sig$ via the Hamiltonian flow. We also allow for the
transversality condition (\ref{eqn:homoclinic_channel}) to hold only almost everywhere,
which prove to be useful to be enough in our subsequent constructions
of drifting orbits)\footnote{
In particular, we will be able to deal with points $x^-\in \jC$ and $S(x^-)=x^+\in \jC$ with
$\xi\in W^-(x^-)\cap W^+(x^+)$ such that
$W^+(\hc)$, $W^-(\hc)$ are not transverse at $\xi$ in $H\inv(\e)$}.

\vskip2mm

Let us first informally describe our two-dimensional reduction.
Consider an energy  $\e$, and suppose that $\jC$ is a tame cylinder with continuation $\hc$ and
continued \twist section $\ha\Sig$  in the corresponding level.
Assume that we are given an (in general uncountable) family of scattering maps on~$\hc$
\beq\label{eqn:Si}
S=(S_i)_{i\in I},  \qquad S_i:\Dom S_i\to \Im S_i,
\eeq
where $\Dom S_i\subset\inte\hc$ and $\Im S_i\subset\inte\hc$ are open subsets and where
each $S_i$ is a {\em measure preserving}\footnote{in the normally hyperbolic situation, this comes from the definition of the Borel measure of $\jC$ together with the symplectic properties of its stable and unstable manifolds} $C^1$-diffeomorphism from  $\Dom S_i$ to $\Im S_i$. Moreover, assume that for every $i\in I$
there is an open subset $\Domt\, S_i\subset \Dom S_i$, with full measure in $\Dom S_i$,
such that
\beq\label{eqn:ImSi}
\Im S_i\cap\inte \ha\Sig\neq\emptyset.
\eeq
and
\beq\label{eq:redhom3}
\forall y\in \Domt S_i,
\quad
W^-(y)\cap
W^+\big(S_i(y)\big)\cap
\Homt(\hc)\neq\emptyset.
\eeq

Assume finally that there exists an open subset $D_i\subset \inte \ha\Sig$ and a  non-negative  $C^1$  function
$\tau_i: D_i \to \Dom S_i$ such that for all $x\in D_i$
\beq\label{eqn:tauix}
\Phi^{\tau_i(x)}_H(x)\in \Dom S_i,\qquad   S_i\big(\Phi^{\tau_i(x)}_H(x)\big)\in \Im S_i\cap\inte \ha\Sig
\eeq
(this is not restrictive : such a function always exists if $\Dom S_i$ is small enough,
since $\ha\Sig$ is a global Poincar\'e section).
%, for each point $x'\in \Dom S_i$ the backwards orbit $\Phi_H^{-t}(x')$, with $t>0$,  intersects $\ha\Sig$  for the first time at some point $x=\Phi_H^{-\tau}(x')$. The time $\tau$ and the point $x'$ depend on $x\in \ha\Sig$ in a $C^1$ fashion. Conversely, there exists an open set in $\ha\Sig$, which we denote by $\Dom \xi_i$, and a  non-negative  $C^1$  function $\tau_i:\Dom \xi_i \to \Dom S_i$ such that for all $x\in \Dom \xi_i$
%\beq\label{eqn:tauix}
%\Phi^{\tau_i(x)}_H(x)\in \Dom S_i.
%\eeq
Note that the function $\tau_i$ together with is open domain $\Dom \tau_i$ need not be uniquely determined.
This description leads to the following precise definition.

\begin{figure}[h]
\centering
\includegraphics[width=0.45\textwidth]{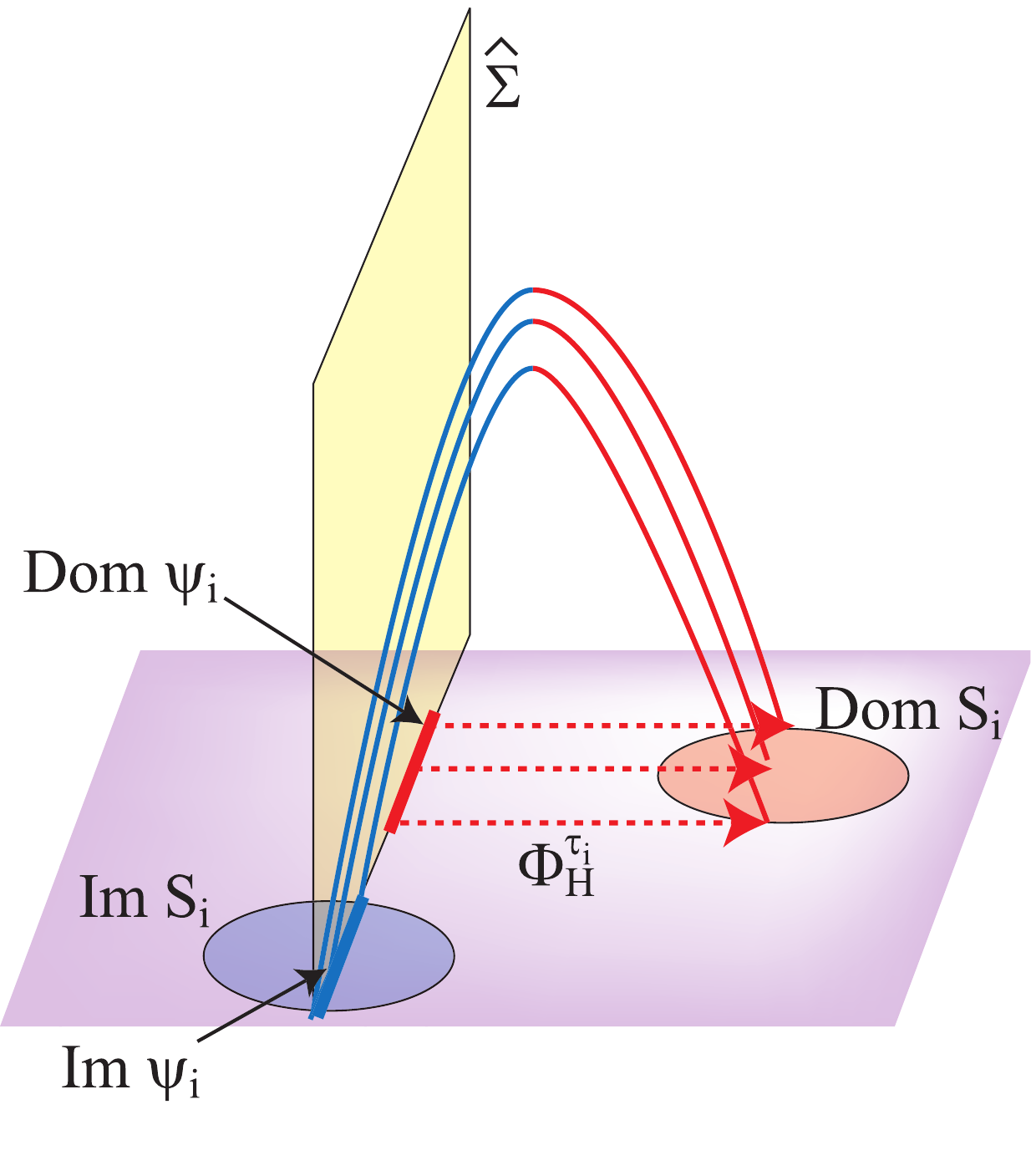}
\caption{Homoclinic correspondence  on the cylinder and the induced homoclinic correspondence  on the Poincar\'e section.}
\label{fig:Homoclinic_map_section}
\end{figure}

\begin{Def}\label{def:scattering}
We consider a $C^2$ Hamiltonian function $H$ on $\A^3$, and we fix an energy
$\e$ and a tame cylinder $\jC$ at energy $\e$, with continuation $\hc$ and
continued \twist section $\ha\Sig$.

A {\em homoclinic correspondence} associated with these data is an (in general uncountable) family
of $C^1$ local diffeomorphisms of $\inte\ha\Sig$:
\beq
\psi=(\psi_i)_{i\in I},  \qquad \psi_i:\Dom \psi_i\to \Im\psi_i,
\eeq
where $\Dom \psi_i\subset \inte\ha\Sig$ and $\Im\psi_i\subset \inte\ha\Sig$ are open subsets, for which there exists
a family of measure preserving homoclinic maps
$S=(S_i)_{i\in I}$, $S_i:\Dom S_i\to \Im S_i$ on $\hc$ such that
$S_i$ and $\psi_i$ satisfy  \eqref{eqn:Si}, \eqref{eq:redhom3}, \eqref{eqn:ImSi}, \eqref{eqn:tauix},
where
\beq\label{eqn:xiSi}
\forall x\in\Dom\psi_i,\qquad  \psi_i(x):=S_i\Big(\Phi_H^{\tau_i(x)}(x)\Big).
\eeq
\end{Def}
%See Figure~\ref{fig:Homoclinic_map_section}.

Homoclinic correspondences are not uniquely defined, and the domains $\Dom \psi_i$ (resp. $\Dom S_i$)
are not necessarily pairwise disjoint. In the following we
identify $\ha\Sig$ with $\ha\bA=\mathbb{T}\times[\ha a,\ha b]$ via the embedding $\ha\chi$, and we
indifferently consider our homoclinic correspondences as defined on $\ha\Sig=\ha\chi(\ha\bA)$ or on $\ha\bA$.

\paraga {\em Splitting arcs and associated domains.} An {\em arc} in $\ha\bA$ is a continuous {injective} map
$\ze:[0,1]\to\ha\bA$, and we write $\til\ze=\ze([0,1])\subset\ha\bA$ for its image.
Denote by  $\pi:\ha\bA\to\T$ the projection onto the $\th$-coordinate.
Given two points $\th,\th'$ of $\T$ at a distance less than  $1/2$, we write $[\th,\th']\subset\T$ for the unique segment
bounded by $\th$ and $\th'$ according to the natural orientation of $\T$.
In the following definition we implicitly assume that the pairs points in $\T$ that we consider are always close
enough so that the previous convention  applies.

\begin{Def}\label{def:splittingarc} We consider a tame cylinder $\jC$ with section $(\Sig,\bA,\chi,\ph)$
and continuations $\hc$, $(\ha\Sig,\ha\bA,\ha\chi,\ha\ph)$, respectively. We let  $\psi=(\psi_i)_{i\in I}:\ha\Sig\righttoleftarrow$ be a homoclinic
correspondence. %We identify $\ha\Sig$ with $\ha\bA=[\ha a,\ha b]\times \mathbb{T}$ via the embedding $\ha\chi$.
Fix $\Ga\in\Ess(\ha\ph)$ contained in $\ha\bA\setm\Ga(\ha a)$ and let $\al_0$ be a point in $\Ga$.
\begin{itemize}
\item A {\em splitting arc} based at $\al_0$ %for these data
is an arc $\ze$ in $\ha\bA$ for which
$$
\ze(0)=\al_0,\quad \ze(]0,1])\subset \Ga^-;\quad \exists i\in I,\ \ze(]0,1])\subset \Dom\psi_i,\quad \psi_i(\ze(]0,1]))\subset \Ga.
$$
\item A {\em non-vertical splitting arc}
is a splitting arc $\ze$  for which there exists a sequence $(s_n)_{n\in\N}$
in $]0,1]$ with $\lim_{n\to\infty} s_n=0$ such that
$$
\pi\big(\ze(s_n)\big)\neq \pi\big(\ze(0)\big),\qquad \forall n\in\N.
$$

\item A {\em right (left) splitting arc}
is a splitting arc $\ze$
for which there exists a sequence $(s_n)_{n\in\N}$
in $]0,1]$ with $\lim_{n\to\infty} s_n=0$ such that
$$
\pi\big(\ze(s_n)\big)>\pi\big(\ze(0)\big),\qquad (\textrm{resp. } \pi\big(\ze(s_n)\big)<\pi\big(\ze(0)\big) ),\qquad \forall n\in\N.
$$

\end{itemize}
\end{Def}

%%%%%%%%%%%%SplittingArc%%%%%%%%%%%%%%

\begin{figure}[h]
\centering
\includegraphics[width=0.7\textwidth]{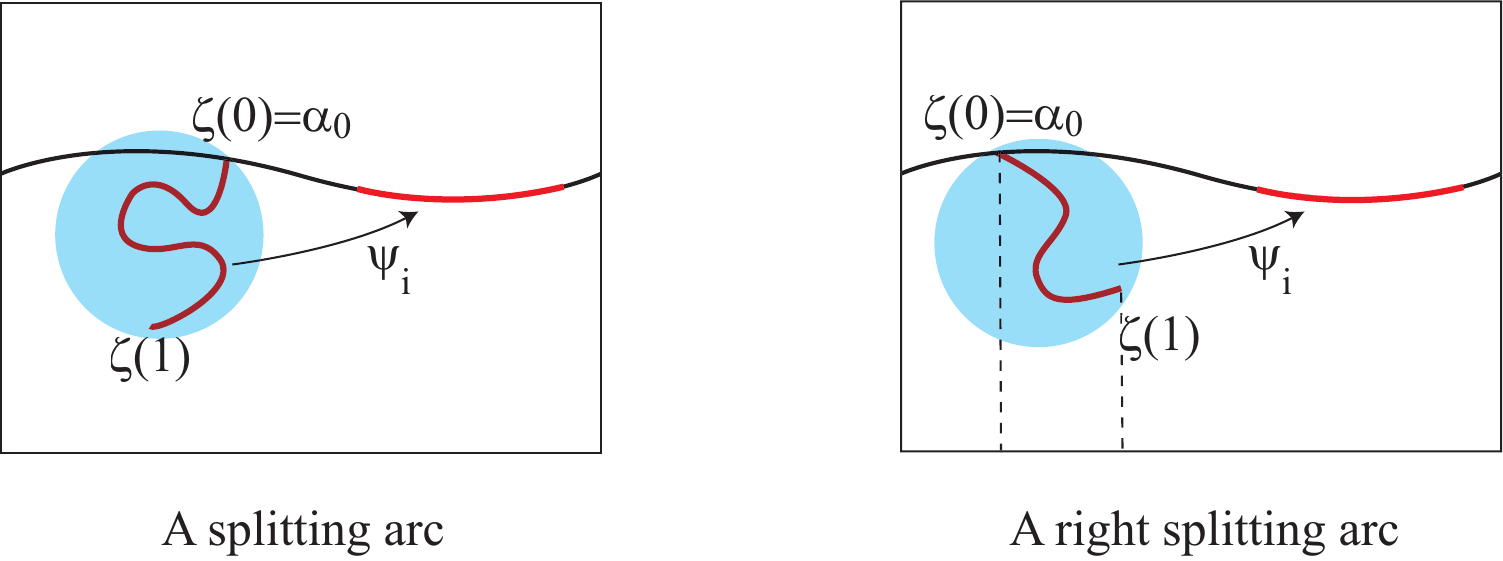}
\caption{Splitting arcs.}
\label{fig:splitting_arcs}
\end{figure}

See Figure \ref{fig:splitting_arcs}. We say that a splitting arc is based on $\Ga$ when it is based at some point of $\Ga$.

\vskip1mm

\begin{rem}\label{rem:nonvert} Here we collect several remarks on splitting arcs.
\begin{itemize}

\item If $\ze$ is a right (resp. left) splitting arc, then the restriction $\ze_{\vert [0,s_*]}$ (up to reparametrization) is also a right (resp. left) splitting arc, for any $0\leq s^*\leq 1$.

\item A sufficient condition for a splitting arc $\ze$ to be a right (left) splitting arc   is that it admits a derivative $\ze'(0)=(u,v)$ with $u>0$ (resp. $u<0$).

\item A non-vertical splitting arc is either a right splitting arc, or a left splitting arc, or both.

\item For a vertical splitting arc  (that is, a splitting arc for which there exists $s_*>0$ such that  $\pi\big(\ga(s)\big)=\pi\big(\ga(0)\big)$  for all $0\leq s\leq s_*$),   its image under the \twist map $\ha\ph$ is a left arc (resp., its image under $\ha\ph^{-1}$ is a right arc); this is not necessarily a splitting arc, since for instance $\ha\ph(\tilde\ze)$ (resp. $\ha\ph^{-1}(\tilde\ze)$) may not be contained in any
$\Dom\psi_i$, $i\in I$.
\end{itemize}

\end{rem}

\vskip1mm

Given a point $x_0=(\th_0,r_0)$ in $\ha\bA$, we denote
by
$$
V(x_0)=\big\{(\th_0,r)\mid r\in[\ha a,\ha b]\big\}
$$
the vertical through  $x_0$ in $\ha\bA$,  by
$$
V^-(x_0)=\big\{(\th_0,r)\mid r\in[\ha a,r_0]\big\}
$$
the semi-vertical below $x_0$ in $\ha\bA$, and by
$$
V^+(x_0)=\big\{(\th_0,r)\mid r\in[r_0,\ha b]\big\}
$$
the semi-vertical above $x_0$ in $\ha\bA$.

\begin{Def} \label{def:domain} Let $\Ga\in\Ess(\ha\ph)$ contained in $\ha\bA\setm\Ga(\ha a)$, be the graph of the continuous function
 $\ga:\T\to[\ha a,\ha b]$ and $\al_0\in\Ga$.

Let $\ze$ be a right splitting arc based at $\al_0=\ze(0)$, let  $\al_*$ be a point in $\Ga$ such that $$\pi (\al_0)<\pi (\al_*)<\textrm{max}_{s\in[0,1]}\pi (\ze(s)),$$
and let $\be_*=\ze(s_*)$ be the point in $V^-(\al_*)\cap  \til\ze$ with the maximal  $r$-coordinate. Let $C$ be the Jordan curve
formed by the concatenation of the arcs $\ze([0,s_{*}])$, $[\be_*,\al_*]$, and $[\al_*,\al_0]\subseteq \Ga$.
We denote by $D(\ze_{\mid [0,s_{*}]})$ the bounded connected component of the complement of $C$ in $\T\times \R$.
We say that $D(\ze_{\mid [0,s_{*}]})$ is   a {\em triangular domain associated with $\ze$}.
We define a triangular domain associated with a left splitting arc similarly.
\end{Def}

\begin{figure}[h]
\centering
\includegraphics[width=0.5\textwidth]{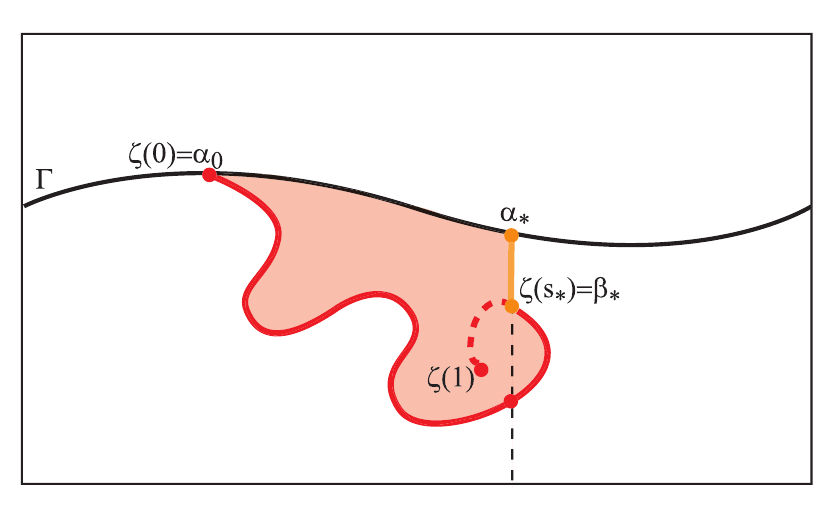}
\caption{A triangular domain associated to a right splitting arc.}
\label{fig:splitting_arc_domain}
\end{figure}

Let us now define the notion of {\em tilted arcs} (see Appendix~\ref{App:twistmaps} for more details).
Given $u,v$ two vectors in $\R^2$, let $\angle(u,v)$ be the  angle from $u$ to $v$ in $[0,2\pi[$, measured counterclockwise. Denote $\partial_r=(0,1)\in\mathbb{R}^2$.

\vskip1mm

A  $C^1$ arc emanating from $\Ga$ with $\eta'(s)\neq 0 $ for $s\in[0,1]$ is said to be
{\em positively tilted} (resp. {\em negatively tilted}) when \begin{itemize}
\item $\angle\big(\partial_r,\eta'(0)\big)\in\,]0,\pi[$
(resp. $\angle\big(\partial_r,\eta'(0)\big)\in\,]-\pi,0[$),
and \item the continuous lift  $s\in[0,1]\mapsto\angle\big(\partial_r,\eta'(s)\big)\in\R$
is positive (resp. negative) for all $s\in [0,1]$.
\end{itemize}

\vskip1mm

The main property of the triangular domains defined above is the following (see Appendix~\ref{App:twistmaps}
for a proof).

\begin{lemma}\label{lem:tiltarc}  Consider an essential invariant circle $\Ga\in\Ess(\ha\ph)$
contained in $\ha\bA\setm\Ga(\ha a)$, and $\Ga_\bu\subset \ha\bA$ an essential circle in $\Ga^-$.
Assume that $\ze$ is a right (resp. left) splitting arc based
on $\Ga$ such that    $\til\ze$ is contained in the
region $\Ga_\bu^+$ above $\Ga_\bu$, and let $D(\ze_{\mid[0,s_*]})$ be a triangular domain associated to $\ze$ as in Definition \ref{def:domain}.
Assume that $\eta$ be a negatively (resp. positively) tilted arc with $\eta(0)\in \Ga_\bu$,  $\eta(]0,1])\subset \Ga_\bu^+ \cap \Ga^-$,
and $\eta(1)\in D(\ze_{\mid[0,s_*]})$.

Then \[\eta(]0,1[) \cap \ze(]0,1])\neq \emptyset.\]
\end{lemma}

In the previous lemma, the circle $\Ga_\bu$ is not assumed to be invariant by the return map.

%%%%%%%%%%%%%%%%%%%%%%%%%%%%%%%%%%%%%%%%%%%%%%%%%%%%%%%%%%%%%%%%%%%%
%%%%%%%%%%%%%%%%%%%%%%%%%%%%%%%%%%%%%%%%%%%%%%%%%%%%%%%%%%%%%%%%%%%%

\subsection{Good cylinders.}\label{Sec:goodcyl}
 We now define the main first notion of this section.

\begin{Def}\label{def:goodcyl} Consider a tame cylinder $\jC$ equipped with a \twist section $(\Sig,\ph)$,
with continuations $\hc$, $(\ha\Sig,\ha\ph)$,
and homoclinic correspondence $\psi$.
We say that $\jC$ is a good cylinder when for every essential invariant circle $\Ga \in \Ess(\ha\ph)$ in $\inte\ha\Sig$
 there exists a splitting arc based on $\Ga$.
\end{Def}

An additional quantitative definition will be necessary to produce $\de$-admissible orbits.
We denote by $\cl(E)$ the closure of  a subset $E\subset \ha\Sig$.

\begin{Def}\label{def:bounded}
With the same assumptions as in the previous definition, fix $\de>0$.
We say that $\psi=(\psi_i)_{i\in I}$ is {\em $\de$-bounded} when for each essential invariant circle $\Ga\in\Ess(\ph)$
\beq\label{eqn:delta-bounded}
\Sup\big\{\dist(x,\Ga)\mid x\in\Ga^-\cap\psi\inv\big(\Ga^+\big)\big\}<\de,
\eeq
where $\dist$ is the point-set distance on $\ha\Sig$ induced by the canonical distance on $\A^3$
and where $\psi\inv(\Ga^+)$ is the set of all points $x\in\ha\Sig$ such that there exists $i\in I$ with
$x\in\Dom\psi_i$ and $\psi_i(x)\in\Ga^+$.
We say that $\jC$ is $\de$-good when
$\psi$ is $\de$-bounded.
\end{Def}

The following Lemma proves that assuming that a good cylinder admits a $\de$-bounded correspondence is not restrictive.

%\begin{lemma}\label{rem:restriction}
%Consider a good cylinder $\jC$ equipped with a \twist section $(\Sig,\ph)$,
%with continuations $\hc$, $(\ha\Sig,\ha\ph)$,
%and homoclinic correspondence $\psi=(\psi_i)_{i\in I}$ on $\Sig$. Then given $\rho>0$, there is a homoclinic
%correspondence $\til\psi=(\til \psi_i)_{i\in I}$ such that
%$$
%\diam \Dom \til\psi_i<\rho,\qquad \forall i\in  I,
%$$
%and such that $\jC$ is still a good cylinder when endowed with $\til \psi$.
%\end{lemma}

\begin{lemma}\label{lem:deltagood}
Consider a good cylinder $\jC$ equipped with a \twist section $(\Sig,\ph)$,
with continuations $\hc$, $(\ha\Sig,\ha\ph)$,
and homoclinic correspondence $\psi=(\psi_i)_{i\in I}$ on $\Sig$. Then given $\rho>0$, there is a homoclinic
correspondence $\til\psi=(\til \psi_i)_{i\in J}$ such that
%$$
%\diam \Dom \til\psi_i<\rho,\qquad \forall i\in  I,
%$$
%and such that
$\jC$ is still a good cylinder when endowed with $\til \psi$.
\end{lemma}

\begin{proof} For each essential invariant circle $\Ga$, select one homoclinic map $\psi:=\psi_{i(\Ga)}$ which produces a splitting arc. Then one can consider the intersection $D$ of the domain $\Dom \psi$ with the $\de$-ball centered at the base point of the arc. As a consequence, the diameter of $\psi\inv(\Ga^+)$ is less than $\de$.
Let $\til\psi_i=\psi_{\vert D}$
Doing this process for each essential circle gives rise to a new homoclinic correspondence $(\til\psi_j)_{j\in J}$, with a new set of indices $J\subset I$, such that $\jC$ equipped with $(\til\psi_j)_{j\in J}$ is still a good cylinder.
\end{proof}

\vskip2mm

\noindent{\bf Convention.} {\em In the following, given a good cylinder $\jC$, we implicitly choose once and for all a section
$(\Sig,\bA,\chi,\ph)$, a continuation $\hc$, a continuation $(\ha\Sig,\ha\bA,\ha\chi,\ha\ph)$ and a homoclinic correspondence
$\psi$, which will always be denoted this way.}

%%%%%%%%%%%%%%%%%%%%%%%%%%%%%%%%%%%%%%%%%%%%%%%%%%%%%%%%%%%%%%%%%%%%
%%%%%%%%%%%%%%%%%%%%%%%%%%%%%%%%%%%%%%%%%%%%%%%%%%%%%%%%%%%%%%%%%%%%

\subsection{Good chains} Recall we want to be able to prove the existence of drifting orbits along chains of cylinders,
and not only along a single one. Again, this requires slightly more involved definitions.
We consider a $C^2$ Hamiltonian function $H$ on $\A^3$, and fix an energy $\e$.

\begin{Def}  Let $\jC_1$ and $\jC_2$ be good cylinders at energy $\e$ for $H$, with continuations
$\hc_i$ and characteristic projections $\Pi_i^\pm:W^\pm(\hc_i)\to\hc_i$.
We define
the {\em transverse heteroclinic intersection of $\hc_1$ and $\hc_2$} as the set
\beq
\Hett(\hc_1,\hc_2)\subset W^-(\inte\hc_1)\cap W^+(\inte\hc_2)
\eeq
formed by the points $\xi$ satisfying
\beq
W^-\big(\Pi_1^-(\xi)\big)\trans_\xi W^+(\inte\jC_2)\quad\textrm{and}\quad W^+\big(\Pi_2^+(\xi)\big)\trans_\xi W^-(\inte\jC_1).
\eeq
\end{Def}

As in Lemma \ref{lemma:scatteringmap}, when the characteristic projections are $C^1$,
for $\xi\in \Hett(\hc)$ we set $x^-=\Pi_1^-(\xi)\in\hc_1$ and
$x^+=\Pi_2^+(\xi)\in\hc_2$, so $W^+(\hc_1)$
and $W^-(\hc_2)$ intersect transversely
at $\xi$ in $H\inv(\e)$. Then there exist  a $3$-dimensional open neighborhood $\jO$ of
$\xi$ in $W^+(\inte \hc_1)\cap W^-(\inte\hc_2)$,
and open neighborhoods $O^-_1$ of $x^-$ in $\inte \hc_1$ and $O^+_2$ of
$x^+$ in $\inte \hc_2$, such that the restrictions $(\Pi_1^-)_{\vert \jO}$ and $(\Pi_2^+)_{\vert \jO}$
are $C^1$ {\em measure preserving} diffeomorphisms from $\jO$ onto their images $O^-_1,O^+_2$, respectively.
Hence one can define the  local diffeomorphism
$$
S_{1\to2}=\Pi^+_2\circ\big((\Pi_1^-)_{\vert \jO}\big)^{-1}:O^-_1\to O^+_2.
$$
This motivates the following definition for the notion of heteroclinic correspondences\footnote{It turns out
that in the subsequent constructions of \cite{Mar2} the projections are locally $C^1$ in a neighborhood of the
heteroclinic points we are interested in.}.

\begin{Def}
A {\em heteroclinic map} associated with $(\hc_1,\hc_2)$ is a $C^1$ diffeomorphism
\beq
\psi_{1\to2}:\Dom \psi_{1\to2} \to \Im\psi_{1\to2}
\eeq
where $\Dom \psi_{1\to2}$ is open in $\inte\ha\Sig_1$ and $\Im\psi_{1\to2}$ is open in $\inte\ha\Sig_2$,
for which there exists a $C^1$ measure preserving diffeomorphism
\beq
S_{1\to2}:\Dom S_{1\to2} \to \Im S_{1\to2}
\eeq
where $\Dom S_{1\to2}$ and $\Im S_{1\to2}$ are open in $\inte\hc_1$ and $\inte\hc_2$ respectively,
which satisfies the following conditions:
\begin{itemize}
\item there is an open subset $\Domt S_{1\to2}\subset \Dom S_{1\to2}$, with full measure in $\Dom S_{1\to2}$,
such that
\beq\label{eq:compahetero}
\forall y\in \Domt\, S_{1\to2},\qquad  W^-(y)\cap W^+\big(S_{1\to2}(y)\big)\cap\Hett(\hc_1,\hc_2)\neq\emptyset;
\eeq
\item the range of $S^1_2$ intersects the interior of $\ha\Sig_2$, i.e.,
\beq
\Im S_{1\to2}\cap \inte\ha\Sig_2\neq\emptyset;
\eeq
\item  there exists a non-negative $C^1$ function $\tau:\Dom\psi_{1\to2}\to \R$ such that
\beq
\forall x\in\Dom\psi_{1\to2},\qquad \Phi_H^{\tau(x)}(x)\in \Dom S_{1\to2} \quad \textit{and}\quad
\psi_{1\to2}(x)=S_{1\to2}\Big(\Phi_H^{\tau(x)}(x)\Big).
\eeq

\end{itemize}
\end{Def}

For purely technical reasons, we finally have to consider the case where the cylinders
$\jC_i$ are subsets of the same cylinder and they share a common boundary. The motivation, which appears in \cite{Mar1}, is to avoid the construction of sophisticated normal forms along ``long cylinders'' and make much easier the use of symplectic methods to detect homoclinic intersections (see also \cite{Mar4}).

\vskip2mm

%have one boundary component in common(and their union itself is a cylinder).
In this case we introduce a {\em transition map} between the sections,  defined as the natural flow-induced map.

\begin{Def}\label{def:adjacent}
Assume that the good cylinders $\jC_1$ and $\jC_2$ are contained in a cylinder $\jC$ and satisfy
$\d^\bu\jC_1=\d_\bu\jC_2$, in which case we say they are {\em adjacent}.
A {\em transition map} from $\hc_1$ to $\hc_2$ is a homoclinic map (relative to $\jC$)
\beq
\psi_{1\to2}:\Dom \psi_{1\to2} \to \Im\psi_{1\to2}
\eeq
where $\Dom \psi_{1\to2}$ is an open neighborhood of the circle $\Ga=\d^\bu\Sig_1=\d_\bu\Sig_2$ in $\ha\Sig_1$, and
$\Im \psi_{1\to2}$ is an open neighborhood of $\Ga$ in $\ha\Sig_2$. More precisely,
there exists a non-negative $C^1$ function $\tau:\Dom \psi_{1\to2}\to\R$ such that
$$
\psi_{1\to2}(x)=\Phi_H^{\tau(x)}(x),\qquad \forall x\in \Dom \psi_{1\to2}.
$$
\end{Def}

In the following our heteroclinic maps and transition maps will play the same role in our constructions and
we will not make any distinction between them (including terminology).
This yields our final definition.

\begin{Def}\label {def:goodchains} Fix $\de>0$. A {\em $\de$-good chain of cylinders} at energy $\e$
is a finite ordered family $(\jC_k)_{1\leq k\leq k_*}$ of $\de$-good cylinders at energy $\e$,
such
%that for $1\leq k\leq k_*-1$, either $\jC_k$ and $\jC_{k+1}$ are adjacent, or
there exists a heteroclinic map $\psi_{k\to{k+1}}$
from $\hc_k$ to $\hc_{k+1}$ which satisfies the following condition:
\begin{itemize}
\item for any open neighborhood $O$ of $\d^\bu\Sig_k$ in $\ha \Sig_k$, the
image $\psi_{k\to{k+1}}(O)$
intersects a dynamically minimal essential invariant circle in $\Ess(\ph_{k+1})$ located in a $\de$-neighborhood of
$\d_\bu\Sig_{k+1}$ in $\ha\Sig_{k+1}$.
\end{itemize}
\end{Def}

Again, the $\de$-neighborhoods are defined relatively to the induced distance on the section $\ha\Sig_{k+1}$.
%Note that, in the case where $\jC_k$ and $\jC_{k+1}$ are adjacent, the transition map $\psi_{k\to{k+1}}$
%can be chosen to send  $\d^\bu\Sig_k$ onto $\d_\bu\Sig_{k+1}$, so that the previous property is also satisfied by~$\psi_{k\to{k+1}}$.

%%%%%%%%%%%%%%%%%%%%%%%%%%%%%%%%%%%%%%%%%%%%%%%%%%%%%%%%%%%%%%%%%%%%
%%%%%%%%%%%%%%%%%%%%%%%%%%%%%%%%%%%%%%%%%%%%%%%%%%%%%%%%%%%%%%%%%%%%

\subsection{Polysystems}
\label{sec:polysystems}
Let us first make the definition of an orbit of a polysystem more precise (see \cite{M02,Mar08} for a formal definition).
Let $A$ be some set and consider a set $f=\{f_i\mid i\in I\}$ of locally defined maps  $f_i:\Dom f_i\subseteq A \to A$.
We say that a finite sequence $(x_n)_{0\leq n\leq n_*}$, $n_*\geq1$, of points of $A$ is a {\em finite orbit of $f$, of length
$n_*+1$,}
provided that there exists
$\om=(i_n)_{0\leq n\leq n_*-1}\in I^{n_*}$ such that for $0\leq n\leq n_*-1$:
\[
x_{n+1}= f_{i_n}(x_n), \textrm{ with } x_n\in \Dom f_{i_n} \textrm { and } x_{n+1}\in \Dom f_{i_{n+1}}.
\]
We use the following short-hand notation
\[
x_{n_*}=f^\om(x_0).
\]
We consider the point $x_0$ as being the length-$1$ orbit of $x_0$.

Given a subset $B\subset A$, we set
$$
f^\om(B)=\bigcup_{x\in B_\om}f^\om(x)
$$
where $B_\om$ is the subset of $B$ formed by the points $x$ such that $f^\om(x)$ is well-defined.

The {\em full orbit of  $B\subset A$} under $f$  is the subset of $A$ formed by the union of all
$f^\om(B)$ for all sequences (of any length) $\om$.
In particular $B$ is contained in its full orbit under~$f$.

\vskip2mm

\noindent{\bf Conventions.} {\em
Given locally defined maps $f_i:\Dom f_i\to \Im f_i$ on $A$, we write $f_i\circ f_j$ for the map defined by composition on the
subset $\Dom f_i\cap \Im f_j$ and for a subset $B\subset A$, we write $f_i(B)$ for $f_i\big(\Dom f_i\cap B\big)$.

Given a finite set of polysystems $f,g,\ldots$ on $A$, we write $\{f,g,\ldots\}$ or $(f,g)$
for the polysystem formed by their union.

Given a polysystem $f=\{f_i\mid i\in I\}$ on a set $A$, and a subset $A_*$ of $A$, we define the
restriction $f_{\vert A_*}$ of $f$ to $A_*$ as the polysystem formed by the maps
$$
f_i^*:\Dom f_i\cap A_*\cap f_i\inv(A_*)\to A_*.
$$

Note that $f_i^*$ is not necessarily open, even when $f_i$ is.}

\begin{rem}
By Lemma~\ref{lem:deltagood}, any pseudo-orbit
of $(\ph,\til\psi)$ is a pseudo-orbit of  $(\ph,\psi)$. As a consequence, given $\de>0$,
one can always assume the homoclinic correspondence of a good cylinder to be $\de$-bounded.
\end{rem}

\vskip2mm

\begin{Def}\label{def:pseudoorbit}
Let $H$ be a $C^2$ Hamiltonian on $\A^3$ and fix  an energy~$\e$.
Let $(\jC_k)_{1\leq k\leq k_*}$ be a good chain of cylinders in $H\inv(\e)$.
Its {\em associated polysystem} is the following set of (locally) defined diffeomorphisms:
\beq\label{def:polysystem}
\mathscr{F}=\big\{\ha\ph_k\mid 1\leq k\leq k_*\big\}\cup \big\{\psi_{k,i}\mid 1\leq k\leq k_*,\ i\in I_k\big\}\cup\big\{\psi_{k\to k+1}\mid 1\leq k\leq k_*-1\big\},
\eeq
where $\ha\ph_k$ is the continuation of $\ph_k$, $\psi_k=(\psi_{k,i})_{i\in I_k}$ is the homoclinic correspondence associated to
$\jC_k$, and $\psi_{k\to k+1}$ is either a heteroclinic map or a transition map from $\hc_k$ to $\hc_{k+1}$.

This polysystem will be indifferently considered to be defined on
$
\bigcup_{1\leq k\leq k_*}\ha\bA_k,
$
or on
$
\bigcup_{1\leq k\leq k_*}\ha\Sig_k
$.

The finite orbits of the polysystem will be called {\em pseudo-orbits along the chain}.
\end{Def}

The pseudo-orbits above are therefore generated by successive applications of the maps
$\ha\ph_k$, $\psi_{k,i}$, $\psi_{k\to k+1}$ in any possible order.

%%%%%%%%%%%%%%%%%%%%%%%%%%%%%%%%%%%%%%%%%%%%
%%%%%%%%%%%%%%%%%%%Section3%%%%%%%%%%%%%%%%%%%%%
%%%%%%%%%%%%%%%%%%%%%%%%%%%%%%%%%%%%%%%%%%%%

\section{Pseudo orbits along a good chain of cylinders}\label{Sec:pseudoorb}
\setcounter{paraga}{0}

Our aim in this section is to prove the following result.

\begin{thm}\label{thm:pseudoorb1}
Let $H$ be a $C^2$ Hamiltonian on $\A^3$ and fix an energy $\e$.
Let $(\jC_k)_{1\leq k\leq k_*}$ be a $\de$-good chain of cylinders at energy $\e$ for some $\de>0$.

Then there exists a pseudo-orbit $(x_n)_{0\leq n\leq n_*}$ of the polysystem $\mathscr{F}$ along the chain such that for
any essential invariant circle $\Ga$ in $\bigcup_{1\leq k\leq k_*}\Ess(\ph_k)$, there
is a $\nu\in\{0,\ldots,n_*\}$ with $d(x_\nu,\Ga)<\de$.
\end{thm}

We first prove in Section~\ref{sec:pseudoorbcyl} the existence of pseudo-orbits in the case where the chain is reduced to a single
cylinder. Then {\bf Theorem~\ref{thm:pseudoorb1}} is easily deduced from this preliminary result in Section~\ref{sec:pseudoorbchain}.
The main remark in our proof comes from Lemma~\ref{lem:dense} of Appendix~\ref{app:symmetrization}: it is enough
to prove the existence of pseudo-orbits for the ``symmetrized'' polysystem
\beq\label{def:sympolysystem_sym}
\mathscr{G}=\big\{\ha\ph_k,\ \ha\ph_k\inv \mid 1\leq k\leq k_*\big\}\cup \big\{\psi_{k,i}\mid 1\leq k\leq k_*,\ i\in I_k\big\}
\cup\big\{\psi_{k\to k+1}\mid 1\leq k\leq k_*-1\big\}
\eeq
and obtain the conclusion by the density of the pseudo-orbits of $\mathscr{F}$ relative to the set of the  pseudo-orbits of $\mathscr{G}$.

%%%%%%%%%%%%%%%%%%%%%%%%%%%%%%%%%%%%%%%%%%%%%%%%%%%%%%%%%%%%%%%%%%%%
%%%%%%%%%%%%%%%%%%%%%%%%%%%%%%%%%%%%%%%%%%%%%%%%%%%%%%%%%%%%%%%%%%%%

\subsection{Pseudo-orbits in the case of a single cylinder}\label{sec:pseudoorbcyl}

\def\len{{\rm length\,}}
%We denote by  $\CC(A,B)$ the connected component of  $A$ containing $B$.
Given $\nu>0$, we define a {\em $\nu$-ball}
of $\T\times\R$ as a subset $B=B_\th\times B_r$ where $B_\th$ and $B_r$ are intervals of $\T$ and~$\R$
respectively, such that
\beq\label{eq:cball}
\len B_r>\nu\,\len B_\th.
\eeq
The {\em center of $B$} is $(a_\th,a_r)$, where $a_\th$, $a_r$ are the mid-points of $B_\th$ and $B_r$, respectively.
We keep the notation and conventions of the last section for the good cylinders and their continuations.

\begin{prop}\label{prop:tamecyl}
Let $H$ be a $C^2$ Hamiltonian on $\A^3$ and fix an energy $\e$.
Let $\jC$ be a good cylinder at energy $\e$ for $H$.
Let $g=(\ha\ph,\ha\ph\inv,\psi)$ be the associated symmetrized) polysystem on $\ha\bA$,
with $\psi=(\psi_i)_{i\in I}$.
Fix $\Ga_\bu\in\Ess(\ha\ph)$ with $\Ga_\bu\subset \bA\setm\Ga(b)$.
Fix a neighborhood $V$ of $\Ga_\bu$ in $\bA$ (see Figure \ref{fig:annulus}).

Then the full orbit of $V$ under $g$ contains $\Ga(b)$.
\end{prop}

\begin{figure}[h]
\centering
\includegraphics[width=0.3\textwidth]{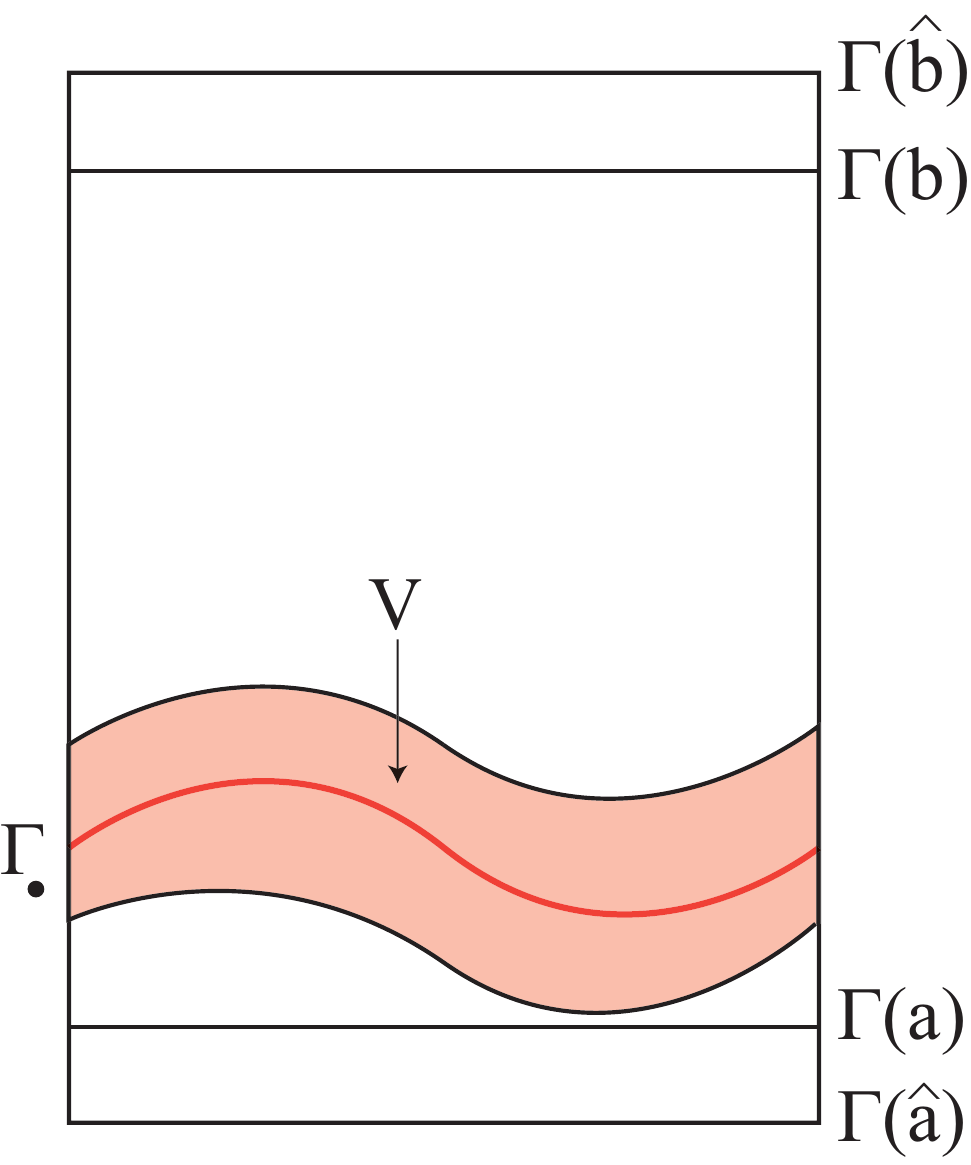}
\caption{The setting of Proposition \ref{prop:tamecyl}}\label{fig:annulus}
\end{figure}

The assumption $\Ga_\bu\subset \bA\setm\Ga(b)$ is a matter of simplification
and can easily be relaxed in practical cases.

\begin{proof} [Proof of Proposition~\ref{prop:tamecyl}]
We assume for example that $\ha\ph$ tilts the vertical to the right, the other case being exactly similar.
%We will prove a slightly stronger property than our initial statement. Namely,
Let $\bA_*$ be the
$\ha\ph$--invariant subannulus of~$\ha\bA$
limited by the (disjoint) circles $\Ga_\bu$ and $\Ga(\ha b)$.

Without loss of generality, one can assume
that for every circle $\Ga\in\Ess(\ha\ph)$ contained in $\bA_*\setm\{\Ga_\bu\}$, for every map $\psi_i$
such that $\Im\psi_i\cap\Ga\neq\emptyset$, we have that $\Dom \psi_i\subset \bA_*$ and
$\Im\psi_i\subset\bA_*$, therefore the  good cylinder condition is still satisfied on $\bA_*$.
\vskip1mm

\paraga We consider the restricted polysystem
$g_{\vert \bA_*}$ instead of $g$ and we consider a neighborhood $V$ of $\Ga_\bu$ {\em in $\bA_*$}.
For simplicity, we still denote by $\ha\ph$ and $\psi=(\psi_i)_{i\in I}$ the restrictions of $\ha\ph$ and $\psi$ to
$\bA_*$.
{Thanks to our previous remark, the maps $\psi_i$ are {\em open}.}

We also assume without loss of generality that the set $V$ is open in $\bA_*$  and connected.
Let $U$ be the full orbit of  $V$ under the polysystem
$g=\big(\ha\ph,\ha\ph\inv,\psi=(\psi_i)_{i\in I}\big)$ on $\bA_*$.
Note that $\ha\ph(U)= U$ and
%$\ha\ph\big(O(V)\big)\subset O(V)$,
$\psi_i(U)\subset U$.
%and $\psi_i\big(O(V)\big)\subset O(V)$ for $i\in I$, by construction.
Set $U_c=\CC(U,\Ga_\bu)$.
Then $U_c$ is open and contains $V$, so  $\ha\ph(U_c)= U_c$.
%Similarly,
%$O_c(V)$ is open and contains $V$, so
%$\ha\ph\big(O_c(V)\big)\subset O_c(V)$.
It is therefore enough to prove that $U_c$ intersects $\Ga(b)$: since $\Ga(b)$ is dynamically minimal,
this yields $\Ga(b)\subset U_c\subset U$.

\vskip2mm

Let us assume by contradiction that $U_c\cap\Ga(b)=\emptyset$, so that $U_c$ is contained in the
lower connected component of $\bA_*\setm \Ga(b)$. We want to apply the Birkhoff graph theorem to
describe the frontier of $U_c$ (see Appendix\ref{App:twistmaps}), however $U_c$ is not homeomorphic to an annulus.
The following paragraph describes a standard way to ``fill the holes'' of $U_c$ without altering its frontier,
making it possible to apply the Birkhoff theory.

\paraga Set $O=\bA_*\setm \ov{U_c}$, so that $O$ is open, contains $\Ga(\ha b)$, and $O\cap V=\varnothing$. Moreover,
since  $\ha\ph(\ov{U_c})= \ov{U_c}$,
$$
\ha\ph(O)=\bA_*\setm \ha\ph(\ov{U_c})=\bA_*\setm \ov{U_c}=O.
$$
Then $\ha\ph\big(\CC(O,\Ga(\ha b)\big)= \CC(O,\Ga(\ha b))$ and so $\ha\ph(\ov{\CC(O,\Ga(\ha b)})= \ov{\CC(O,\Ga(\ha b))}$.
Let
$$
\ha U=\bA_*\setm \ov{\CC(O,\Ga(\ha b))},
$$
 so that $\ha U$ is open and $\ha\ph(\ha U)= \ha U$, and set finally
\beq\label{eq:defbU}
\bU=\CC(\ha U,\Ga_\bu),
\eeq
hence $\bU$ is open, connected and $\ha\ph(\bU)=\bU$. Moreover clearly
\beq\label{eq:inclusion}
\ov\bU\subset \bA_*\setm \CC(O,\Ga(\ha b)),
\eeq
and
\beq\label{eq:U}
U_c\subset\bU,
\eeq
since $\ov O=\bA_*\setm\Int(\ov{U_c})\subset \bA_*\setm U_c$, so
$\ov{\CC(O,\Ga(\ha b))}\subset \bA_*\setm U_c$ and $U_c\subset \bA_*\setm \ov{\CC(O,\Ga(\ha b))}=\ha U$,
which proves (\ref{eq:U}) since $\Ga_\bu\subset U_c$.

\paraga We denote by $\Fr A$ the frontier $\ov A\setm \Int A$ of a subset in a topological space. Let $\Ga:=\Fr\bU\cap \Ga_\bu^+$.
We shall prove that $\Ga$ is a Lipschitz graph over $\T$, invariant under $\ha\ph$, by the
Birkhoff theorem (see Appendix~\ref{App:twistmaps}).
By local connectedness of $\bA_*$, one readily proves that $\Int\ov{\bU}=\bU$, since $\bU$ is a connected component
of the complement of the closure of an open set.
Moreover $\ha\ph(\bU)=\bU$.
Let now $\bf S$ be the quotient of $\bA_*$ by the identification of each boundary circle to one point, so that $\bf S$ is homeomorphic
to $S^2$. Up to this quotient,  $\bU$ is a connected component of the complement in $\bf S$ of a compact connected subset, so
is homeomorphic to a disk. Going back to the initial space $\bA_*$ proves
that $\bU$ is homeomorphic to
$\T\times [0,1[$. So by the Birkhoff theorem, $\Ga=\d\bU$ is a Lipschitz graph over $\T$, invariant under $\ha\ph$
(see \cite{M02} for more details).

\paraga Let us now prove that $\Ga\subset \ov{U_c}$, and so $\Ga\subset \Fr(U_c)=\cl(U_c)\setm U_c$.
Assume that $x\in\Ga$ is not in $\ov{U_c}$, so that there exists
a small ball $B(x,\eps)$ with $B(x,\eps)\cap \ov U_c=\varnothing$. Let $z$ be some point on the vertical through $x$,
located under $\Ga$ and inside $B(x,\eps)$.  Let us show that the semi-vertical $\sig=V^+(z)$ over $z$ in $\bA_*$ is disjoint from $\ov{U_c}$.
First $\Ga\cap\sig=\{x\}$, since $\Ga$ is a graph, so that $\sig=[z,x]\cup[x,\xi]$, with $\xi\in \Ga(\ha b)$.
Clearly $[z,x]\subset B(x,\eps)$ so $[z,x]\cap \ov{U_c}=\varnothing$, and $]x,\xi]\cap\ov{\bU}=\varnothing$ since
$\Ga=\d\bU$ is a graph. Since $U_c\subset\bU$, this proves that $\sig\cap\ov U_c=\vide$.

As a consequence $\sig\cup  \Ga(\ha b)$ is a connected set which satisfies
$
(\sig\cup \Ga(\ha b))\cap \ov{U_c}=\varnothing.
$
Therefore $(\sig\cup \Ga(\ha b))\subset \CC(O,\Ga(\ha b))$ and thus $(\sig\cup \Ga(\ha b))\cap\ov{\bU}=\varnothing$
by (\ref{eq:inclusion}).
This is a contradiction since
$x\in\Ga\subset\ov{\bU}$. Therefore $\Ga\subset\ov{U_c}$.

\paraga Since $\Ga$ is an invariant essential circle for the good \twist map $\ha\ph$, there are only two possibilities:
\vskip1mm
\noindent -- either $\Ga$ is the upper boundary of a Birkhoff zone,
\vskip1mm
\noindent -- {or} $\Ga$ is accumulated from below by essential invariant circles in the Hausdorff topology.
\vskip1mm
\noindent We will prove that both possibilities yield a contradiction with the initial assumption that
$U_c\cap\Ga(b)=\emptyset$.

We call attention that the argument from this point on is different from \cite{M02}, who considers a polysystem consisting of two globally defined mappings on the annulus. In our case the polysystem contains not only  the  maps $\ha \ph,\ha\ph^{-1}$ which are globally defined,   but also the maps $\psi_i$, which, by the definition of a good cylinder,  are only  locally defined around points on essential invariant circles.

\paraga Assume first that $\Ga$ is the upper boundary of a Birkhoff zone $\jZ$ and let $\Ga_*$ be the lower
boundary of $\jZ$. Let $\nu$ be the Lipschitz constant of $\Ga_*$.
Since $\Ga_*$ is a graph and $U_c$ is open, connected, contains $V$ and satisfies $\ov{U_c}\cap\Ga\neq\varnothing$,
then $U_c\cap \Ga_*\neq\varnothing$. So there exists a $\nu$-ball $\bB\subset U_c$ centered on $\Ga_*$.

Since $\jC$ is a good cylinder,
there exists a  splitting arc $\ze$ based at some point $\al_0$ of $\Ga$
(see Section~\ref{Sec:setting}). Either $\ze$ is a non-vertical arc, in which case it is a right splitting arc or a left splitting arc, or it is a vertical arc.

If $\ze$ is a right (resp. left) splitting arc,
there exists a triangular domain $D:=D(\ze_{\mid [0,s_*]})$
associated to $\ze$ such that  $D\subset \jZ$.

First, assume  that $\til\ze$ is a right splitting arc.
By Proposition~\ref{thm:extensionbirkhoff2}, there exist  $z_0\in \bB$ and
$n\in\N$ such that $z_n:=\ha\ph^n(z_0)\in D$.
By Lemma~\ref{lem:negtilt} there exists a {\em negatively tilted arc} $\eta$ with image
in $\bB$ emanating from  $\Ga_*$ and ending at $z_0$.
Therefore, by Lemma~\ref{lem:iteratetilt}, $\eta_n:=\ha\ph^n\circ \eta$ is a {\em negatively} tilted arc
emanating from $\Ga_*$ and ending at $z_n\in D$.
By Lemma \ref{lem:tiltarc}, the image $\til\eta_n$ intersects $\til\ze$.

Second, assume that $\til\ze$ is a left splitting arc. We use $\ha\ph\inv$ instead of $\ha\ph$. This yields a point $z_0\in \bB$
such that $z_{-n}:=\ha\ph^{-n}(z_0)\in D$, and a
{\em positively  tilted arc} $\eta$  emanating from $\Ga_*$ and ending at $z_{0}$.
Then $\eta_{n}:=\ha\ph^{-n}\circ \eta$ is a  {\em positively  tilted arc} emanating from $\Ga_*$ and ending at $z_{-n}\in D$, whose image $\til\eta_{n}$ intersects $\til\ze$.

Third, assume that $\til\ze$ is a vertical splitting arc. Then $\ha\ph^{-1}(\til\ze)$ is a right arc (not necessarily splitting), see Remark \ref{rem:nonvert}.
Let $D$ be a triangular domain associated to  $\ha\ph^{-1}(\til\ze)$. As above, there  exists  $z_0\in \bB$, a {\em negatively tilted arc} $\eta$ with image in $\bB$, and  $m\in\N$ such that the negatively tilted arc   $\ha\ph^m(\til\eta)$, which emanates from $\Ga_*$ and ends at $z_m$, intersects  $\ha\ph^{-1}(\til\ze)$. Since $\ha\ph^m(\til\eta)\cap\ha\ph^{-1}(\til\ze)\neq\varnothing$, we obtain
$\ha\ph^{m+1}(\til\eta)\cap\til\ze\neq\varnothing$. Setting  $n=m+1$ and $\til\eta_n=\ha\ph^{n}(\til\eta)$ we obtain $\til\eta_n\cap\til\ze\neq \varnothing$.
%Note that while $\ha\ph^{-1}(\til\ze)$ is not necessarily a splitting arc,
%ct$\til\ze$ is a splitting arc.

As a consequence,  in either case we have that $U_c\cap \til\ze\neq\varnothing$ since $\til\eta_n\subset U_c$, and therefore there is a small open ball
$B\subset U_c$ centered on $\ze(]0,1])$ and %, by definition,
an index $i\in I$ such that  $B\subset\Dom \psi_i$.
Thus $\psi_i(B)$ is an open set which intersects $\Ga$, and therefore also
$U_c$ since $\Ga\subset \ov{U_c}$. This proves that $\psi_i(B)\subset U_c$ by connectedness,
so that $U_c$ contains points strictly above the circle $\Ga$.
This is a contradiction with the construction of
$\Ga=\Fr\bU$ and the inclusion $U_c\subset\bU$, which ensures that all points of $U_c$ are located below $\Ga$.

\paraga  Assume now that $\Ga$ is accumulated from below by an increasing sequence $(\Ga_m)_{m\geq1}$ of essential
invariant circles for $\ha\ph$. Let $\ze$ be a splitting arc based on $\Ga$.
Let $S_m$ be the closed strip limited by $\Ga_m$ and $\Ga_{m+1}$.
For $m$ large enough,  $S_m\cap \til\ze$ contains a $C^0$ curve $\ell$ which intersects both $\Ga_m$ and $\Ga_{m+1}$.
Now $\Ga\subset \ov{U_c}$, so that $U_c\cap S_m$ contains a $C^0$ curve $\ell'$ which  also
intersects both $\Ga_m$ and $\Ga_{m+1}$. Therefore, by Lemma~\ref{lem:torsion}, there exists an integer $n$ such that
$\ha\ph^n(U_c)\cap \ell\neq\vide$, { and so  by invariance of $U_c$ under $\ha\ph$,  $U_c\cap \ell\neq\vide$.
Since $\ell\subset \til\ze\subset \Dom \psi_i$
for some~$i\in I$, there exists a ball $B\subset U_c$ centered on $\ell\subset \til\ze$ and contained in $\Dom\psi_i$.
This yields the same contradiction as in the previous paragraph.}

\paraga We have therefore proved that the two possibilities of our alternative yield a contradiction,
which proves that our initial
assumption $U_c\cap \Ga(b)=\emptyset$ is false.  This concludes the proof of the proposition.
\end{proof}

As a consequence of Lemma~\ref{lem:dense}, under the assumptions of Proposition~\ref{prop:tamecyl},
for any $\de>0$, the polysystem $f=(\ha\ph,\psi)$ admits orbits which intersect the $\de$-neighborhoods
of $\Ga_\bu$ and $\Ga(b)$. We now take into account the additional assumption that $\psi$ is $\de$-bounded.

\begin{cor}\label{cor:deltaad}
With the same assumptions and notation as in {\rm Proposition~\ref{prop:tamecyl}}, if $\psi$ is $\de$-bounded,
any orbit of $f=(\ha\ph,\psi)$ which intersects the $\de$--neighborhoods of $\Ga_\bu$ and $\Ga(b)$ intersects
the $\de$-neighborhood  of every element of $\Ess(\ph)$ contained in $\bA$.
\end{cor}

\begin{proof} Consider an orbit $(x_n)_{0\leq n\leq n_*}$ of $f$
with $x_0$ in the $\de$-neighborhood of $\Ga_\bu$ and $x_{n_*}$ in the $\de$-neighborhood of $\Ga(b)$.
Let $\Ga\in\Ess(\ph)$ be contained in the complement of these neighborhoods.
{Then there exists $\nu\in\{0,\ldots,n_*\}$ such that $x_\nu\in \Ga^-\cap \psi\inv(\Ga^+)$
(if the orbit has no point in $\Ga^-\cap \psi\inv(\Ga^+)$, an immediate
induction shows that $x_n\in\Ga^-$ for $0\leq n\leq n_*$).
Now, since $\psi$ is $\de$--bounded,
$$
\Sup\big\{\dist(x,\Ga)\mid  x\in\Ga^-\cap \psi\inv(\Ga^+)\big\}<\de,
$$
which proves that $\dist(x_\nu,\Ga)<\de$.}
\end{proof}

%%%%%%%%%%%%%%%%%%%%%%%%%%%%%%%%%%%%%%%%%%%%%%%%%%%%%%%%%%%%%%%%%%%%
%%%%%%%%%%%%%%%%%%%%%%%%%%%%%%%%%%%%%%%%%%%%%%%%%%%%%%%%%%%%%%%%%%%%
\subsection{End of proof of Theorem~\ref{thm:pseudoorb1} and a remark on pseudo-orbits}\label{sec:pseudoorbchain}

Let $\bA_k=\T\times[a_k,b_k]$ and $\ha \bA_k=\T\times[\ha a_k,\ha b_k]$. Fix $\de>0$.
Let $V_1=\T\times [a_1,a_1+\de]\subset\bA_1$. We proceed by induction.
Given $k\in\{1,\ldots,k_*\}$,
we set the following condition:
\begin{itemize}
\item $(H_k)$:  the full orbit $O:=O(V_1)$ of $V_1$ under the polysystem $\mathscr{G}$
contains $\Ga(b_k)\subset \bA_k$,
\end{itemize}
where $\mathscr{G}$ was defined in \eqref{def:sympolysystem_sym}.

\vskip1mm

By Proposition~\ref{prop:tamecyl}, the connected component of $O$ containing $\Ga(a_1)$
also contains $\Ga(b_1)$, so that $(H_1)$ is satisfied.

\vskip1mm

Assume that $(H_k)$ is satisfied for some $k\in\{1,\ldots,k_*-1\}$. Then by our definition of a $\de$-good chain,
since $O$ is a neighborhood of $\Ga(b_k)$, $\psi_{k\to{k+1}}(O)$ intersects some dynamically
minimal circle $\Ga_{k+1}\in \Ess(\ph_{k+1})$ located in the $\de$-neighborhood of $\Ga(a_{k+1})$.
As a consequence, the invariance of $O$ under $\ha\ph_{k+1}$ proves that $O$ contains a neighborhood $V_{k+1}$ of
$\Ga_{k+1}$, and Proposition~\ref{prop:tamecyl} applied to  $(\ha \ph_{k+1},\ha \ph_{k+1}\inv,\psi_{k+1})$ proves that $O$
contains a neighborhood of $\Ga(b_{k+1})$, so that $(H_{k+1})$ is satisfied.

\vskip1mm

By finite induction, $(H_k)$ is satisfied for $k\in\{1,\ldots,k_*\}$.

\vskip1mm

The previous statement holds without any boundedness assumption on the homoclinic or heteroclinic maps.
Assuming now that the full assumptions for good chains are satisfied, then
Corollary~\ref{cor:deltaad} (and the remark before) proves the existence of an orbit of the initial polysystem
$\mathscr{F}$ defined in~(\ref{def:polysystem}) which intersects the $\de$-neighborhood
of each element of $\Ess(\ph_k)$, $1\leq k\leq k_*$, which concludes the proof of {\bf Theorem~\ref{thm:pseudoorb1}}.
$\hfill\Box$

\vskip3mm

We can finally take advantage of the Poincar\'e recurrence theorem and prescribe one possible form for  the pseudo-orbits. This is the content of the next remark.

\begin{rem}\label{rem:formorbits}
 With the same assumptions as in {\bf Theorem~\ref{thm:pseudoorb1}},
one can choose the pseudo-orbit $(x_n)_{1\leq n\leq n_*}$ so that it is $2\de$-admissible and is the
concatenation of segments of the form
$$
\psi\circ\ph^m(x)
$$
with $m>0$, where
$$
\ph\in\big\{\ph_k\mid 1\leq k\leq k_*\big\}
\quad\text{and}\quad
\psi\in\big\{\psi_{k,i}\mid 1\leq k\leq k_*,\ i\in I_k\big\}\cup\big\{\psi_{k\to k+1}\mid 1\leq k\leq k_*-1\big\}.
$$
In other words, in the pseudo-orbit $(x_n)$ each application of a   homoclinic/heteroclinic maps is always interspersed with a positive number of applications of the flow-induced map on the cylinder.
\end{rem}

The proof follows the same lines as that of Lemma~\ref{lem:dense}. Assume that a pseudo-orbit involves the consecutive use of two homoclinic/heteroclinic maps $\psi$ and $\psi'$ (we get rid of the indices), so that say $x_{n+2}=\psi'\circ\psi (x_n)$. Using the Poincar\'e recurrence theorem,  one can slightly perturb the point $x_n$ into $x'_n$ in such a way that $\psi(x'_n)$ is positively recurrent, and
produce a pseudo orbit starting arbitrarily close to $x_n$ and ending arbitrarily close to $x_{n+2}$, of the following form
$$
\begin{array}{lll}
x'_{n},\  x'_{n+1}=\psi(x'_n),\\
x'_{n+2}=\ph(x'_{n+1}),\ldots, x'_{n+m+1}=\ph^m(x'_{n+1})\\
x'_{n+m+2}=\psi'(x'_{n+m+1}).
\end{array}
$$
Indeed, it suffices to choose $m$ is such a way that $x'_{n+m+1}$ is close enough to the initial image
$\psi(x'_n)$ so as to be able to control the whole orbit, using the continuity of the maps involved in the polysystem.

\subsection{An addendum for singular cylinders}\label{sec:singular}
We now state an additional result which will enable us in \cite{Mar2} to deal with ``singular cylinders''
exactly in the same way as if they were good cylinders.
More precisely,  singular cylinders are $3$-dimensional manifolds diffeomorphic to the product of $\mathbb{T}^1$
with a $2$-sphere with three open disks with disjoint closures
cut off\footnote{Such singular cylinders appear near double-resonances in near-integrable systems}.
Our way to deal with singular cylinders is to ``fill-in'' one of the cut-off disks, and extend the dynamics to a tame
$3$-dimensional cylinder. Since the  extended dynamics on the filled disk
does not correspond to ``true'' orbits of the system, one is interested in
obtaining pseudo-orbits that avoid the ``filled'' disk.
See \cite{Mar2} for details.

\begin{cor}\label{cor:singcyl}
Assume the conditions from Proposition~\ref{prop:tamecyl}. Assume that {$K\subset\T\times\,]a+\de,b-\de[$}
is a compact subset of $\bA$, invariant under $\ph$ and contained in the interior of some Birkhoff zone  of $\ph$.
%Let $\rho=\dist(K,\d Z)$.
%Assume that for each $i\in I$,
%\beq\label{eq:diamdom1}
%\diam\Dom \psi_i<\rho/2.
%\eeq
%and that $\psi$ is $\de$-bounded.
Then $f=(\ph,\psi)$ admits a pseudo-orbit which does not intersect~$K$. If moreover $\psi$ is $\de$-bounded,
then there exists a $\de$-admissible pseudo-orbit.
\end{cor}

\begin{proof} Since the interior of the Birkhoff zone contains no element of $\Ess(\ha\ph)$, one can reduce
the domains of the heteroclinic correspondence in such a way that $\Dom\psi_i\cap K=\emptyset$
for $i\in I$, still preserving the good cylinder condition. Therefore there exists a  pseudo-orbit
$(x_n)_{0\leq n\leq n_*}$ for  $(\ph,\psi)$ such that $x_0$ is $\de$-close to $\Ga(a)$ and $x_n$ is $\de$-close
to $\Ga(b)$. This orbit cannot intersect~$K$, for if $x_n\in K$, then by induction $x_m\in K$ for $m\geq n$.
The existence of a $\de$-admissible pseudo-orbit is then immediate.
\end{proof}

%Assumption~(\ref{eq:diamdom1}) is not a restriction: one easily sees that
%one can always ``reduce the domains'' of the homoclinic correspondence without altering the properties
%relative to the splitting arcs.

%%%%%%%%%%%%%%%%%%%%%%%%%%%%%%%%%%%%%%%%%%%%
%%%%%%%%%%%%%%%%%%%Section4%%%%%%%%%%%%%%%%%%%%%
%%%%%%%%%%%%%%%%%%%%%%%%%%%%%%%%%%%%%%%%%%%%

\def\vvert{\,\vert\,}

\section{Shadowing of pseudo-orbits and proof of Theorem~\ref{thm:main1}}\label{Sec:shadowing}
\setcounter{paraga}{0}

The aim of this section is to prove that given an arbitrary pseudo-orbit along a good chain
(as introduced in Definition~\ref{def:pseudoorbit}),
there exists an orbit of the Hamiltonian system \eqref{eq:hampert}, which passes arbitrarily close to each point of the  pseudo-orbit.
{\bf Theorem~\ref{thm:main1}} is then an immediate consequence of the existence of pseudo-orbits
intersecting arbitrarily small neighborhoods of any essential invariant circle in the \twist sections of a good chain,
as proved in {\bf Theorem~\ref{thm:pseudoorb1}}. The main result of this section is the following.

\begin{thm}\label{thm:shadowing}
Let $H$ be a Hamiltonian of class $C^2$ on $\A^3$ and fix an energy $\e$.
Let $(\jC_k)_{1\leq k\leq k_*}$ be a good chain of cylinders at energy $\e$.

Then given a pseudo-orbit $(x_n)_{0\leq n\leq n_*}$ along this chain and any $\de>0$, there
exist a solution of the Hamiltonian vector field $X_H$, with initial condition $a_0$,
and a finite  sequence of times
$(\tau_n)_{0\leq n\leq n_*}$, with $\tau_0=0$, such that
\[\dist\big(\Phi_H(\tau_n,a_0),x_n\big)<\de\textrm{ for }0\leq n\leq n_*.\]
\end{thm}

\begin{proof} The proof generalizes those of  \cite{BT99,DLS00}.
We keep the notation from Section~\ref{Sec:pseudoorb}
for the polysystem $\{\ha\ph_k,\psi_k,\psi_{k\to{k+1}}\}$, with $\psi_k=(\psi_{k,i})_{i\in I_k}$,
which we consider to be defined on the union
$$
\bigcup_{1\leq k\leq k_*}\ha\Sig_k.
$$

\def\bC{{\bf C}}
\paraga Let $(x_n)_{0\leq n\leq n_*}$ be a pseudo-orbit along the chain.

First, we shall prove the existence of a sequence $(y_j)_{0\leq j\leq j_*}$ of points of
$\ha\bC=\bigcup_{1\leq k\leq k_*}\hc_k$, which satisfies the following two properties.
\begin{itemize}
\item $(P_1)$ For each $n\in\{0,\ldots,n_*\}$, there exists $j_n\in\{0,\ldots,j_*\}$ such that
$
x_n=y_{j_n}.
$

\item $(P_2)$  Fix $j\in\{0,\ldots, j_*-1\}$ and  let $k$ be such that $y_j\in\hc_k$. Then
one of the following three conditions is satisfied:

 $(i)$
there exists $\tau_j\geq0$ such that $y_{j+1}=\Phi_H^{\tau_j}(y_j)$,

 $(ii)$
$y_{j+1}=S_{k,i}(y_j)$,

 $(iii)$
$y_{j+1}=S_{k\to{k+1}}(y_j)$.
\end{itemize}
Here $S_{k,i}$ with $i \in I_k$ stands  for the family of homoclinic maps associated to $\jC_k$ and $S_{k\to{k+1}}$ stands for the family of
heteroclinic maps from $\jC_k$ to $\jC_{k+1}$.

The existence of the sequence $(y_j)$ is proved by induction, starting with $y_0=x_0$.
Assume that $y_{j_n}=x_n$. By definition
$$
x_{n+1}=\ph_k(x_n)\quad\textrm{or}\quad x_{n+1}=\psi_{k,i}(x_n)\quad\textrm{or}\quad x_{n+1}=\psi_{k\to{k+1}}(x_n).
$$

In the first case, set $y_{j_n+1}=x_{n+1}$.

In the second case, by definition of $\psi_{k,i}$, there exists
$\tau\geq0$ such that $S_{k,i}\big(\Phi_H^{\tau}(x_n)\big)=x_{n+1}$. Set $y_{j_n+1}=\Phi_H^{\tau}(x_n)$ and $y_{j_n+2}=x_{n+1}$.

The last case is similar to the latter one when $\psi_{k\to{k+1}}$ is a heteroclinic map: there exists
 $\tau\in\R$ such that $S_{k\to{k+1}}\big(\Phi_H^{\tau}(x_n)\big)=x_{n+1}$, so set
$y_{j_n+1}=\Phi_H^{\tau} (x_n)$ and $y_{j_n+2}=x_{n+1}$. When $\psi_{k\to{k+1}}$ is a transition map (that is, $\jC_k$ and $\jC_{k+1}$
are adjacent cylinders, see Definition~\ref{def:adjacent}),
by definition there exists $\tau\geq0$ such that $x_{k+1}=\Phi_H^{\tau}(x_n)$, so we set $y_{j_n+1}=x_{n+1}$.

The induction stops after a finite number of steps and yields
a sequence  $(y_j)_{0\leq j\leq j_*}$ with $j^*\leq 2n$, which satisfies $(P_1)$ and $(P_2)$.

\paraga This provides us with a sequence $\g:=(g_0,\ldots,g_{j_*-1})$ of elements from the set
\[
\big\{\Phi^{\tau}_H\mid \tau\in[0,+\infty[\big\}\cap \{S_{k,i}\mid 1\leq k\leq k_*,\, i \in I_k\}\cap\{S_{k\to{k+1}}\mid 1\leq k\leq k_*-1\},
\]
such that
\beq
y_{j+1}=g_j(y_j),\qquad 0\leq j\leq j_*-1.
\eeq
For a  point $z$ of $\ha\bC$ that is {\em close enough to $y_0$}, we call the $\g$-orbit
of $z$  the ordered sequence of points
\beq
\big(z,g_0(z),g_1\circ g_0(z)\ldots,g_{j_*-1}\circ\cdots\circ g_0(z)\big).
\eeq
We similarly define the $\g$-orbit of a subset of $\ha\bC$ that is close enough to $y_0$. These
$\g$ orbits are well-defined, thanks to the continuity of the  maps involved and the openness
of their domains.

\paraga
For each $k$, we fix a negatively invariant neighborhood $\jN_k\subset W^-(\hc_k)$
of $\hc_k$ in
$W^-(\hc_k)$ such that, given $z\in\hc_k$ and $w\in\jN_k\cap W^-(z)$, then
\beq\label{eq:dec}
d\big(\Phi_H^{-\tau}(z),\Phi_H^{-\tau}(w)\big)\leq d(z,w),\qquad \forall \tau\geq0.
\eeq
(see Definition~\ref{def:invcyl}, (\ref{eq:decreasing})). In the case of adjacent cylinders
$\hc_k$ and $\hc_{k+1}$, on can moreover assume
\beq\label{eq:compa}
\jN_{k+1}\cap W^-(\hc_{k})\subset \jN_k.
\eeq
To see this,  note that
the union $\hc_k \cup \hc_{k+1}$ is itself a tame cylinder $C$,
and fix a negatively invariant neighborhood $N\subset W^-(C)$ which satisfies (\ref{eq:dec}).
Then set $\jN_k=N\cap W^-(\hc_k)$ and $\jN_{k+1}=N\cap W^-(\hc_{k+1})$.

\paraga We say that a point $z\in\hc_k$ is {\em recurrent} when it is positively and negatively recurrent
under the restriction of
$\Phi_H$ to $\hc_k$.
Recall that  $\hc_k$ admits a Borel measure which is positive on the open sets and invariant under {$\Phi_H$}.
The set of recurrent points has full  measure in $\hc_k$, by the Poincar\'e recurrence theorem.

\paraga We write $B(u,\rho)$ for the ball in $H\inv(\e)$ centered at $u\in H\inv(\e)$ with radius $\rho$. Fix $\de>0$.
Let $\de_0>0$ be so small that the $\g$-orbit $(D_j)_{0\leq j\leq j_*}$ of $D_0:=B(y_0,\de_0)\cap \ha\bC$ is well-defined
(and so $D_{j+1}=g_{j}(D_{j})$ for $0\leq j\leq j_*-1$) and satisfies
\beq\label{eq:balls}
D_j\subset B(y_j,\de/2),\qquad 0\leq j\leq j_*.
\eeq
We will prove the existence of a sequence $(z_j)_{1\leq j\leq j_*}$, where the point
$z_j$ is in $D_j$ and is recurrent,
together with a sequence of  balls $(B_j)_{0\leq j\leq j_*}$, such that, for $0\leq j\leq j_*$:
\vskip2mm
$\bu$   $B_j$ is centered {at  a point in} $W^-(z_j)\cap\jN$ and  $B_j\subset B(z_j,\de/2)$;  \hfill $(C_j)$
\vskip2mm
%\noindent and:
\vskip2mm$\bu$
there exists a time $\sig_j\geq 0$ with $\Phi_H^{\sig_j}(B_j)\subset  B_{j+1} \textrm{ for } 0\leq j\leq j_*-1.$ \hfill $(T_j)$
\vskip2mm
We will construct these sequences backwards, by finite induction.
More precisely, given some recurrent point
$z_{j+1}\in D_{j+1}$ together with a ball $ B_{j+1}$ satisfying $(C_{j+1})$, we will find a recurrent point $z_j\in D_j$, a ball
$ B_j$ satisfying $(C_j)$ and a time $\sig_j\geq0$ which satisfies $(T_j)$.

\paraga Assume first that $g_j=\Phi^{\tau_j}_H$, in which case the sets $D_j$ and $D_{j+1}$
are either contained in the same cylinder, say $\jC_k$, or consecutive adjacent cylinders $\jC_k$, $\jC_{k+1}$.

In the first case, let $ z_j:=\Phi_H^{-\tau_j}(z_{j+1})$, so   $ z_j\in D_j$
and $z_j$ is a recurrent point of $\hc_k$. Let $w_{j+1}$ be the center of $ B_{j+1}$, and set $w_j=\Phi_{H}^{-\tau_j}(w_{j+1})$.
Then $w_j\in W^-(z_j)\cap\jN_k$ by equivariance of the unstable foliation,
and satisfies $d(w_j,z_j)\leq d(w_{j+1},z_{j+1})<\de/2$,
so that there exists a ball $ B_j$ centered at $w_j$ and contained in $B(z_j,\de/2)$, which satisfies
$$
\Phi_H^{\tau_j}( B_j)\subset  B_{j+1}
$$
Hence Conditions $(C_j)$ and  $(T_j)$ with $\sig_j=\tau_j$ are satisfied.

In the case of adjacent cylinders,  $z_j:=\Phi_H^{-\tau_j}(z_{j+1})$ and $z_{j+1}$ both
belong to the intersection $\hc_k\cap\hc_{k+1}$ and $z_j$ is again recurrent in $\hc_k$ by invariance
of $\ha\jC_k$ and the definition of adjacent cylinders.
Let again $w_{j+1}$ be the center of $ B_{j+1}$, so that $w_{j+1}\in \jN_{k+1}\cap W^-(\hc_k)$.
By (\ref{eq:compa}), this proves that $w_j=\Phi_{H}^{-\tau_j}(w_{j+1})$ is in $\jN_k$, since this
latter neighborhood is negatively invariant. As above, $d(w_j,z_j)\leq d(w_{j+1},z_{j+1})<\de/2$
and the conclusion follows.

\paraga Assume now that $g_j=S_{\ell,i}$,
so that
$z_{j+1}\in\hc_\ell$. Let $R_j$ and $R_{j+1}$ be the full-measure subsets of $D_j$ and $D_{j+1}$ formed by
the recurrent points.
Since $S_{\ell,i}$ is $C^1$ and hence preserves the full-measure property:
$$
R_{j+1}\cap S_{\ell,i}\Big(R_j\cap\Domt\, S_{\ell,i}\Big)
$$
is a full measure subset of $D_{j+1}$. Therefore,
there exists a  point
$$
\ov z_j\in R_j\cap\Domt\,S_{\ell,i}
$$
such that $\ov z_{j+1}:=S_{\ell,i}(\ov z_j)$ is recurrent and so close to
$z_{j+1}$ that (by continuity of the unstable foliation)   $W^-(\ov z_{j+1})$ intersects the ball $ B_{j+1}$.

\vskip2mm

By definition of $\Domt\,S_{\ell,i}$, the submanifold $W^-(\ov z_j)$  intersects $W^+(\inte\hc_\ell)$ transversely at some point
$\xi\in W^+(\ov z_{j+1})$. Apply the $\la$-property (Definition~\ref{def:invcyl}) to $W^-(\ov z_j)$ in the neighborhood of $\xi$, together
with the positive recurrence property of $\ov z_{j+1}$:
there exists a (large) time $\tau$ such that $\Phi_H^\tau\big(W^-(\ov z_j)\big)=W^-\big(\Phi_H^\tau(\ov z_j)\big)$ intersects $B_{j+1}$.
Fix
$$
\ze\in W^-\big(\Phi_H^\tau(\ov z_j)\big)\cap B_{j+1}
$$
and note that $\Phi_H^\tau(\ov z_j)$ is  negatively recurrent since $\ov z_j$ is. Hence there exist arbitrarily large times $\tau'$
such that $z_j:=\Phi_H^{-\tau'+\tau}(\ov z_j)\in D_j$.
In particular, one can choose $\tau'$ such that this latter condition is satisfied and
$$
\Phi_H^{-\tau'}(\ze)\in W^-(z_j)\cap B(z_j,\de/2)\cap \jN.
$$
Therefore, there exists a ball $ B_j$ centered at $\Phi_H^{-\tau'}(\ze)$ and contained in $B(z_j,\de/2)$
such that $\Phi_H^{\tau'}(B_j)\subset B_{j+1}$. This proves  $(C_j)$  and  $(T_j)$ with $\sig_j=\tau'$.

\paraga The case where $g_j=S_{\ell-1}^\ell$ is completely similar to the previous one.

\paraga This yields our sequences $(z_j)_{0\leq j\leq j_*}$ and $( B_j)_{0\leq j\leq j_*}$.
To end the proof, observe that for any point
$a_0\subset  B_0$, the point $a_{j}=\Phi_H^{\sig_{j-1}+\cdots+\sig_0}(a_0)$ is in $ B_{j}$, $1\leq j\leq j_*$.
Therefore
$$
\dist(a_{j},y_{j})\leq\dist(a_{j},z_{j})+\dist(z_{j},y_{j})\leq \de,\qquad 0\leq j\leq j_*.
$$
In particular, by $(P_1)$
$$
\dist(a_{j_n},x_n)\leq \de,\qquad 0\leq n\leq n_*.
$$
Therefore the solution $\Phi_H(\cdot,a_0)$ satisfies our requirement, with
$\tau_n=\sig_{j_n-1}+\cdots+\sig_0$ for $1\leq n\leq n_*$.
\end{proof}

\setcounter{paraga}{0}
%%%%%%%%%%%%%%%%%%%%%%%%%%%%%%%%%%%%%%%%%%%%
%%%%%%%%%%%%%%%%%%%Section5%%%%%%%%%%%%%%%%%%%%%
%%%%%%%%%%%%%%%%%%%%%%%%%%%%%%%%%%%%%%%%%%%%

\section{An algorithmic construction of pseudo-orbits}\label{Sec:constructive}
%In this section we introduce stronger compatibility conditions for the polysystem
%attached to a good cylinder, under which we are able to give a different proof of

In this section we provide a  second proof, of   algorithmic nature,  of {\bf Proposition~\ref{prop:tamecyl}}.

Let $\jC$ be a  good cylinder with continuations $\hc$,  let $(\Sig,\ph)$ and $(\ha\Sig,\ha\ph)$ be the corresponding \twist sections, respectively, and $\psi$ be a homoclinic correspondence.
As before, we identify $\Sig\sim\T\times [a,b]$, and $\ha\Sig\sim\T\times[\ha a,\ha b]$.
The  good cylinder condition assumes that for essential invariant circle $\Ga \in \Ess(\ha\ph)$ in $\inte\ha\Sig$
there exists a splitting arc based on $\Ga$.

We will first prove a version of {\bf Proposition~\ref{prop:tamecyl}}, under some compatibility condition between the homoclinic correspondence and the \twist map, formulated below:

\vskip2mm$\bu$ if  $\ha\ph$ tilts the verticals to the right
(resp. left), then  for each element $\Ga\in\Ess(\ha\ph)$ that is the upper boundary of a Birkhoff zone in $\T\times\,]\ha a,\ha b[$,
there exists a right (resp. left) splitting arc
based on $\Ga$. \hfill(VG)

\vskip2mm

The compatibility between the orientation of the splitting arcs and the direction of the \twist of the  map $\ha\ph$  enables us, via the following result, to prove the existence of diffusing orbits without the use of the Poincar\'e recurrence theorem. %opening up the possibility of quantitative estimates.

\begin{prop} \label{prop:tamecyl2}
Let $\jC$ be a  good cylinder.  Assume the setting  of Proposition~\ref{prop:tamecyl}, and the compatibility condition (VG).

Fix $\Ga_\bu\subset \mathbf{A}\setminus\Ga(b)$ in  $\Ess(\ph)$ and
a neighborhood $U_\bu$ of $\Ga_\bu$ in $\bA$.
Then, there exist a point $z_0\in U_\bu$, an integer $N\geq 1$, and, for $n\in\{0,\ldots, N-1\}$,   an index $i_n\in I$
and  an integer $m_n>0$ such that the orbit given by
\[z_{n+1}=\psi_{i_n}\circ\ha\ph^{m_n}(z_n), \quad n=0,\ldots,N-1,\]
has the following properties: \begin{itemize} \item[(i)] $z_{N}\in \Gamma(b)^+$,
\item[(ii)] for every essential
invariant circle  $\Ga$ between $\Ga_\bu$ and $\Ga(b)$, there exists a point $z_n$ that is
$\delta$-close to $\Ga$.
\end{itemize}
\end{prop}

\begin{proof}[Proof of Proposition~\ref{prop:tamecyl2}.]
In the sequel we will assume that $\ha\ph$ tilts the verticals to the right (the case when $\ha\ph$ tilts the
verticals to the left can be dealt with similarly).
For the sake of notational simplicity, we   assume that $\Ga_\bu=\Ga(a)$ and we denote $\bA_*=\T\times[a,\ha b]$.
We still write $\ha\ph$, $\psi$ for the restrictions of the corresponding maps to $\bA_*$.
The proof immediately extends to the case where $\Ga_\bu$ is any element of $\Ess(\ha\ph)$ in the same
way as in the proof of {\bf Proposition~\ref{prop:tamecyl}}.
Following Remark~\ref{rem:restriction}, we will also assume that $\diam\Dom\psi_i<\de$ for $i\in I$.

As in the proof of {\bf Proposition~\ref{prop:tamecyl}}, we want
to prove that given any connected open neighborhood
$U_\bu$
of $\Ga(a)$ in $\bA_*$, the orbit
${U}$
of
$U_\bu$
under the restricted polysystem $(\ha\ph,\psi)_{\vert \bA_*}$ satisfies
\beq
\Ga(b)\subset{\CC({U},\Ga(a))}.
\eeq

The following constructions and lemmas are part of the proof of {\bf Proposition~\ref{prop:tamecyl2}}.

\paraga {\em The Birkhoff procedure.}
%Let us first give a formal definition, inspired by the first step of the proof of Proposition~\ref{prop:tamecyl}.
Consider {\em any} connected positively $\ha\ph$-invariant open subset ${\bf U}\subset \bA_*$ containing $\Ga(a)$ and
such that $\ov {\bf U}\cap \Ga(\ha b)=\emptyset$.
Set $O=\bA_*\setm \ov {\bf U}$.

We define the
{\em filled subset $\bU$} associated with ${\bf U}$ as the open set
\beq\label{eq:orbitU}
\bU=\CC\big(\AA_*\setm \ov {\CC(O,\Ga(\ha b))},\Ga(a)\big).
\eeq
One proves as in {\bf Proposition~\ref{prop:tamecyl}} that
$\Int\cl (\bU)=\bU$, and that
   $\bU$ is homeomorphic to $\T\times[0,1[$  and satisfies
$\ha\ph(\bU)=\bU$. As a consequence,
by the Birkhoff theorem, the frontier $\Fr\bU$ is in $\Ess(\ha\ph)$, and moreover one easily
deduces from the proof of {\bf Proposition~\ref{prop:tamecyl}} that
$$
\Fr\bU\subset\Fr {\bf U}.
$$
We can now introduce our
procedure. Let $\nu$ be a uniform Lipschitz constant for the circles of $\Ess(\ha\ph)$.  We define
$\nu$-balls of $\ha\bA$ as the intersections with $\ha\bA$ of open rectangles satisfying~(\ref{eq:cball}).
Given $\Ga\in\Ess(\ha\ph)$, we call a neighborhood $V$ of $\Ga$ in $\ha\bA$  a $\nu$-neighborhood if $V$ is a union of $\nu$-balls.

Recall that for two elements $\Ga_1,\Ga_2\in\Ess(\ha\ph)$, we write $\Ga_1>\Ga_2$
(resp. $\Ga_1\geq\Ga_2$) if $\gamma_1(\theta)>\gamma_2(\theta)$  (resp. $\gamma_1(\theta)\geq \gamma_2(\theta)$) for all $\theta\in \T$, where $\gamma_i$ is the Lipschitz function whose graph is $\Ga_i$, for $i=1,2$.

%\begin{Def}\label{def:Birkhoffproc}
Denote by $\jP$ the set of pairs $(\Ga,V)$, where  $\Ga\subset \T\times [a,\ha b[$ is in $\Ess(\ha\ph)$, and
 $V$ is  either a $\nu$-neighborhood of $\Ga$ in $\bA_*$,
or a $\nu$-ball centered on $\Ga$.

We define the {\em Birkhoff procedure} as the map
\beq
\B:\jP\longrightarrow \Ess(\ha\ph),\qquad
\B(\Ga,V)=\Fr\bU,
\eeq
where $\bU$ is the filled subset associated with the connected  $\ha\ph$-invariant open set
$$
U=\bigcup_{n\geq0} \ha\ph^n(\Ga^-\cup V)=\Ga^-\cup\big(\bigcup_{n\geq0} \ha\ph^n(V)\big).
$$
%\end{Def}

Recall that $\Ga^-$
%(resp. $\Ga^+$)
stands for the connected component of $\bA_*\setm\Ga$ located below
%(resp. above)
$\Ga$.
Observe that since  $\ha\ph$ is a good \twist map, $\B(\Ga,V)>\Ga$.
Note that the region between $\Ga$ and $\B(\Ga,V)$ is not necessarily a Birkhoff zone.
However, if $\Ga$ is the lower boundary of a Birkhoff zone and if $V\cap \Ga^+$
is contained in the Birkhoff zone, then $\B(\Ga,V)$ is the upper boundary of that Birkhoff zone.
One main property of $\B$ that we will use is the following easy transition lemma.

\begin{lemma}\label{lem:transfert}
Fix $(\Ga,V)\in\jP$ and let $\Ga'=\B(\Ga,V)$.
Then, for every open set $V'$ intersecting~$\Ga'$,
there exists $n\geq 0$ such that $\ha\ph^{n}(V)\cap  V'\cap(\Ga')^- \neq\emptyset$.
%Hence there exists  $z\in V$ such that $\ha\ph^{n}(z)\in V'\cap(\Ga')^-$.
\end{lemma}

\begin{proof}
By definition, $\B(\Ga,V)=\Fr\cU$ so $\B(\Ga,V)\subset \Fr U\subset \cl U$.
Hence $V'\cap U\neq\emptyset$, so that there exists $n\geq0$ with $\ha\ph^n(V)\cap V'\neq\emptyset$.
Moreover $V\subset (\Ga')^-$ by construction, so $\ha\ph^n(V)\subset (\Ga')^-$.
\end{proof}

The following result is a crucial step in our subsequent construction of diffusing orbits. Recall
that we write $\psi=(\psi_i)_{i\in I}$ for the homoclinic correspondence.

\begin{lemma}\label{lem:crux}
Consider a pair $(\Ga,V)\in\jP$ and let $\Ga'=\B(\Ga,V)$. Then there exist $n\geq0$, $i\in I$ and a $\nu$-ball $V'$ centered
on $\Ga'$ such that
$$
V'\subset \psi_i\big(\ha\ph^{n}(V)\big).
$$
%\textcolor{red}{The  integer $m$ can be chosenV_{n_*}
%arbitrarily large.}
\end{lemma}

\begin{proof} By the assumption that $\ha\ph$ is a good \twist map, either $\Ga'$ is the upper boundary of a Birkhoff zone, or  is accumulated from below by a  sequence of essential  circles.

First, assume that  $\Ga'$ is the upper boundary of a Birkhoff zone.
By  condition (VG) there exists a right splitting arc $\ze$  based at a point  $\al_0\in\Ga'$
such that $\til\ze\setminus\{\al_0\}=\ze(]0,1])\subset \Ga^+\cap (\Ga')^-$,
and a homoclinic map $\psi_i$  such that $\til\ze\setm\{\al_0\}\subset\Dom\psi_i$.

\vskip1mm

Let $D$ be a triangular domain associated
with $\til\ze$.
By Lemma~\ref{lem:transfert}, there exist  $z_0\in V$ and
$n\geq0$ such that $z_n=\ha\ph^n(z_0)\in D$. Since $V$ is a $\nu$-ball or a union of $\nu$-balls, by Lemma~\ref{lem:negtilt}, one can
find a negatively tilted arc $\eta$ contained in $V$, emanating from some point of $\Ga$ and ending at $z_0$.
Therefore, by Lemma~\ref{lem:iteratetilt}, $\eta_n=\ha\ph^n\circ \eta$ is a negatively tilted arc
emanating from $\Ga$ and ending at $z_n\in D$.
By Lemma \ref{lem:tiltarc}, $\til\eta_n \cap\big(\til\ze\,\setm\{\al_0\}\big)\neq\emptyset$.

\vskip1mm

As a consequence, $\ha\ph^n(V)\cap (\til\ze\,\setm\{\al_0\})\neq\emptyset$, and therefore there exists a ball
$B$ centered on $\til\ze$ satisfying
$$
B\subset \big(\Dom \psi_i\cap \ha\ph^n(V)\big).
$$
So, by definition of $\til\ze$,
$\psi_i(B)\subset \psi_i\big(\ha\ph^n(V)\big)$ contains a $\nu$-ball $V'$ centered on $\Ga'$, which proves our claim.
%\textcolor{red}{Since in Lemma~\ref{lem:transfert} the integer $n$ can be chosen arbitrarily large, it immediately
%follows that the integer $m$ can be chosen arbitrarily large.}

\vskip1mm
Second, assume that $\Ga'$ is   accumulated from below by a  sequence of essential circles $\Ga_n$, with $\Ga_n<\Ga_{n+1}$ for all $n$.
Let $\ze$ be a splitting arc based at a point in $\Ga'$. There exists $\Ga_{n_*}<\Ga$ such that $\til\ze\subset(\Ga')^-$ intersects both $\Ga_{n_*}$ and $\Ga_{n_*+1}$. There exists a  segment $\ze[\sigma,\sigma']$   of $\til\ze$ such that $\ze(\sigma)\in\Ga_{n_*+1}$, $\ze(\sigma')\in\Ga_{n_*}$, and $\ze]\sigma,\sigma'[$ is contained in the region between $\Ga_{n_*}$ and $\Ga_{n_*+1}$. Lemma \ref{lem:transfert} applied to $\Ga_{n_*+1}^+$, which is a neighborhood of $\Ga'$,  implies that there exists $m\geq 0$ such that $\ha\ph^m(V)$ intersects both  $\Ga_{n_*}$ and $\Ga_{n_*+1}$. Hence, there exists an arc $\til\xi=\xi([a,a'])\subset \ha\ph^m(V)$  that is contained in the region between $\Ga_{n_*}$ and $\Ga_{n_*+1}$, with $\xi(a)\in  \Ga_{n_*+1}$ and $\xi(a')\in\Ga_{n_*}$.
By Lemma \ref{lem:torsion}, for some $m'\geq m$ sufficiently large we have $\ha\ph^{m'-m}(\til\xi)\cap \ze(]\sigma,\sigma'[)\neq \emptyset$. As $\ha\ph^{m'-m}(\til\xi)\subset \ha\ph^{m'}(V)$, it follows that $\ha\ph^{m'}(V)\cap \ze(]\sigma,\sigma'[)\neq \emptyset$.
Since $\ze$ is a splitting arc based at a point in $\Ga'$, there exists a homoclinic map $\psi_i$ such that  $\psi_i(\ha\ph^{m'}(V))$ contains a $\nu$-ball centered on $\Ga'$.
This ends the proof of the lemma.
\end{proof}

\paraga {\em Coherent sequences of circles and the existence of pseudo-orbits.}
%We now  prove (\ref{eq:orbitU}).

Recall that $\Ess(\ha\ph)$ is endowed with an order relation given by the natural order on graphs.

We say that an increasing sequence
$(\Ga_n)_{0\leq n\leq n_*}$
of elements of $\Ess(\ha\ph)$ is {\em coherent} if there exists a sequence
$(V_n)_{0\leq n\leq n_*}$, where each
$V_n$ is either $\nu$-ball centered on $\Ga_n$ or a $\nu$-neighborhood of $\Ga_n$,
%(which is also a neighborhood of a point in $\Ga_n$),
and for every $n=0,\ldots, n_*-1$ we have
\begin{equation}\label{eqn:coherent}
V_{n+1}\subset \psi_{i_n}\big(\ha\ph^{m_n}(V_n)\big),
\end{equation}
for some $i_n\in I$ and $m_n\geq 0$.

We say that  $(V_n)_{0\leq n\leq n_*}$ is  the system of {\em neighborhoods associated with}
$(\Ga_n)_{0\leq n\leq n_*}$ and that $\Ga_{n_*}$ is the {\em upper} circle of the coherent sequence.

\begin{lemma}\label{lem:iterationBirkhoff}
Fix a neighborhood $U_\bu$ of $\Ga(a)$ in $\bA_*$. Choose an essential invariant circle $\Ga_*$ such that $\Ga(a)<\Ga_*\leq\Ga(\ha b)$.

Then there exists a coherent sequence $(\Ga_n)_{0\leq n\leq n_*}$, with associated
neighborhoods $(V_n)_{0\leq n\leq n_*}$, such that
\[\Ga_0=\Ga(a),\qquad V_0\subseteq U_\bu,\qquad \Ga_{n_*}=\Ga_*.\]

Moreover, we can choose $(\Ga_n)_{0\leq n\leq n_*}$ with the additional property that any two consecutive circles $\Ga_n< \Ga_{n+1}$ in the sequence are either $\delta$-close, i.e., $d_H(\Ga_n,\Ga_{n+1})<\delta$, or
 they are at the boundary of a Birkhoff zone.
\end{lemma}

\begin{proof}
We recall that, since $\ha\ph$  is good, for any   $\Ga\leq \Ga_*$, either $\Ga<\Ga_*$ or $\Ga=\Ga_*$.
Let $\Adm$ be the set of all coherent sequences of circles originating at
$\Ga_0$, with $V_0=V$, and such that $\Ga_n\leq \Ga_*$ for all $n=0,\ldots,n_*$.
This set is nonempty by Lemma~\ref{lem:crux}  (while  reduced to the a single coherent sequence
when $\Ga_0$ and $\Ga_*$ bound a Birkhoff zone).
Let $\Adm^\bu$ be the set of upper circles of all elements of $\Adm$. Let
$$
\Ga^\bu=\Sup\Adm^\bu,
$$
(see Lemma~\ref{lem:accum}) so $\Ga^\bu>\Ga(a)$. Assume by contradiction that $\Ga^\bu<\Ga_*$.

Since $\ha\ph$ is good, one proves as in Lemma~\ref{lem:accum} that either $\Ga^\bu\in\Adm^\bu$,
or $\Ga^\bu\notin\Adm^\bu$ and is accumulated
from below in the Hausdorff topology by elements of $\Adm^\bu$.

In the former case, there exists a coherent sequence $(\Ga_n)_{1\leq n\leq n_*}$
with $\Ga_{n_*}=\Ga^\bu$, and with associated neighborhoods  $(V_n)_{0\leq n\leq n_*}$.
Let  $\Ga:=\B(\Ga_{n_*},V_{n_*})$.
By Lemma~\ref{lem:crux} there exist $m\geq 0$, $i\in I$ and a $\nu$-ball $V$
centered on $\Ga:=\B(\Ga_{n_*},V_{n_*})$ such that
\[V\subset \psi_i\big(\ha\ph^{m}(V_{n_*})\big).\]

The neighborhood $V_{n_*}$ can be chosen arbitrarily small, so that we can assume $V_{n_*}\cap \Ga_*=\emptyset$. Hence,
 since $\Ga_{n_*}=\Ga^\bu<\Ga_*$, we also have $\Ga=\B(\Ga_{n_*},V_{n_*})\leq \Ga_*$.
Therefore one can extend the sequence $(\Ga_n)_{0\leq n\leq n_*}$ to a longer coherent sequence, by  adjoining to the previous sequence
$\Ga_{n_*+1}=\Ga$ and $V_{n_*+1}=V$. The resulting  coherent sequence is in $\Adm$, and
has its upper circle located in $(\Ga^\bu)^+$,
which is a contradiction.

Assume now that $\Ga^\bu\notin\Adm^\bu$  and is accumulated
from below by elements of $\Adm^\bu$. In particular, $\Ga^{\bu}$ is not the upper boundary of a Birkhoff zone.
By the good cylinder condition, there exists a  splitting arc $\ze$  based on $\Ga^\bu$,
with $\psi_i\big(\ze(]0,1])\big)\subset\Ga^\bu$ for some $i\in I$.
%Let $[\al,\be]\subset \Ga$, $\til\ze$ and $[\be,\ga]$ be the three parts of the frontier of $T$ as in
%Definition~\ref{def:splittingarc}.
So, by accumulation, there exists a coherent sequence $(\Ga_n)_{1\leq n\leq n_*}$ such that
\[
\ze(1)\in (\Ga_{n_*})^-.
\]
Consider the circle $\Ga=\B(\Ga_{n_*},V_{n_*})$, which satisfies $\Ga_{n_*}<\Ga<\Ga^\bu$
(otherwise $\Ga=\Ga^\bu$, which contradicts our initial assumption $\Ga^\bu\notin\Adm^\bu$).
The set $K=\ze\inv(\Ga)$ is compact and contained in $]0,1[$, set $\sig=\Max K$.
So $\ze(\sig)\in\Ga$ and $\ze(]\sig,1])\subset \Ga^-$ by continuity.
Moreover, there exists $\sig'>\sig$ such that $\ze(\sig')\in\Gamma_{n_*}$; let $\sig'$ be the smallest
one with such property. Thus $\ze([\sig,\sig'])$ is an arc contained in the closed annulus bounded by
$\Ga_{n_*}$ and $\Ga$, and has one endpoint on $\Gamma_{n_*}$ and the other one on $\Gamma$.

By Lemma \ref{lem:torsion2} below, there exists a point  $z\in V_{n_*}$ with the property that for the
vertical segment  $V^-(z)$ below $z$  there exists $m\geq 0$ such that  $\ph^m(V^-(z))$ intersects
$\ze(]\sig,\sig'[)$.
%The power $m$ can be chosen arbitrarily large.
Since $V_{n_*}$ is a $\nu$-ball,
the vertical segment $V^-(z)$ is contained in $V_{n_*}$. It follows that
$\ph^m(V_{n_*})\cap \ze(]\sig,\sig'[)\neq\emptyset$.

%Take some small $\nu$-ball  $V'$ centered at a point in  $\Gamma$ that is disjoint from $V_{n_*}$
%and also from $\ze([\sig,\sig'])$. By Proposition \ref{thm:extensionbirkhoff2}, there exists a point
%$z\in V_{n_*}$ such that $\ph^m(z)\in V'$, where $m$ can be chosen arbitrarily large.
%Denoting by $V^-(z)$ the vertical segment below $z$, which is contained in $V_{n_*}$, it
%follows that $\ph^m(V^-(z))$ intersects $\ze(]\sig,\sig'[)$, provided $m$ is large enough.
%(This follows similarly to Lemma \ref{lem:torsion}). That is, $\ph^m(V_{n_*})\cap \ze(]\sig,\sig'[)\neq\emptyset$.

Therefore,
since $\psi_i(\til\ze)\subset\Ga^\bu$, the image
$
\psi_i\big(\ha\ph^{m}(V_{n_*})\big)
$
contains a $\nu$-ball $V^\bu$ centered on $\Ga^\bu$. We can choose the neighborhood $V^\bu$
small enough so that $V^\bu\cap \Ga_*=\emptyset$, hence $\B(\Ga^\bu,V^\bu)\leq \Ga_*$.
%\textcolor{red}{The integer $m$ can be chosen arbitrarily large.}
Thus one can construct a coherent sequence in $\Adm$ given by
$$
\Ga_0,\ldots,\Ga_{n_*},\Ga^\bu,\B(\Ga^\bu,V^\bu)
$$
whose upper circle is located in $(\Ga^\bu)^+\cap (\Ga_*)^-$,
which is a contradiction. This proves the first claim of the lemma\footnote{in which, as indicated at the beginning of the
proof, one can eventually choose any  $\Ga_0\in \Ess(\ha\ph)$, $\Ga_0<\Ga(\ha b)$, instead of $\Ga(a)$ as first element
of the coherent sequence.}.

Now we show  the second claim of the lemma. We will restrict ourselves to coherent sequences with the property that the neighborhood $V$ of $\Ga(a)$ is contained in a $\delta$-neighborhood of $\Ga(a)$, and that each neighborhood $V_n$ associated to the coherent sequence is contained in a $\delta$-neighborhood of the corresponding $\Gamma_n$. Restricting to coherent sequences satisfying this condition  is possible due to \eqref{eqn:coherent}.

We first construct recursively a finite sequence  of essential circles $\Ga^j$, $j=0,1,\ldots,N$, with the following properties:
\begin{itemize}
\item[(i)]   $\Ga(a)=\Ga^0<\Ga^1<\cdots<\Ga^{N-1}<\Ga^N=\Ga_*$;
\item [(ii)] For each $j\in\{0,\ldots,N-1\}$, either
\begin{itemize}
\item [(ii.a)] the Hausdorff distance between $\Ga^j$ and $\Ga^{j+1}$ is less than $\delta$, or
\item [(ii.b)] the region between $\Ga^j$ and $\Ga^{j+1}$ is a Birkhoff zone.
\end{itemize}
\end{itemize}

We start by setting $\Ga^0=\Ga(a)$. Let $0<\delta'<\delta$. If $\Ga^j$ has been constructed, we define  the set
\begin{equation}\label{eqn:GaSet}\mathscr{T}_{\Ga^j}=\{\Ga\in\Ess(\ph)\,|\, \Ga^j<\Ga<\Ga_*, \textrm{ and } d_H(\Ga^j,\Ga)<\delta'\},\end{equation}
where $d_H$  stands for the Hausdorff distance.

If $\mathscr{T}_{\Ga^j}\neq \emptyset$ we define \begin{equation}\label{eqn:GaSup}\Ga^{j+1}=\textrm{Sup}\mathscr{T}_{\Ga^j}.\end{equation}
By construction,  $d_H(\Ga^j,\Ga^{j+1})\leq \delta'<\delta$.

If $\mathscr{T}_{\Ga^j}=\emptyset$, it follows from  Lemma \ref{lem:accum} that
$\Ga^j$ is the lower boundary of some Birkhoff zone. In this case we define $\Ga^{j+1}$ as the upper boundary
of that Birkhoff zone. In this case we must have $d_H(\Ga^j,\Ga^{j+1})\geq \delta'$.
%Of course, it may happen that
%$d_H(\Ga^j,\Ga^{j+1})<\delta$.

This recursive construction ends after finitely many steps at $\Ga_*$. Otherwise, we would obtain an infinite sequence
$\Ga^j$, $j=0,1,\ldots$, satisfying (ii), with $\Ga^j<\Ga_*$ for all $j$. Take $\Ga^*=\Sup \Ga^j$. Since the sequence $\Ga^j$ is infinite, it accumulates on $\Ga^*$ from below.  Take a $\delta'$-neighborhood of $\Ga^*$ in the Hausdorff topology. There exists
$\Ga^{j_*}<\Ga^*$ that is contained in that neighborhood. It follows that
$$
\Ga^{j_*+1}=\Sup\mathscr{T}_{\Ga^{j_*}}
$$
is also contained in the same neighborhood, with moreover $\Ga^{j_*}<\Ga^{j_*+1}<\Ga^*$.  Thus, both $\Ga^{j_*+1}$ and $\Ga^*$ are
$\delta'$-close to $\Ga^{j_*}$. This is in contradiction with the definition of $\Ga^{j_*+1}$.  Hence the recursive construction ends after finitely many steps, therefore it must end at $\Ga_*$.

At this point we have constructed a  finite sequence  of essential circles $\Ga^j$, $j\in\{0,\ldots,N\}$ satisfying properties (i) and (ii). The first statement of Lemma \ref{lem:iterationBirkhoff}, applied successively to the pairs $\Ga^0$ and $\Ga_*=\Ga^1$, then to
 $\Ga_0=\Ga^1$ and $\Ga_*=\Ga^2$,..., and finally  $\Ga_0=\Ga^{N-1}$ and $\Ga_*=\Ga^N$, yields coherent sequences between each consecutive pair of essential circles. Concatenating these coherent sequences yields a coherent sequence satisfying the second claim of Lemma~\ref{lem:iterationBirkhoff}, which concludes its proof.
\end{proof}

The following  result (which can be viewed as a loose extension of Lemma \ref{lem:torsion}) has been used in the above argument.
\begin{lemma}\label{lem:torsion2} Let $f:\bA\to\bA$ be an area-preserving good \twist map. Let {$\Ga_\bu<\Ga^\bu$} be two nonintersecting
essential invariant circles contained in $\bA$.
Let $\ga:[0,1]\to\bA$ be a  continuous function such that $\ga(0)\in\Ga_\bu$ and $\ga(1)\in\Ga^\bu$.
Let $U_\bu$ be a neighborhood of some point in {$\Ga_\bu$}.
Then there exists a vertical segment $\upsilon$ originating at $\Ga_\bu$ and contained in $U_\bu$, such that  the positive orbit
of $\upsilon$ under $f$ intersects~$\til\ga=\ga([0,1])$.
\end{lemma}
\begin{proof}[Proof of Lemma \ref{lem:torsion2}]
Without loss of generality, we can assume that $U_\bu$ is a $\nu$-ball, that $U_\bu\cap \til\ga=\emptyset$, and that
$\til\ga$ is contained  in the strip bounded by $\Ga_\bu$ and $\Ga^\bu$.
Lifts to $\R\times[a,b]$ of all objects at hand will be denoted by boldface letters. We fix a lift $\bph$ of $\ph$
and a lift $\bga$ of $\ga$. Let $\ha \th_1$ and $\ha \th_2$ in $\R$ be such that the image $\til\bga$ is contained
in the strip $[\ha\th_1,\ha\th_2]\times[a,b]$. Set $\de=1+\ha\th_2-\ha\th_1$.
\vskip1mm
Let $\rho_\bu$ the rotation number of $\Ga_\bu$.  Apply the Birkhoff procedure
 to $\Ga_\bu$ and $U_\bu$, obtaining $\Ga'=\B(\Ga_\bu, U_\bu)\leq \Ga^\bu$.
Since $\ph$ is a good \twist map, $\Ga'\cap \Ga_\bu=\emptyset$, so the rotation number
$\rho'$ of $\Ga'$ satisfies $\rho'>\rho_\bu$.
Therefore,  there exists $k>0$ such that, for every $\bxi\in\bGa_\bu$ and $\bxi'\in\bGa'$ such that
$\pi(\bxi)-\pi(\bxi')\geq-1$:
\begin{equation}\label{eqn:angle_shift}
\pi[\bph^k(\bxi')]-\pi[\bph^k(\bxi)]>\de+1,
\end{equation}
where $\pi$ stands for the first projection.
\vskip1mm
Given $u=(\th,r)\in\bA$, let $B(u,\eps)$ be the $\nu$-ball $]\th-\eps,\th+\eps[\,\times\,]r-\nu\eps,r+\nu\eps[\,\cap\,\bA$.
Choose a point $x'\in\Ga'$, fix $\eps\in\,]0,1/4[$ such that $B(x',\eps)\in\Ga_\bu^+$,  and set
\begin{equation}\label{eqn:neighbor}
V'=\bigcap_{i=0}^{k}\ph^{-i}\Big(B\big(\ph^i(x'),\eps\big)\Big),
\end{equation}
so that $V'\subset\Ga_\bu^+$ is a open neighborhood of $x'$.
\vskip1mm
By Lemma \ref{thm:extensionbirkhoff2}, there exists $z\in U_\bu$ and $n\geq 0$ such that $\ph^n(z)\in  V'$.
Let $x$ be the intersection point between the vertical line through $z$ and $\Ga_\bu$.
Since $U_\bu$ is a $\nu$-ball, the vertical segment $\upsilon=[x,z]$ is contained in $U_\bu$.
The image $\ph^n(\upsilon)$ is a positively tilted arc emanating from $\ph^n(x)\in \Ga_\bu$
and ending at $\ph^n(z)\in V'$. Fix a lift $\bsig=[\bx,\bz]$ of $\upsilon$ and set $\bxi_n=\bph^n(\bx)$  and
$\bz_n=\bph^n(\bz)$. Since $\ph^n(\upsilon)$ is positively tilted,
$$
\pi(\bz_n)\geq\pi(\bxi_n).
$$
By \eqref{eqn:neighbor}, there is a unique lift $\bxi'_n$ of $x'$ such that $\bz\in B(\bxi_n',\eps)$, and moreover
$$
\bph^k(\bz_n)\in B\big(\bph^k(\bxi_n'),\eps\big).
$$
In particular, $\pi(\bxi'_n)-\pi(\bxi_n)\geq-1$ and by \eqref{eqn:angle_shift}:
$$
\pi[\bph^{k}(\bxi'_n)]-\pi[\bph^{k}(\bxi_n)]>\de+1.
$$
Therefore
$$
\pi[\bph^{k}(\bz_n)]-\pi[\bph^{k}(\bxi_n)]>\de.
$$
As a consequence, there is an integer $m$ such that $\bph^k\big(\bph^n(\bsig)\big)$ intersects $m+\til\bga$,
and so $\ph^{k+n}(\upsilon)$ intersects $\til\ga$, which proves our claim.
\end{proof}

\paraga{\em End of proof of Proposition \ref{prop:tamecyl2}.}
By Lemma~\ref{lem:iterationBirkhoff},
applied to $\Ga_*=\Ga(b)$, there exists a
coherent sequence of circles $(\Ga_n)_{0\leq n\leq n_*}$, with associated
neighborhoods $(V_n)_{0\leq n\leq n_*}$, such that
$$
V_0\subseteq U_\bu,\qquad \Ga_{n_*}=\Ga(b).
$$
Then, by continuity, there exists a sequence $(W_n)_{0\leq n\leq n_*}$ of open sets
with $W_{n_*}\subset V_{n_*}$, such, for $0\leq n\leq n_*-1$:
$$
W_n\subset V_n, \qquad \psi_{i_n}\circ \ha\ph^{m_n}(W_n)=W_{n+1}
$$
for suitable $i_n\in I$ and $m_n\geq0$. As a consequence, since $W_{n_*}\cap \Ga(b)^+\neq\emptyset$,
there exists a point $z_0\in W_0$  such that  $z_{n+1}=\psi_{i_n}(\ha\ph^{m_n}(z_n))\in W_{n+1}$, for $n=0,\ldots,N-1$.
In particular $z_N\in \Ga(b)^+$.

At the beginning of the proof we assumed that the maps $\psi_i$ are chosen so that the domain of each map has diameter less than $\delta$. The coherent sequence obtained at the end of Lemma \ref{lem:iterationBirkhoff}
has been constructed so that, for any two consecutive circles $\Ga_n<\Ga_{n+1}$, either $d_H(\Ga_n,\Ga_{n+1})<\delta$, or the region between $\Ga_n$ and $\Ga_{n+1}$ is a Birkhoff zone.  Each essential invariant circle is between a pair $\Ga_n$ and $\Ga_{n+1}$ as above, or is at the boundary of a Birkhoff zone between  $\Ga_n$ and $\Ga_{n+1}$, for some $n$. Hence the pseudo-orbit $(z_n)_{0\leq n\leq N}$ obtained above   gets $\delta$-close  to every essential invariant circle.

This concludes the proof of Proposition~\ref{prop:tamecyl2}.
\end{proof}

The same type of argument is used for the following corollary dedicated to the ``singular cylinders'', for which a constructive
proof can also be deduced from the previous one.

\begin{cor}\label{cor:singcyl2}
We   assume the same condition as in Proposition~\ref{prop:tamecyl2}. Assume that {$K\subset\T\times\,]a+\de,b-\de[$}
is a compact subset of $\bA$, invariant under $\ph$ and contained in the interior of some Birkhoff zone of $\ph$.
%Let $\rho=\dist(K,\d Z)$.
%Assume that for each $i\in I$,
%\beq\label{eq:diamdom1}
%\diam\Dom \psi_i<\rho/2.
%\eeq
%and that $\psi$ is $\de$-bounded.

Then $f=(\ph,\psi)$ admits a pseudo-orbit which does not intersect~$K$. If moreover $\psi$ is $\de$-bounded,
then there exists a $\de$-admissible pseudo-orbit.
\end{cor}

The proof is exactly the same as that of Corollary~\ref{cor:singcyl}.

\begin{proof}[Second proof of Proposition~\ref{prop:tamecyl}]

If the compatibility condition (VG) is satisfied, then the result follows immediately from Proposition~\ref{prop:tamecyl2}.

Assume that the compatibility condition (VG) is not satisfied for some $\Ga\in\Ess(\ha\ph)$, for instance $\ha\ph$ tilts the vertical to the right,  $\Ga$ is the upper boundary of a Birkhoff zone in $\T\times\,]\ha a,\ha b[$,
and there exists a left splitting arc based on $\Ga$.

The key observation is that $\ha\ph^{-1}$ tilts the verticals to the left, and the compatibility condition (VG) is locally satisfied for $\ha\ph^{-1}$ and $\psi$.

Then, as in the proof of Proposition~\ref{prop:tamecyl2} we obtain an orbit $(z_n)$ such that, for each $n$,
either
\[z_{n+1}=\psi_{i_n}\circ\ha\ph^{m_n}(z_n)\]
or
\[z_{n+1}=\psi_{i_n}\circ\ha\ph^{-m_n}(z_n).\]

Using the symmetrization result on polysystems, Lemma~\ref{lem:dense} in Appendix~\ref{app:symmetrization},
yields the result in Proposition~\ref{prop:tamecyl}.

We note here that the proof of Proposition~\ref{prop:tamecyl2} does not require to use
the Poincar\'e recurrence theorem, while this proof Proposition~\ref{prop:tamecyl} requires to use
the Poincar\'e recurrence theorem when the compatibility condition (VG) is not satisfied.
\end{proof}

\paraga{\bf Remarks.} The compatibility condition (VG) between the orientation of the splitting arcs and the direction of the
\twist  may  not be generic, even in specific classes of examples
as \cite{LM}.

We describe an alternative way to eliminate the condition (VG) and obtain an algorithmic proof of {\bf Proposition~\ref{prop:tamecyl}},
without using Poincar\'e recurrence.
Instead we use the equivariance properties of the homoclinic  correspondences (see Appendix~\ref{App:equiv}).

Assume, for instance, that $\ph$ tilts  the verticals to the right, and  that the homoclinic
correspondence $\psi$ is equivariant.

For the key part of the proof of {\bf Proposition~\ref{prop:tamecyl2}}, let $(\Ga,V)\in\jP$ and
$\Ga'=\B(\Ga,V)$. Assume that there is  a \emph{left-splitting arc} $\til\ze$ for $\Ga'$ based at a point
$\al_0\in\Ga'$, such that $\til\ze\setminus \{\al_0\}\subset \textrm{Dom}(\psi_i)$ for some homoclinic
map $\psi_i$, and $\psi_i(\til\ze\setminus\{\al_0\})\subset\Ga$.
Since the  inverse map $\ph^{-1}$ \emph{tilts the verticals to the left},  the argument in the proof of
Lemma \ref{lem:crux}  implies that $\ph^{-n}(V) \cap (\til\ze\setminus\{\al_0\})\neq \emptyset$ for some $n>0$.
Choose $m>n$. Then $\ph^{m-n}(V)\cap \ph^m(\til\ze\setminus\{\al_0\})\neq \emptyset$.
By the definition of an equivariant homoclinic correspondence, there exists a homoclinic map
$\psi_{i(m)}$  satisfying
$$
\Dom \psi_{i(m)}=\ha\ph^m\big(\Dom \psi_i),\qquad \Im\psi_{i(m)}=\ha\ph^m\big(\Im \psi_i),
$$
and
$$
\psi_{i(m)}\circ\ha\ph^m=\ha\ph^m\circ\psi_i.
$$
Therefore
\[\psi_{i(m)}(\ph^m(\til\ze\setminus\{\al_0\}))=\ph^m(\psi_i(\til\ze\setminus\{\al_0\}))\subset \ph^m(\Ga)\subset\Ga,\]
since $\Ga$ is $\ph$-invariant. The rest of the proof goes as before.

%\MG{To produce $\de$-admissible orbits, we only need to  use in the above argument coherent sequences of circles that satisfy the additional property: if for  $\Ga_n< \Ga_{n+1}$ there exists $\Ga\in \Ess (\ha\ph)$ with $\Ga_{n}<\Ga<\Ga_{n+1}$, then  $d_H(\Ga_n,\Ga_{n+1})<\delta$.}
%To produce $\de$-admissible orbits, one moreover needs to assume that the homoclinic correspondence is $\de$-bounded. This is possible when using Method 1, while Method 2 precludes this possibility since equivariant correspondences cannot be bounded by arbitrarily small factors in general.

\section*{Acknowledgements}
Research of M.G. was partially supported by NSF grants  DMS-0635607, DMS-1515851 and  DMS-1814543, and by  the  Alfred P. Sloan Foundation grant G-2016-7320.

%%%%%%%%%%%%%%%%%%%%%%%%%%%%%%%%%%%%%%%%%%%%%%
%%%%%%%%%%%%%%%%%%%%%%%%%%%%%%%%%%%%%%%%%%%%%%
%%%%%%%%%%%%%%%%%%%%%%%%%%%%%%%%%%%%%%%%%%%%%%

%\newpage
\appendix

\section{A reminder on twist and  \twist maps}\label{App:twistmaps}
Here we review some basic concepts on twist and \twist maps, following
\cite{Bo,Ha,HF,Hu,Ma84,Y92,LC,HK}, and prove some auxiliary results that are used in the paper.

Let $a<b$ be fixed. We set
$$
\AA=\T\times[a,b],\qquad
\Ga(a)=\T\times\{a\},\qquad
\Ga(b)=\T\times\{b\}.
$$
The closure of a subset $E\subset \bA$
will  denoted by $\cl E$, and its interior will be denoted by $\Int E$.
The set $\Fr E=\cl E\setm \Int E$ is the frontier of $E$.
A disk is an open, connected and simply connected subset of $\AA$.

An essential circle in $\AA$ is a $C^0$ curve which is homotopic to $\Ga(a)$. We denote by $\Ess(f)$ the set of all
essential circles in $\AA$.

Here we say that $f:\AA\to\AA$ is a {\em twist map} when it is a $C^1$ diffeomorphism, preserves $\Ga(a)$ and $\Ga(b)$
and twists  the verticals  {\em to the right} (resp., {\em to the left}), that is, for $f(\th,r)=(\Th,R)$ we have
\[
\d_r\Th(\th,r)>0 \quad\text{(resp.}\quad \d_r\Th(\th,r)<0\textrm{)},\qquad \forall (\th,r)\in\AA.
\]

We now recall the definition of tilted arcs.

Fix a circle $\Ga\in\Ess(f)$. An {\em arc  emanating from $\Ga$} is a $C^0$ {injective} function $\eta:[0,1]\to\bA$  such that $\eta(0)\in\Ga$
and $\eta(]0,1])\in \Ga^+$. We say that such an arc {\em starts at $\eta(0)$} and {\em ends at $\eta(1)$}.

Given $u,v$ two vectors in $\R^2$, let $\angle(u,v)$ be the  angle from $u$ to $v$ in $[0,2\pi[$, measured counterclockwise. Denote $\partial_r=(0,1)\in\mathbb{R}^2$.

A  $C^1$ arc emanating from $\Ga$ with $\eta'(s)\neq 0 $ for $s\in[0,1]$ is said to be
{\em positively tilted} (resp. {\em negatively tilted}) when \begin{itemize}
\item $\angle\big(\partial_r,\eta'(0)\big)\in\,]0,\pi[$
(resp. $\angle\big(\partial_r,\eta'(0)\big)\in\,]-\pi,0[$),
and \item the continuous lift  $s\in[0,1]\mapsto\angle\big(\partial_r,\eta'(s)\big)\in\R$
is positive (resp. negative) for all $s\in [0,1]$.
\end{itemize}

\begin{figure}[h]
\centering
\includegraphics[width=0.35\textwidth]{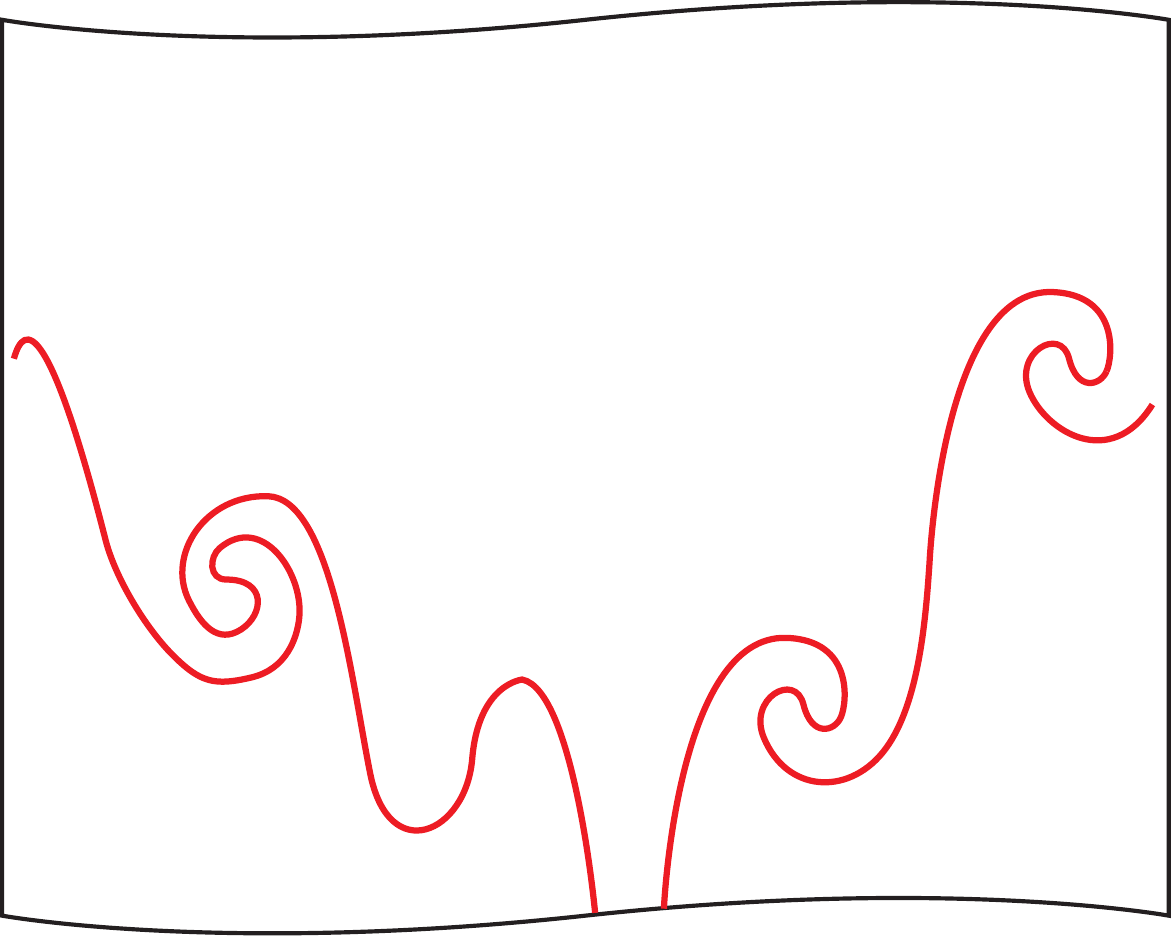}
\caption{Positively and negatively tilted arcs}
\end{figure}

We  now define tilt maps.

The map $f:\AA\to\AA$ is a {\em tilt map} when it is a $C^1$ diffeomorphism, preserves $\Ga(a)$ and $\Ga(b)$
and maps each vertical   into a   {\em positively tilted} (resp., {\em negatively tilted}) arc, that is, for each  map
$\ell:[0,1]\to \AA$ of the form  $\ell(s)=(\theta_0,a+s(b-a))$,  for some $\theta_0\in\mathbb{T}$, we have
\[
\angle(\partial_r, (f\circ \ell)'(s))<0 \quad\text{(resp.}\quad \angle(\partial_r, (f\circ \ell)'(s))>0 \textrm{)},\qquad \forall s\in[0,1].
\]
We also say that the map tilts the verticals to right (resp. left) when it maps verticals onto   {\em negatively tilted} (resp.  {\em positively tilted}) arcs.

A right (resp. left) twist map is a map satisfying $-\pi < \angle(\partial_r, (f\circ \ell)'(s)) < 0$ (resp. $0< \angle(\partial_r, (f\circ \ell)'(s)) < \pi$)
for every vertical $\ell$.
That is, any twist map is a tilt map (but not every tilt map is a twist map).

Any composition of twist maps of the same type is a tilt map (but not necessarily a twist map).
Any composition of tilt maps of the same type is a tilt map of the same type.

A continuous map $f:\AA\to\AA$ is said to be {\em area-preserving} when it leaves invariant a Radon measure
which is positive on the open subsets of $\AA$.

In the sequel we will always assume that $f$ is a twist map or a tilt map  that tilts the verticals to the right.

Below we list several results from  the Birkhoff theory, which are valid for both twist and tilt maps.

\vskip2mm
\noindent
{\bf Theorem (Birkhoff).} {\em Let $f:\AA\to\AA$ be an area-preserving twist map (resp. tilt map).

Then there exists
$\nu>0$ such that every essential circle invariant under $f$ is the graph of some
$\nu$--Lipschitz function $\gamma:\T\to[a,b]$.}

\vskip2mm

The second result from Birkhoff's theory we need is the following.

\vskip2mm

\noindent
{\bf Theorem (Birkhoff).} {\em Let $f:\AA\to\AA$ be an area-preserving twist map (resp. tilt map). Assume that
$U$ is an open subset of $\AA$ homeomorphic to $\T\times [0,1[$,  with $\Ga(a)\subset U$,
such that $f(U)\subset U$ and such that $U$ is the interior of its closure. Then the frontier
$\Fr U$ is an invariant essential circle.}

\vskip2mm

One easily deduces from the first Birkhoff theorem that the set $\Ess(f)$ of essential invariant circles  of $f$,
endowed with the Hausdorff
topology, is compact. Given $\Ga\in\Ess(f)$ with $\Ga=\Graph(\gamma)$, we set
\beq\label{eq:gapm}
\Ga^+=\big\{(\th,r)\in\bA\mid r>\gamma(\th)\big\},\quad
\Ga^-=\big\{(\th,r)\in\bA\mid r<\gamma(\th)\big\}.
\eeq
By the Poincar\'e theory, every $\Ga\in\Ess(f)$  admits a rotation number in $\T$ for $f_{\vert \Ga}$.
One can choose a common lift to $\R$ for the rotation number of all circles, which yields a function $\rho:\Ess(f)\to\R$.
This function is continuous and increasing, in the sense that if $\Ga_i=\Graph\gamma_i$, $i=1,2$ are invariant with
$\gamma_1\leq\gamma_2$, then $\rho(\gamma_1)\leq\rho(\gamma_2)$. Moreover,  $\rho(\gamma_1)<\rho(\gamma_2)$ when $\gamma_1<\gamma_2$.

\begin{Def}\label{def:birkzone} Let $f:\bA\to\bA$ be an area-preserving twist map (resp. tilt map) of the annulus $\bA$.
Let $\gamma_\bu,\gamma^\bu:\T\to\,]a,b[$ be two functions with $\gamma_\bu<\gamma^\bu$,
whose graphs $\Ga_\bu$ and $\Ga^\bu$
are in $\Ess(f)$.

The set
$$
\jB=\{(\theta,r)\,|\, \theta\in \mathbb{T},\ \gamma_\bu(\theta)\leq r\leq \gamma^\bu(\theta)\}
$$
is said to be a {\em Birkhoff zone}  when that there is no element
$\Gamma=\Gr\gamma \in\Ess(f)$ such that $\gamma_\bu\leq \gamma\leq \gamma^\bu$ and $\gamma\neq\gamma_\bu$, $\gamma\neq\gamma^\bu$.
\end{Def}

We  now   prove  Lemma~\ref{lem:accum}, for which we refer to Definition~\ref{def:good}.

\begin{proof}[Proof of Lemma~\ref{lem:accum}]
The main property of a good twist (resp. tilt) map  $f$, coming from the fact that no
element of $\Ess(f)$ has rational rotation number, is that two distinct elements of $\Ess(f)$ are disjoint (see \cite{HK}, Section 13.2).
As a consequence, the rotation number $\rho:\Ess(f)\to\R$ is a homeomorphism onto its image $\cR=\rho\big(\Ess(f)\big)$,
by the compactness of $\Ess(f)$, {where $\Ess(f)$ is endowed with the uniform topology}.
The boundaries of the Birkhoff zones are mapped by $\rho$ onto the boundaries of the maximal intervals in the complement
$\Rot\setm\rho\big(\Ess(f)\big)$, where $\Rot=\big[\rho(\Ga(a)),\rho(\Ga(b))\big]$ is the rotation interval of $f$.
The other claims easily follow.
\end{proof}

The following easy result on negatively tilted arcs is used several times in our constructions.

\begin{lemma}\label{lem:negtilt}
Let $\Ga$ be an essential circle of $\bA$ which is the graph of a $\nu$-Lipschitz function $\gamma:\T\to[0,1]$,
and let $B$ be a $\nu$-ball centered on $\Ga$.
Then for every $z\in\Ga^+\cap B$, there exists a negatively tilted arc emanating from $\Ga$ and ending at $z$,
whose image is contained in $B$.
\end{lemma}

\begin{proof}
Set $B=B_\th\times B_r$. Then since $B$ is a $\nu$-ball the graph of $\gamma_{\vert B_\th}$ is contained in $B$.
Given $z=(\th,r)\in \Ga^+\cap B$ and $\th_0\in B_\th$ with $\th_0<\th$ close enough to $\th$ so that
$$
\frac{r-\gamma(\th_0)}{\th-\th_0}>\nu,
$$
and setting
$z_0=\big(\th_0,\gamma(\th_0)\big)$, the segment
$[z_0,z]$ satisfies our requirements.
\end{proof}

The proof of the following lemma is immediate.

\begin{lemma}\label{lem:torsion} Let $f:\bA\to\bA$ be an area-reserving twist map (resp. tilt map). Let $\Ga^\pm$ be two
nonintersecting essential invariant circles contained in $\bA$.

Then for any continuous curves $C$ and $C'$ which intersect both circles $\Ga^\pm$, the positive orbit
of $C$ under $f$ intersects~$C'$.
\end{lemma}

We refer to \cite{LC86} and \cite{LC87} for  the proofs of the following two results from Birkhoff's theory.

\begin{lemma}\label{lem:iteratetilt} Let $f:\bA\to\bA$ be an area-preserving twist map (resp. tilt map), and let $\Ga$ be an essential invariant
circle for $f$.

Then, the inverse image $f\inv\circ\eta$ of a positively tilted arc $\eta$ emanating from $\Ga$ is a positively tilted arc emanating from $\Ga$.
Also, the direct image $f\circ\eta$ of a negatively tilted arc $\eta$ emanating  from $\Ga$ is a negatively tilted arc emanating from $\Ga$.
\end{lemma}

Given a point $x\in\bA$, we define the lower vertical $V^-(x)$ as the vertical segment joining  $x$ to the corresponding point of $\dma$.

\begin{lemma}\label{lem:posneg} Let $f:\bA\to\bA$ be an area-preserving twist map (resp. tilt map). Let $\Ga$ be an essential invariant
circle for $f$.
Let $X$ be a connected closed subset of $\bA$ which disconnects the annulus $\bA$ and such that $X\subset \Ga^+$.
Let $x\in\bA$ be such that there exists a positively tilted arc $\ga$  and a negatively tilted arc $\eta$, both emanating from
$\Ga$ and ending at $x$, such that the images of $\ga$ and $\eta$ do not intersect $X$.

Then the vertical $V^-(x)$
does not intersect $X$.
\end{lemma}

We can now prove our second lemma on good \twist maps and triangular domains associated with right or left
splitting arcs, stated in Section~\ref{Sec:setting}.

\begin{proof}[Proof of Lemma~\ref{lem:tiltarc}]
We provide the proof only for the case when $\ze$ is a right splitting arc and $\eta$ is a negatively tilted arc; the other case can be proved similarly.
Let $D=D_{\ze\mid[0,s_{*}]}$ as in Definition \ref{lem:tiltarc}. With the notation from the definition $D$  is bounded by the arcs $\ze([0,s_{*}])$, $[\be_*,\al_*]$, and $[\al_*,\al_0]$.
We claim that in fact $\eta(]0,1[)\cap \ze([0,s_{*}])\neq\emptyset$.  Assume by contradiction that the arc $\eta$ is such that $\eta(1)\in D$ but $\tilde\eta\cap  \ze([0,s_{*}])=\emptyset$. Note that $\ze([0,s_{*}])$ is the image of a Lipschitz arc (a piece of an  essential circle)  by a $C^1$ local diffeomorphism $\psi^{-1}_i$. There exists $0<t<1$ such that $\eta(t)\in ]\be_*,\al_*[$;  let $t_x$ be the {largest} $t$ with this property, and let $x=\eta(t_x)$.
Note that if $t'>t_x$ then $\eta(t')\in D$.  Let $x'=\eta(t_{x'})=(\theta',r')$, with $t'_x>t_x$ sufficiently  close to $t_x$  so that the vertical segment  from $x'=(\theta',r')$ to $\al'=(\theta',\gamma(\theta'))\in \Gamma$ is contained in $D$. Note that the vertical arc $V^{-}(x')$ below $x'$ intersects $\ze([0,s_{*}])$, since this curve bounds $D$ from below.

The arc $\eta_{\mid[0,t_{x'}]}$ is a negatively tilted arc from $\eta(0)\in \Ga_\bu$ to $x'\in D$, by assumption.

There exists  a positively tilted curve from $\Ga_\bu$ to $x'$ that does not intersect $X$, as we will shown now.
Choose $\th''>\textrm{max}_{s\in[0,s_{*}]}\pi(\zeta(s))$, where $\pi$ denotes the projection onto the $\th$-coordinate.
Let $y'\in\Ga$ be the intersection point of the vertical through $x'$ with $\Ga$,
and   $y''\in\Ga$ be given by $y''=(\theta'',\gamma(\theta''))$.
Let $[y'',y']_\Ga$ be the segment of $\Ga$ between   $y'$ and  $y''$, traversed from $y''$ to $y'$.
Choose  a point $z\in \Ga_\bu$ and a positively tilted arc from $z$ to $y''$, which we denote by $[z,y'']$.
We define a $C^0$ arc emanating from $z\in \Ga_\bu$ and ending at $x'$ as the concatenation of the
arcs $[z,y'']$, $[y'',y']_\Ga$ and $[y',x']$. Approximating the segment $[y'',y']_\Ga$ from below
and ``rounding the corners''  yields a positively tilted $C^1$ arc  $\xi$ from $z$ to $x'$.

Let $X=\Gamma \cup \ze([0,s_{*}])$. Observe that $\eta[0,t_{x'}]$ is a negatively tilted curve from $\eta(0)\in \Ga_\bu$ to $x'$, which does not intersect $X$ by assumption, and that $\xi$ is a positively tilted arc from $z\in \Ga_\bu$ to $x'$ that does not intersect $X$.  Therefore, Lemma \ref{lem:posneg} implies $V^{-}(x')$ does not intersect $X$. This is a contradiction, since $V^{-}(x')$ must intersect $\ze([0,s_{*}])$ as previously noted.
\end{proof}

The following strong connecting lemma appeared with a different proof in \cite{GR13}.

\begin{prop}\label{thm:extensionbirkhoff2} Let $f:\bA\to\bA$ be a area-preserving twist map (resp. tilt map)\footnote{Here the map is not necessarily assumed to be good}.
Let $\Ga_\bu$ and $\Ga^\bu$ be the boundary components of some Birkhoff zone of instability
for $f$.  Fix a  pair of open sets $V_\bu,V^\bu$ which intersect $\Ga_\bu$ and $\Ga^\bu$ respectively,
with moreover $V_\bu\subset (\Ga^\bu)^-$.
Then there exist a point $z\in V_\bu$ and an integer $n\geq0$ such that $f^{n}(z)\in V^\bu$.
Moreover the integer $n$ can be chosen arbitrarily large.
\end{prop}

\begin{proof}
Set
$$
U=\bigcup_{n\geq 0} f^n(\Ga_\bu^-\cup V_\bu)=\Ga_\bu^-\cup\big(\bigcup_{n\geq 0} f^n(V_\bu)\big),
$$
so that $U$ is a connected and $f$-invariant neighborhood of $\dma$, which satisfies
$$
U\subset (\Ga^\bu)^-.
$$
Hence the frontier $\Ga:=\Fr\cU$ of its associated filled subset (see Section~\ref{Sec:constructive})
is in $\Ess(f)$ and satisfies $\Ga_\bu\leq\Ga\leq\Ga^\bu$. Therefore
$\Ga=\Ga_\bu$ or $\Ga=\Ga^\bu$. The former equality is impossible by construction, so $\Ga=\Ga^\bu$.

\vskip2mm

As a consequence, $\Ga^\bu\subset\Fr \cU\subset \Fr U$, so there exists an integer $n\geq 0$ such that
$$
f^{n}(V_\bu)\cap V^\bu\neq\emptyset,
$$
which proves our claim.
Finally, observe that by choosing arbitrarily small open subsets $W_\bu\subset V_\bu$, $W^\bu\subset V^\bu$
and applying the previous result to the pair $W_\bu$, $W^\bu$, one can ensure that the  integer
$n$ can be chosen arbitrarily large.
\end{proof}

%%%%%%%%%%%%%%%%%%%%%%%%%%%%%%%%%%%%%%%%%%%%%%
%%%%%%%%%%%%%%%%%AppendixB%%%%%%%%%%%%%%%%%%%%%%%%
%%%%%%%%%%%%%%%%%%%%%%%%%%%%%%%%%%%%%%%%%%%%%%

\section{Equivariance properties for homoclinic correspondences}\label{App:equiv}
We describe here some equivariance properties of the homoclinic correspondences with
respect to the Hamiltonian flow, see \cite{DLS08} for related equivariance results.
These properties reveal themselves to be useful in the constructive framework of
Section~\ref{Sec:constructive}.

\begin{Def}
Consider a $C^2$ Hamiltonian function $H$ on $\A^3$,  an energy
$\e$ and a tame cylinder $\jC$ at energy $\e$, with continuation $\hc$ and
continued \twist section~$\ha\Sig$ with return map $\ha\ph$. We say that a homoclinic correspondence $\psi=(\psi_i)_{i\in I}$
is {\em equivariant} when for each $i\in I$ and for each $m\in\Z$ there exists an index $i(m)$
such that
\beq
\Dom \psi_{i(m)}=\ha\ph^m\big(\Dom \psi_i),\quad \Im\psi_{i(m)}=\ha\ph^m\big(\Im \psi_i)
\eeq
and
\beq
\psi_{i(m)}\circ\ha\ph^m=\ha\ph^m\circ\psi_i.
\eeq
\end{Def}

The main result of this section is the following.

\begin{lemma} With the same assumptions as in the previous definition,
any homoclinic correspondence $(\psi_i)_{i\in I}$ can be extended to an equivariant homoclinic
correspondence, in the sense that there exists an equivariant homoclinic correspondence
$(\bpsi_i)_{i\in \bI}$ with $\bI\supset I$ and $\bpsi_i=\psi_i$ for $i\in I$.
\end{lemma}

\setcounter{paraga}{0}
\begin{proof}
Given a subset $A\subset \hc$ and a function $T:A\to\R$, we denote by $\bPhi^T:A\to \hc$
the map defined by
$$
\bPhi^T(x)=\Phi_H\big(T(x),x\big),\quad x\in A.
$$

\paraga Fix an element $\psi=\psi_i:\Dom\psi\to\Im\psi$ of the initial homoclinic correspondence.
By definition, there exist a $C^1$ map $S:\Dom S\to\Im S$ and a $C^1$ function
$\tau:\Dom\psi\to\R^+$ such that the following diagram
$$
\begin{CD}
\Dom\psi@>\psi>>\Im\psi\\
@VV\bPhi^{\tau}V@VVjV\\
\Dom S@>S>>\Im S\\
\end{CD}
$$
commutes (where we denoted by $j$ the canonical inclusion). The map $S$ moreover satisfies
condition (\ref{eq:redhom3}) of Definition~\ref{def:scattering}.

\paraga Since $\Dom S$ is an open subset of $\inte \hc$ which contains $\bPhi^\tau(\Dom\psi)$, one can find
$C^1$ functions $\al,\be:\Dom \psi\to\R$ such that $\al <\tau<\be$ and
$$
D(S):=\bigcup_{x\in\Dom\psi} \Phi_H\big(]\al(x),\be(x)[\,\times\{x\}\big)\subset \Dom S.
$$
By equivariance of $S$ and commutativity of the previous diagram
$$
I(S):=\bigcup_{x\in\Im\psi} \Phi_H\big(]\al(x)-\tau(x),\be(x)-\tau(x)[\,\times\{x\}\big)
$$
is an open subset of $\Im S$.
We still denote by $S$ the induced map from $D(S)$ to $I(S)$. Note that the previous diagram
still commutes if $S:\Dom S\to\Im S$ is replaced with $S:D(S)\to I(S)$. Note moreover that if
$D_t(S):=D(S)\cap \Homt \hc$, then $D_t$ is an open subset of full measure of $D(S)$ such that
\beq\label{eq:dts}
\forall x\in D_t(S),\quad W^-(x)\cap W^+\big(S(x)\big)\cap \Homt\hc\neq\emptyset.
\eeq

\paraga For $m\in\Z$, we denote by $T^{(m)}:\ha\Sig\to\R^+$ the $m^{th}$ return time map associated
with the Hamiltonian flow,
so that $\ha\ph^m=\bPhi^{T^{(m)}}:\ha\Sig\to\ha\Sig$. We still denote by $T^{(m)}$ the ``fiberwise continuation''
of $T^{(m)}$ to the domains $D(S)$ and $I(S)$, that is:
$$
\begin{array}{lll}
&\forall x\in\Dom\psi,\quad \forall z\in \Phi_H\big(]\al(x),\be(x)[\,\times\{x\}\big),\qquad T^{(m)}(z)=T^{(m)}(x),\\[5pt]
&\forall y\in\Im\psi,\quad  \forall z\in\Phi_H\big(]\al(y)-\tau(y),\be(y)-\tau(y)[\,\times\{y\}\big),\qquad T^{(m)}(z)=T^{(m)}(y).
\end{array}
$$
Then, observe that the map $S^{(m)}:\Phi^{T^{(m)}}\big(D(S)\big)\to \Phi^{T^{(m)}}\big(I(S)\big)$ defined by
$$
S^{(m)}=\Phi^{T^{(m)}}\circ S\circ \Phi^{-T^{(m)}}
$$
satisfies
$$
\forall z\in \Phi^{T^{(m)}}\big(D_t(S)\big), \quad W^-(z)\cap W^+(S^{(m)}(z))\cap\Homt\hc\neq\emptyset,
$$
by equivariance of the characteristic foliations and preservation of the transversality. Moreover, $\Phi^{T^{(m)}}\big(D_t(S)\big)$
is an open subset with full measure of $\Phi^{T^{(m)}}\big(D(S)\big)$.

\paraga For $m\in\Z$, let $\psi^{(m)}$ be the unique map such that the following diagram
$$
\begin{CD}
\ha\psi^{(m)}(\Dom\psi)@>\psi^{(m)}>>\ha\psi^{(m)}(\Im\psi)\\
@VV\bPhi^{\tau}V@VVjV\\
\bPhi^{T^{(m)}}\big(D(S)\big)@>S^{(m)}>>\bPhi^{T^{(m)}}\big(I(S)\big)\\
\end{CD}
$$
commutes. Then clearly $\Dom\psi^{(m)}=\ha\ph^m(\Dom \psi)$, $\Im\psi^{(m)}=\ha\ph^m(\Im\psi)$
and
\beq\label{eq:equiv}
\forall x\in \Dom\psi,\qquad  \psi^{(m)}(x)\circ \ha\ph^m =\ha\ph^m\circ\psi.
\eeq
\paraga Let us consider the new index set $\bI=I\times\Z$ and set
$$
\psi_{(i,m)}=\psi_i^{(m)}.
$$
Then the previous construction proves that $(\psi_{(i,m)})_{(i,m)\in\bI}$ is a homoclinic correspondence, with associated
family $(S_i^{(m)})_{(i,m)\in \bI}$.

\paraga It suffices now to show that it is equivariant. Fix ${\bf i}=(i,k)\in\bI$, set $\bpsi:=\psi_{\bf i}=\psi_i^{(k)}$,  fix $m\in\Z$ and
set $\bpsi^{(m)}=\psi_i^{(m+k)}$.
Then
$$
\begin{array}{lll}
&\Dom \bpsi^{(m)}=\ha\ph^{m+k}(\Dom \psi_i)=\ha\ph^m(\Dom\bpsi),\\[5pt]
&\Im \bpsi^{(m)}=\ha\ph^{m+k}(\Im \psi_i)=\ha\ph^m(\Im\bpsi),
\end{array}
$$
and, by (\ref{eq:equiv}):
$$
\bpsi^{(m)}\circ\ha\ph^{m+k}=\ha\ph^{m+k}\circ \psi_i=\ha\ph^m\circ\ha\ph^k\circ\psi_i=\ha\ph^m\circ\bpsi\circ\ha\ph^k
$$
so that
$$
\bpsi^{(m)}\circ\ha\ph^{m}=\ha\ph^m\circ\bpsi.
$$
This proves the equivariance condition for $\psi_{\bf i}$, with ${\bf i}(m)=(i,m+k)$.
\end{proof}

Observe that while homoclinic correspondences giving rise to good cylinders can always be modified
in order to obtain $\de$-bounded correspondences, this is no longer the case for equivariant correspondences,
due to the extension process.

%%%%%%%%%%%%%%%%%%%%%%%%%%%%%%%%%%%%%%%%%%%%%%
%%%%%%%%%%%%%%%%%AppendixC%%%%%%%%%%%%%%%%%%%%%%%
%%%%%%%%%%%%%%%%%%%%%%%%%%%%%%%%%%%%%%%%%%%%%%

\section{Symmetrization of polysystems}\label{app:symmetrization}

The following general lemma has been used several times is specific contexts.

\begin{lemma}\label{lem:dense}
Let $A$ be a metric space endowed with a finite Borel measure, positive on the nonempty open subsets of $A$.
Let $\ph$ be a measure-preserving homeomorphism of $A$  and let $(\psi_i)_{i\in I}$
be a polysystem on $A$, where $\Dom\psi_i$ is open and the map $\psi_i:\Dom \psi_i\to\Im \psi_i$ is a homeomorphism,
for all $i\in I$. Fix a nonempty open subset $V\subset A$. Let $U_f$ and $U_g$ be the full orbit of
 $V$ under the polysystems $f=\big(\ph,\psi=(\psi_i)_{i\in I}\big)$ and $g=\big(\ph,\ph\inv,\psi=(\psi_i)_{i\in I}\big)$
respectively. Then $U_f$ is contained and dense in $U_g$.
\end{lemma}

\begin{proof}
Since $\ph$, $\ph\inv$ and $\psi_i$ are open maps,
the full orbits $U_f$ and $U_g$ are open in $A$, and clearly $U_f\subset U_g$.
We assume that the index set $I$ does not contain $\{-1,1\}$ and we write $\bI=\{-1,1\}\cup I$.
Set $\tau_{-1}=\ph\inv$, $\tau_1=\ph$, $\tau_i=\psi_i$ for $i\in I$, so that
$$
f=(\tau_1, (\tau_i)_{i\in I}),\qquad
g=(\tau_{-1},\tau_1, (\tau_i)_{i\in I}).
$$
Fix a nonempty open subset $\bW \subset U_g$.
Then (by continuity of the maps) there exist a nonempty open subset $W$ of $V$ and a sequence
$\om=(\om_n)_{0\leq n\leq n_*-1}\in\bI^{n_*}$ such that $g^\om(W)\subset \bW$.
We will iteratively modify the set $W$ and the
sequence $\om$ so that we get a nonempty open set
$\til W\subset W$ and a sequence $\til\om\in\big(1\cup I\big)^{m_*}$ such that
$$
f^{\til\om}(\til W)=g^{\om}(\til W)\subset \bW.
$$
At each step of the process, the number of occurrences of the index $-1$ in the sequence $\om$ will decrease
by $1$, and so the process stops after a finite number of steps. Let us describe the first one.

We set, for $0\leq n\leq n_*-1$:
$$
W_n=g^{(\om_0,\ldots,\om_{n-1})}(W)
$$
so that $W_{n_*}=g^{\om}(W)\subset\bW$.
Let $\ov n\geq 0$ be the smallest index such that $\om_{\ov n}=-1$. Since $\tau_1=\ph$ is measure-preserving, there
exists a nonempty open subset $O\subset W_{\ov n}$ and an integer $\nu\geq 1$ such that
$$
\tau_1^{\nu+1}(O)\subset W_{\ov n},
$$
by the Poincar\'e recurrence theorem.  Consider  the new sequence {$\om'=(\om'_0,\ldots,\om_{n'_*+\nu-2})$} with
%$$
%\begin{array}{ll}
%\om'_n=\om_n,\quad &0\leq n\leq \ov n-1,\\
%\om'_n=1,\quad &\ov n\leq n\leq\ov n+\nu,\\
%\om'_n=\om_{n-\nu},\quad &\ov n+\nu+1\leq n\leq n_*+\nu-1,
%\end{array}
%$$
{$$
\begin{array}{ll}
\om'_n=\om_n,\quad &0\leq n\leq \ov n-1,\\
\om'_n=1,\quad &\ov n\leq n\leq\ov n+\nu-1,\\
\om'_n=\om_{n-\nu+1},\quad &\ov n+\nu\leq n\leq n_*+\nu-2,
\end{array}
$$
}
where the first line has to be omitted when $\ov n=0$.
Set
$$
W'=\tau_{\om_0}\inv\circ\cdots\circ\tau_{\om_{\ov n}-1}\inv(O)\subset W,
$$
and note that $W'\neq\emptyset$.
Set
$$
W'_0=W',\qquad W'_n={g}^{(\om'_0,\ldots,\om'_{n-1})}(W'),\quad 1\leq n\leq n_*+\nu-1.
$$
Therefore, by construction:
$$
W'_n\subset W_n,\quad 0\leq n\leq \ov n-1, \quad\text{and}\quad  W'_{\ov n}=O\subset W_{\ov n}.
$$
Hence
$$
W'_{\ov n+\nu}=\tau_1^{\nu}(W'_{\ov n})=
\tau_1\inv\big(\tau_1^{\nu+1}(W'_{\ov n})\big)=
\tau_1\inv\big(\tau_1^{\nu+1}(O)\big)\subset
\tau_{-1}(W_{\ov n})=
W_{\ov n+1}.
$$
As a consequence, by definition of the sequence $\om'$:
$$
W'_{n}\subset W_{n-\nu+1},\quad \ov n+\nu+1\leq n\leq n_*+\nu-1.
$$
In particular
$$
W'_{n_*+\nu-1}\subset W_{n_*}
$$
and if $\ell$ is the number of occurrences of $-1$ in $\om$,
then $-1$ occurs $\ell-1$ times in $\om'$. Iterating this process therefore yields a nonempty open set
$\til W\subset W$ and a sequence $\til\om\in\{0,1\}^{m_*}$ such that
$
f^{\til\om}(\til W)\subset W_{n_*}\subset \bW
$,
so that the full orbit $U_f$ of $V$ under $f$ intersects $\bW$. Since $\bW$ is arbitrary in $U_g$,
this proves that $U_f$ is dense in $U_g$.
\end{proof}

%%%%%%%%%%%%%%%%%%%%%%%%%%%%%%%%%%%%%%%%%%%%%%
%%%%%%%%%%%%%%%%%Bibliographie%%%%%%%%%%%%%%%%%%%%%%%
%%%%%%%%%%%%%%%%%%%%%%%%%%%%%%%%%%%%%%%%%%%%%%

\end{document}